\input ./style/arxiv-general.cfg
\documentclass[aop,MSNbibl,seceqn,dvips]{arximspdf}
\makeatletter
   \@ifpackageloaded{graphicx}{}{\usepackage{graphicx}}
\makeatother
\usepackage{mathbh}

\doi{10.1214/15-AOP1023}
\volume{44}
\issue{3}
\pubyear{2016}
\firstpage{2349}
\lastpage{2425}
\docsubty{FLA}

\makeatletter
\def\sfrac#1#2{#1/#2}

\def\afrac#1#2{#1/(#2)}

\def\sklfrac#1#2{(#1/#2)}
\def\sklvfrac#1#2{((#1)/#2)}

\newcommand{\llbracket}{[\![}
\newcommand{\rrbracket}{]\!]}
\newcommand{\dllbracket}{\bigl[\!\bigl[}
\newcommand{\drrbracket}{\bigr]\!\bigr]}
\newcommand{\rrvert}{\vert}
\newcommand{\rrVert}{\Vert}
\newcommand{\llvert}{\vert}
\newcommand{\llVert}{\Vert}
\renewcommand{\mid}{|}
\newtheorem{theorem}{Theorem}[section]
\newtheorem{lemma}[theorem]{Lemma}
\newtheorem{corollary}[theorem]{Corollary}
\newtheorem{proposition}[theorem]{Proposition}
\newproclaim{definition}[theorem]{Definition}
\newproclaim{assumption}[theorem]{Assumption}
\newproclaim{remark}[theorem]{Remark}
\newcommand{\R}{\mathbb{R}}
\newcommand{\C}{\mathbb{C}}
\newcommand{\ls}{\mathcal{L}}
\newcommand{\N}{\mathbb{N}}
\newcommand{\Z}{\mathbb{Z}}
\newcommand{\E}{\mathbb{E}}
\newcommand{\diag}{\operatorname{diag}}
\newcommand{\Tr}{\operatorname{Tr}}
\newcommand{\supp}{\operatorname{supp}}
\newcommand{\re}{\operatorname{Re}}
\newcommand{\im}{\operatorname{Im}}
\newcommand{\caA}{{\mathcal A}}
\newcommand{\caB}{{\mathcal B}}
\newcommand{\caC}{{\mathcal C}}
\newcommand{\caD}{{\mathcal D}}
\newcommand{\caG}{{\mathcal G}}
\newcommand{\caH}{{\mathcal H}}
\newcommand{\caK}{{\mathcal K}}
\newcommand{\caL}{{\mathcal L}}
\newcommand{\caO}{{\mathcal O}}
\newcommand{\caR}{{\mathcal R}}
\newcommand{\bbE}{{\mathbb E}}
\newcommand{\bbP}{{\mathbb P}}
\newcommand{\bbR}{{\mathbb R}}
\newcommand{\bflambda}{\bolds{\lambda}}
\newcommand{\bfy}{\mathbf{y}}
\newcommand{\bsm}{\mathbf{m}}
\newcommand{\bsI}{\mathbf{I}}
\renewcommand{\d}{{\mathrm d}}
\newcommand{\lone}{\mathbh{1}}
\newcommand{\dd}{\mathrm{d}}
\newcommand{\ii}{\mathrm{i}}
\newcommand{\bfx}{\mathbf{x}}
\def\Xint#1{\mathchoice
 {\XXint\displaystyle\textstyle{#1}}%
 {\XXint\textstyle\scriptstyle{#1}}%
 {\XXint\scriptstyle\scriptscriptstyle{#1}}%
 {\XXint\scriptscriptstyle\scriptscriptstyle{#1}}%
 \!\int}
 \def\XXint#1#2#3{{\setbox0=\hbox{$#1{#2#3}{\int}$}
 \vcenter{\hbox{$#2#3$}}\kern-.5\wd0}}
 
 \def\dashint{\Xint-}
\makeatother

\begin{document}
\begin{frontmatter}

\title{Bulk universality for deformed Wigner matrices}
\runtitle{Bulk universality for deformed Wigner matrices}

\begin{aug}
\author[A]{\fnms{Ji Oon}~\snm{Lee}\thanksref{M1,T1}\ead[label=e1]{jioon.lee@kaist.edu}},
\author[B]{\fnms{Kevin}~\snm{Schnelli}\corref{}\thanksref{M2,T2,T5}\ead[label=e2]{kevin.schnelli@ist.ac.at}},
\author[C]{\fnms{Ben}~\snm{Stetler}\thanksref{M3,T3}\ead[label=e3]{bstetler@math.harvard.edu}}
\and
\author[C]{\fnms{Horng-Tzer}~\snm{Yau}\thanksref{M3,T4,T5}\ead[label=e4]{htyau@math.harvard.edu}}
\runauthor{Lee, Schnelli, Stetler and Yau}
\affiliation{KAIST\thanksmark{M1}, IST Austria\thanksmark{M2} and
Harvard University\thanksmark{M3}}
\address[A]{J. O. Lee\\
Department of Mathematical Sciences\\
KAIST\\
291 Daehak-ro Yuseong-gu\\
Daejeon 305-701\\
Republic of Korea\\
\printead{e1}}
\address[B]{K. Schnelli\\
Institute of Science and Technology Austria\\
Am Campus 1\\
3400 Klosterneuburg\\
Republic of Austria\\
\printead{e2}}
\address[C]{B. Stetler\\
H.-T. Yau\\
Department of Mathematics\\
Harvard University\\
One Oxford Street\\
Cambridge, Massachusetts 02138\\
USA\\
\printead{e3}\\
\phantom{E-mail:\ }\printead*{e4}}
\end{aug}
\thankstext{T1}{Supported in part by National Research Foundation of
Korea Grant 2011-0013474 and TJ Park Junior Faculty Fellowship.}
\thankstext{T2}{Supported by ERC Advanced Grant RANMAT, No. 338804,
and the ``Fund for Math.''}
\thankstext{T3}{Supported by NSF GRFP Fellowship DGE-1144152.}
\thankstext{T4}{Supported in part by NSF Grant DMS-13-07444 and
Simons investigator fellowship.}
\thankstext{T5}{This work was completed while the authors were
visiting the Institute for Advanced Study.}

%
\received{\smonth{6} \syear{2014}}
%
\revised{\smonth{2} \syear{2015}}

%
\begin{abstract}
We consider $N\times N$ random matrices of the form $H = W + V$
where~$W$ is a real symmetric or complex
Hermitian Wigner matrix and~$V$ is a random or deterministic, real,
diagonal matrix whose entries are
independent of~$W$. We assume subexponential decay for the matrix
entries of~$W$, and we choose~$V$
so that the eigenvalues of~$W$ and~$V$ are typically of the same order.
For a large class of diagonal
matrices~$V$, we show that the local statistics in the bulk of the
spectrum are universal in the limit of large~$N$.
\end{abstract}

%
\end{frontmatter}

\section{Introduction}\label{sectionintroduction}

A prominent class of random matrix models is the\break Wigner ensemble,
consisting of $N\times N$ real symmetric
or complex Hermitian matrices, $W = (w_{ij} )$, whose matrix entries
are random variables that are
independent up to the symmetry constraint $W = W^*$. The first rigorous
result about the spectrum of random matrices of this type is Wigner's
global semicircle law~\cite{Wigner}, which states that the empirical
distribution of the rescaled eigenvalues, $(\lambda_i)$, of a Wigner
matrix~$W$ is given by
%
\begin{equation}
\label{convergencetosemicirclelaw} \frac{1}{N}\sum_{i=1}^N
\delta_{\lambda_i}(E)\longrightarrow\rho_{\mathrm{sc}}(E): =\frac{1}{2\pi}\sqrt{\bigl(4-E^2\bigr)_+}\qquad(E\in\R),
\end{equation}
as $N\to\infty$, in the weak sense. The distribution $\rho_{\mathrm{sc}}$ is
called the \emph{semicircle law}.

Let $p_W^N(\lambda_1,\ldots,\lambda_N)$ denote the joint probability
density of the (unordered) eigenvalues of~$W$.
If the entries of the Wigner matrix~$W$ are i.i.d. (independent and
identically distributed) real or complex Gaussian random variables, the
joint density of the eigenvalues, $p_W^N\equiv p_{G}^N$, is given by
%
\begin{equation}
\label{theGUE} p_{G}^N(\lambda_1,\ldots,\lambda_N)=\frac{1}{Z_{G}^N}
\prod_{i<j}
\llvert\lambda_i-\lambda_j\rrvert^{\beta}
\mathrm{e}^{-\beta N\sum_{i=1}^N
\lambda_{i}^2/4},
\end{equation}
with $\beta=1,2$, for the real, complex case, respectively. The
normalization $Z_{G}^N\equiv Z_{G}^N(\beta)$ in~(\ref{theGUE}) can
be computed explicitly. The real and complex Gaussian matrix ensembles
so defined are known as the Gaussian orthogonal ensemble (GOE, \mbox{$\beta
=1$}) and Gaussian unitary ensemble (GUE, $\beta=2$), respectively, and
as noted above we denote the corresponding joint densities as $p_{G}^N$
instead of $p_W^N$.

The \emph{$n$-point correlation functions} are defined by
\[
\varrho_{W,n}^N(\lambda_1,\ldots,
\lambda_n): =\int_{\R^{N-n}}
p_W^N(\lambda_1,\lambda_2,
\ldots,\lambda_N)\,\dd\lambda_{n+1} \,\dd\lambda_{n+2}
\cdots\,\dd\lambda_{N},
\]
$1\le n\le N$. Using orthogonal polynomials the correlation functions
of the GUE and GOE have been explicitly computed by Dyson, Gaudin and
Mehta; see, for example,~\cite{Mbook}. For the Gaussian unitary
ensemble, their results assert that the limiting behavior on small
scales at a fixed energy $E$ in the bulk of the spectrum, that is, for
$\llvert E\rrvert <2$, satisfies
%
\begin{eqnarray}
\label{GUEcorrelationfunctions} && \frac{1}{[\rho_{\mathrm{sc}}(E)]^n} \varrho
_{G,n}^N \biggl( E+
\frac{\alpha
_1}{\rho_{\mathrm{sc}}(E)N},E+\frac{\alpha_2}{\rho_{\mathrm{sc}}(E)N},\ldots,E+\frac
{\alpha_n}{\rho_{\mathrm{sc}}(E)N} \biggr)
\nonumber\\[-8pt]\\[-8pt]\nonumber
&&\qquad \longrightarrow\det\bigl(K(\alpha_i-\alpha_j)
\bigr)_{i,j=1}^n,
\end{eqnarray}
as $N \to\infty$, where $K$ is the sine-kernel
\begin{eqnarray*}
K(x,y) &:=& \frac{\sin\pi(x-y)}{\pi(x-y)}.
\end{eqnarray*}
Note that the limit in~(\ref{GUEcorrelationfunctions}) is
independent of the energy $E$ as long as $E$ is in the bulk of the spectrum.
The rescaling by a factor $1/N$ of the correlation functions in~(\ref
{GUEcorrelationfunctions}) corresponds to the typical separation of
consecutive eigenvalues, and we refer to the law under such a scaling
as \emph{local statistics}. Similar but
more complicated formulas were also obtained for the Gaussian
orthogonal ensemble; see, for example,~\cite{AGZ,Mbook} for reviews.
Note that the limiting correlation functions do not factorize,
reflecting the fact that the eigenvalues remain strongly correlated in
the limit of large $N$.

The Wigner--Dyson--Gaudin--Mehta conjecture, or bulk universality
conjecture, states that the local eigenvalue statistics of Wigner
matrices are universal in the sense that they depend only on the
symmetry class of the matrix, but are otherwise independent of the
details of the distribution of the matrix entries. The bulk
universality can be formulated in terms of weak convergence of
correlation functions or in terms of eigenvalue gap statistics. This
conjecture for all symmetry classes has been established in a series of
papers~\cite{EPRSY,ESY4,EKYY2,EYY,EKYY4,singlegap}. After this work
began, parallel results were obtained
for complex Hermitian matrices and certain symmetric matrices in~\cite
{TV1,TV2}; see \cite{EY} for a more detailed review.

In the present paper, we consider deformed Wigner matrices. A deformed
Wigner matrix, $H$, is an $N\times N$ random matrix of the form
%
\begin{equation}
\label{lamatrice} H = V + W,
\end{equation}
where~$V$ is a real, diagonal, random or deterministic matrix and~$W$
is a real symmetric or complex
Hermitian Wigner matrix independent of~$V$.
The matrices are normalized so that the eigenvalues
of~$V$ and~$W$ are order one. If the entries, $(v_i)$, of~$V$ are
random we may think of~$V$ as a ``random potential''; if the entries
of~$V$ are deterministic, matrices in the form of~(\ref{lamatrice})
are sometimes referred to as \textit{Wigner matrices with external source}.

Assuming that the empirical eigenvalue distribution of~$V$,
\[
\widehat\nu: =\frac{1}{N}\sum
_{i=1}\delta_{v_i},
\]
converges weakly, respectively, weakly in probability, to a nonrandom
measure, $\nu$, it was shown in~\cite{P} that the empirical
distribution of the eigenvalues of $H$ converges weakly in probability
to a deterministic measure. This measure depends on~$\nu$ and is thus
in general distinct from $\rho_{\mathrm{sc}}$. We refer to it as the \emph
{deformed semicircle law}, henceforth denoted by~$\rho_{\mathrm{fc}}$. There is
no explicit formula for~$\rho_{\mathrm{fc}}$ in terms of~$\nu$. Instead,~$\rho
_{\mathrm{fc}}$ is obtained as the solution of a functional equation for its
Stieltjes transform; see~(\ref{mfcequation}) below. It is known
that~$\rho_{\mathrm{fc}}$ admits a density~\cite{B}. Depending on~$\nu
$,~$\rho_{\mathrm{fc}}$ may be supported on several disjoint intervals. For
simplicity, we assume below that~$\nu$ is such that~$\rho_{\mathrm{fc}}$ is
supported on a single bounded interval. Further, we choose $\widehat
\nu$ such that all eigenvalues of $H$ remain close to the support of
$\rho_{\mathrm{fc}}$; that is, there are no ``outliers'' for~$N$ sufficiently large.

If~$W$ belongs to the GUE, $H$ is said to belong to the deformed GUE.
The deformed GUE for the special case when~$V$ has two eigenvalues~$\pm a$, each with equal multiplicity, has been treated in a series of
papers~\cite{BK1,BK2,BK3}. In this setting the local eigenvalue
statistics of~$H$ can be obtained via the solution to a
Riemann--Hilbert problem; see also~\cite{CW} for the case when~$V$ has
equispaced eigenvalues. Bulk universality for correlation functions of
the deformed GUE with rather general deterministic or random~$V$ has
been proved in~\cite{S1} by means of the Brezin--Hikami/Johansson
integration formula.

In the present paper, we establish bulk universality of local averages
of correlation functions for deformed Wigner matrices of the form
$H=V+W$, where~$W$ is a real symmetric or complex Hermitian Wigner
matrix and~$V$ is a deterministic or random real diagonal matrix. We
assume that the entries of~$W$ are centered independent random
variables with variance $1/N$ whose distributions decay
sub-exponentially; see Definition~\ref{assumptionwigner}. If~$V$ is
random, we assume for simplicity that its entries $(v_i)$ are i.i.d.
random variables. We assume that~$\widehat\nu$ converges weakly,
respectively, weakly in probability, to a nonrandom measure~$\nu$; see
Assumption~\ref{assumptionmuVconvergence}. We further assume that
the corresponding deformed semicircle law~$\rho_{\mathrm{fc}}$ is supported on
a single compact interval and has square root decay at both endpoints.
Sufficient conditions for these assumptions to hold have appeared
in~\cite{S2} and are rephrased in Assumption~\ref{assumptionmuV}.
Under these
assumptions, our main results in Theorem~\ref{theorem1} and in
Theorem~\ref{theorem2} assert that the limiting correlation functions
of the deformed Wigner ensemble are universal when averaged over a
small energy window. Note that our results hold for complex Hermitian
and real symmetric deformed Wigner matrices.

Before we outline our proofs, we recall the notion of $\beta
$-\textit{ensemble} or \textit{log-gas} which generalizes the measures
in~(\ref
{theGUE}). Let $U$ be a real-valued potential, and consider the measure
on $\R^N$ defined by the density
%
\begin{equation}
\label{lebetaensemble} \mu_{U}^N(\lambda_1,\ldots,
\lambda_N): =\frac
{1}{Z_{U}^N} \prod
_{i<j} \llvert\lambda_i-\lambda_j
\rrvert^{\beta}\mathrm{e}^{-{\beta N}\sum
_{i=1}^N (\lambda_{i}^2/2+U(\lambda_i) )/2},
\end{equation}
where $\beta>0$ and $Z_{U}^N\equiv Z_{U}^N(\beta)$ is a
normalization. Bulk universality for \mbox{$\beta$-}ensemb\-les asserts that
the local correlation functions for measures in the form of~(\ref
{lebetaensemble}) are universal (for sufficiently regular potentials $U$)
in the sense that for each value of $\beta>0$ they agree with the
local correlation functions of the Gaussian ensemble with $U\equiv0$.

For the classical values $\beta\in\{1,2,4\}$, the eigenvalue
correlation functions of $\mu_U^N$ can be explicitly expressed in
terms of polynomials orthogonal to the exponential weight in~(\ref
{lebetaensemble}).
Thus the analysis of the correlation functions relies on the asymptotic
properties of the
corresponding orthogonal polynomials.
This approach, initiated by Dyson, Gaudin and Mehta (see~\cite{Mbook}
for a review), was the starting point for many results on
the universality for $\beta$-ensemble with $\beta\in\{1,2,4\}$~\cite
{BI,DKMcLVZ,DKMcLVZ2,LL,L,DG,KSh,ShMa1}.

For general $\beta>0$, bulk universality of $\beta$-ensembles has
been established in~\cite{BEYI,BEYII,BEY} for potentials $U\in C^4$.
Recently, alternative approaches to bulk universality for $\beta
$-ensembles with general $\beta$ have been presented in~\cite{ShMa}
and~\cite{BAG} under different conditions on $U$.

We emphasize at this point that the eigenvalue distributions of the
deformed ensemble in~(\ref{lamatrice}) are in general not of the
form~(\ref{lebetaensemble}), even when~$W$ belongs to the GUE or the GOE.

Returning to the random matrix setting, we recall that the general
approach to bulk universality for (generalized) Wigner matrices
in~\cite{EPRSY,ESY4,EYY} consists of three steps:
\begin{enumerate}[(2)]
\item[(1)] establish a local semicircle law for the density of eigenvalues;

\item[(2)] prove universality of Wigner matrices with a small Gaussian
component by analyzing the convergence of Dyson Brownian motion to
local equilibrium;

\item[(3)] compare the local statistics of Wigner ensembles with
Gaussian divisible ensembles to remove the small Gaussian component of step~(2).
\end{enumerate}

For an overview of recent results and this three-step strategy,
see~\cite{EY}. Note that the ``local equilibrium'' in step~(2) refers
to measure (\ref{theGUE}), with $\beta=1,2$, respectively, in the real
symmetric, complex Hermitian case.

For deformed Wigner matrices, the \emph{local deformed semicircle
law}, the analogue of step~(1), was established in~\cite{LS} for
random~$V$. However, when~$V$ is random, the eigenvalues of~$V+W$
fluctuate on scale~$N^{-1/2}$ in the bulk (see~\cite{LS}), but their
gaps remain rigid on scale $N^{-1}$. To circumvent the mesoscopic
fluctuations of the eigenvalue positions, we condition on~$V$,
considering its entries to be fixed. The methods of~\cite{LS} can be extended, as outlined in Section~\ref
{sectionlocallaw}, to prove a local law on the optimal scale for ``typical''
realizations of random as well as deterministic potentials~$V$.

Our corresponding version of step~(2), a proof of bulk universality for
deformed Wigner ensembles with small Gaussian component, is the main
novelty of this paper. The local equilibrium of Dyson Brownian motion
in the deformed case is unknown but may effectively be approximated by
a ``reference'' $\beta$-ensemble that we explicitly construct in
Section~\ref{sectionbetaensemble}. In Section~\ref{sectionDBM}, we
analyze the convergence of the local distribution of the deformed
Wigner ensemble under Dyson Brownian motion to the ``reference'' $\beta
$-ensemble. However, since the ``reference'' $\beta$-ensemble is not
given by the invariant GUE/GOE, it also evolves in time.
Using the rigidity estimates for the deformed ensemble established in
step~(1) and the rigidity estimates for general $\beta$-ensembles
established in~\cite{BEY}, we obtain, in Section~\ref{sectionDBM},
bounds on the time evolution of the relative entropy between the two
measures being compared. The idea to estimate the entropy flow of the
Dyson Brownian motion with respect to the ``global equilibrium state''
given by the GUE/GOE was initiated in~\cite{ESY4} and~\cite{ESYY}. On
the other hand, the idea to use ``time dependent local equilibrium
states'' to control the entropy flow of hydrodynamical equations was
introduced in~\cite{Y}. There it is observed that the change of
relative entropy is negligible provided that the time dependent local
equilibrium is chosen in agreement with the density predicted by the
hydrodynamical equations. In this paper, we combine both methods to
yield an effective estimate on the entropy flow of the Dyson Brownian
motion in the deformed case. This global entropy estimate is
then used in Section~\ref{localequilibriummeasures} to conclude that
the local statistics of the locally-constrained deformed ensemble with
small Gaussian component agree with those of the locally-constrained
reference $\beta$-ensemble. Relying
on the
main technical result of~\cite{singlegap}, we
further conclude that the local statistics of the locally-constrained
reference $\beta$-ensemble agrees with the local statistics of the
GUE/GOE. Once this conclusion is obtained for the locally-constrained
ensembles, it can be extended to the nonconstrained ensembles. This
completes step~(2) in the deformed case.

In Sections~\ref{Fromgapstatisticstocorrelationfunctions}
and~\ref{Proofsofmainresults}, we outline step~(3) for deformed
Wigner matrices; the proof is similar to the argument for Wigner
matrices in~\cite{EYY1}. The main technical input is a bound on the
resolvent entries of $H$ on scales $N^{-1-\varepsilon}$ that can be
obtained from the local law in step~(1). In Section~\ref
{Proofsofmainresults}, we then combine steps~(1)--(3) to conclude the
proof of
our main results, Theorems~\ref{theorem1}~and~\ref{theorem2}.

We remark that our arguments in step~(2) do not rely on $V$ being
diagonal. Step~(3) depends only on the
deformed local semicircle law of step~(1); in principle, step~(3) is
independent of whether or not $V$ is diagonal, as long as a
deformed local semicircle law is given. Currently, our proof for the
deformed local semicircle law uses that $V$ is diagonal.

In Section~\ref{sectionedgeuniversality}, we prove that, in addition
to bulk universality,
the edge universality also holds for our model, that is, that the local
statistics at the spectral edges are given by the Tracy--Widom--Airy
statistics. From the main technical result of~\cite{BEY}, the proof of the edge universality follows the same
three-step program as the proof of bulk universality. A detailed
discussion of our edge universality result, Theorem~\ref
{theedgeuniversalitytheorem}, and related results can be found in
Section~\ref
{sectionresultsonedgeuniversality}.

In the \hyperref[app]{Appendix}, we collect several technical results on the deformed
semicircle law and its Stieltjes transform. Some of these results have
previously appeared in~\cite{S2} and~\cite{LS,LS2}.

\section{Assumptions and main results}\label{sectionassumptionsandmainresults}
In this section, we list our assumptions and our main results.

\subsection{Definition of the model}

We first introduce real symmetric and complex Hermitian Wigner matrices.

%
\begin{definition}\label{assumptionwigner}
A real symmetric Wigner matrix is an $N\times N$ random matrix,~$W$,
whose entries, $(w_{ij})$ $(1\le i,j\le N)$, are independent (up to
the symmetry constraint $w_{ij}={w}_{ji}$) real centered random
variables satisfying
%
\begin{equation}
\label{momentestimates} \E w_{ii}^2 = \frac{2}{N}, \qquad\E
w_{ij}^2 = \frac{1}{N} \qquad(i\neq j).
\end{equation}
In case $(w_{ij})$ are Gaussian random variables,~$W$ belongs to the
Gaussian orthogonal ensemble (GOE).

A complex Hermitian Wigner matrix is an $N\times N$ random matrix,~$W$,
whose entries, $(w_{ij})$ $(1\le i,j\le N)$, are independent (up to
the symmetry constraint $w_{ij}=\bar{w}_{ji}$) complex centered
random variables satisfying
%
\begin{equation}
\label{wignercomplex} \E w_{ii}^2 = \frac{1}{N},\qquad\E
\llvert w_{ij}\rrvert^2 = \frac{1}{N}, \qquad\E
w_{ij}^2 = 0 \qquad(i \neq j).
\end{equation}
For simplicity, we assume that the real and imaginary parts of
$(w_{ij})$ are independent for all $i,j$. This\vspace*{1pt} ensures that $\E
w_{ij}^2=0$ $(i\neq j)$. In case $(\re w_{ij})$ and $(\im w_{ij})$ are
Gaussian random variables,~$W$ belongs to the Gaussian unitary ensemble (GUE).

Irrespective of the symmetry class of~$W$, we assume that the
entries~$(w_{ij})$ have a subexponential decay, that is,
%
\begin{equation}
\label{eqC0} \mathbb{P} \bigl(\sqrt{N} \llvert w_{ij}\rrvert>x
\bigr)\le C_0 \mathrm{e}^{-x^{1/\theta}},
\end{equation}
for some positive constants $C_0$ and $\theta>1$. In particular,
%
\begin{equation}
\E\llvert w_{ij}\rrvert^p\le C\frac{(\theta p)^{\theta p}}{N^{p/2}}
\qquad(p\ge3).
\end{equation}
\end{definition}

Let $V=\operatorname{diag}(v_i)$ be an $N\times N$ diagonal, random or
deterministic matrix, whose entries $(v_i)$ are real-valued. We denote
by $\widehat\nu$ the empirical eigenvalue distribution of the
diagonal matrix~$V=\operatorname{diag}(v_i)$,
%
\begin{equation}
\label{widehatmuV} \widehat\nu: =\frac{1}{N}\sum
_{i=1}^N\delta_{v_i}.
\end{equation}

%
\begin{assumption}\label{assumptionmuVconvergence}
There is a (nonrandom) centered, compactly supported probability
measure $\nu$ such that the following holds:
\begin{longlist}[(2)]
\item[(1)] If~$V$ is a \textit{random} matrix, we assume that $(v_i)$
are independent and identically distributed real random variables with
law $\nu$. Further, we assume that $(v_i)$ are independent of $(w_{ij})$.

\item[(2)] If~$V$ is a \textit{deterministic} matrix, we assume that
there is $\alpha_0>0$, such that for any fixed compact set $\caD
\subset\C^+$ (independent of $N$) with $\operatorname{dist}(\caD,\break \supp
\nu)>0$, there is $C$ such that
%
\begin{equation}
\label{equationassumptionmuVconvergence} \max_{z\in\caD}\biggl\llvert
\int\frac{\dd\widehat\nu(v)}{v-z}-
\int\frac{\dd\nu(v)}{v-z} \biggr\rrvert\le CN^{-\alpha_0},
\end{equation}
for $N$ sufficiently large.
\end{longlist}
\end{assumption}

Note that~(\ref{equationassumptionmuVconvergence}) implies that
$\widehat\nu$ converges to $\nu$ in the weak sense as $N\to\infty
$. Also note that condition~(\ref{equationassumptionmuVconvergence})
holds for large $N$ with high probability for $0<\alpha
_0<1/2$ if $(v_i)$ are i.i.d. random variables.

\subsection{Deformed semicircle law}\label{subsectiondeforemdsemicirlelaw}

The deformed semicircle can be described in terms of the Stieltjes
transform: for a (probability) measure $\omega$ on the real line we
define its Stieltjes transform, $m_{\omega}$, by
\begin{eqnarray*}
m_{\omega}(z): =\int\frac{\dd\omega(v)}{v-z} \qquad\bigl
(z\in\C
^+\bigr).
\end{eqnarray*}
Note that $m_{\omega}$ is an analytic function in the upper half plane
and that $\im m_\omega(z)\ge0$, $\im z> 0$. Assuming that $\omega$
is absolutely continuous with respect to Lebesgue measure, we can
recover the density of $\omega$ from $m_{\omega}$ by the inversion formula
%
\begin{eqnarray}
\label{stieltjesinversionformula} \omega(E)=\lim_{\eta\searrow0}\frac
{1}{\pi}\im
m_{\omega}(E+\ii\eta) \qquad(E\in\R).
\end{eqnarray}
We use the same symbols to denote measures and their densities.
Moreover, we have
\begin{eqnarray*}
\lim_{\eta\searrow0} \re m_{\omega}(E+\ii\eta)=\dashint
\frac
{\omega(v)\,\dd v}{v-E}\qquad(E\in\R),
\end{eqnarray*}
whenever the left-hand side exists. Here the integral on the right is
understood as principal value integral. We denote in the following by
$\re m_{\omega}(E)$ and $\im m_{\omega}(E)$ the limiting quantities
%
\begin{eqnarray}
\label{abuseofnotationI}
\re m_{\omega}(E)&\equiv&\lim_{\eta\searrow0}\re
m_{\omega}(E+\ii\eta),
\nonumber\\[-8pt]\\[-8pt]\nonumber
\im m_{\omega}(E)&\equiv&\lim
_{\eta\searrow0}\im m_{\omega}(E+\ii\eta),
\end{eqnarray}
$E\in\R$, whenever the limits exist.

Choosing $\omega$ to be the standard semicircular law $\rho_{\mathrm{sc}}$,
the Stieltjes transform $m_{\rho_{\mathrm{sc}}}\equiv m_{\mathrm{sc}}$ can be computed
explicitly, and one checks that $m_{\mathrm{sc}}$ satisfies the relation
\begin{eqnarray*}
m_{\mathrm{sc}}(z) &=&\frac{-1}{m_{\mathrm{sc}}(z)+z},\qquad\im m_{\mathrm{sc}}(z)\ge0
\qquad\bigl(z\in\C^+\bigr).
\end{eqnarray*}

The deformed semicircle law is conveniently defined through its
Stieltjes transform. Let $\nu$ be the limiting probability measure of
Assumption~\ref{assumptionmuVconvergence}. Then it is well
known~\cite{P} that the functional equation
%
\begin{equation}
\label{mfcequation} m_{\mathrm{fc}}(z)= \int\frac{\dd\nu
(v)}{v-z-m_{\mathrm{fc}}(z)},\qquad\im
m_{\mathrm{fc}}(z)\ge0\qquad\bigl(z\in\C^+\bigr),
\end{equation}
has a unique solution, also denoted by $m_{\mathrm{fc}}$, that satisfies, for
all $E\in\R$, $\limsup_{\eta\searrow0}\im m_{\mathrm{fc}}(E+\ii\eta
)<\infty$. Indeed, from~(\ref{mfcequation}), we obtain that
%
\begin{equation}
\label{sumrule} \int\frac{\dd\nu(v)}{\llvert v-z-m_{\mathrm{fc}}(z)\rrvert
^2}=\frac{\im m_{\mathrm{fc}}(z)}{\im
m_{\mathrm{fc}}(z)+\eta}\le1\qquad\bigl(z\in
\C^+\bigr),
\end{equation}
thus $\llvert m_{\mathrm{fc}}(z)\rrvert \le1$, for all $z\in\C^+$.

The deformed semicircle law, denoted by $\rho_{\mathrm{fc}}$, is then defined
through its density
\begin{eqnarray*}
\rho_{\mathrm{fc}}(E): =\lim_{\eta\searrow0}
\frac{1}{\pi
}\im m_{\mathrm{fc}}(E+\ii\eta)\qquad(E\in\R).
\end{eqnarray*}
The measure $\rho_{\mathrm{fc}}$ has been studied in detail in~\cite{B}. For
example, it was shown there that the density $\rho_{\mathrm{fc}}$ is an
analytic function inside the support of the measure.

The measure $\rho_{\mathrm{fc}}$ is also referred to as the additive free
convolution of the semicircular law and the measure~$\nu$. More
generally, the additive free convolution of two (probability) measures
$\omega_1$ and $\omega_2$, usually denoted by $\omega_1\boxplus
\omega_2$, is defined as the distribution of the sum of two freely
independent noncommutative random variables, having distributions
$\omega_1$,~$\omega_2$, respectively; we refer, for example, to~\cite
{VDN,AGZ} for reviews. Similar to~(\ref{mfcequation}), the free
convolution measure $\omega_1\boxplus\omega_2$ can be described in
terms of a set of functional equations for the Stieltjes transforms;
see \mbox{\cite{CG,BeB07}}.

Our second assumption on $\nu$ guarantees (see Lemma~\ref{lemmamfc}
below) that~$\rho_{\mathrm{fc}}$ is supported on a single interval and
that~$\rho_{\mathrm{fc}}$ has a square root behavior at the two endpoints of
its support. Sufficient conditions for this behavior have been
presented in~\cite{S2}. The assumptions below also rule out the
possibility that the matrix $H$ has ``outliers'' in the limit of large $N$.

\begin{assumption}\label{assumptionmuV}
Let $I_{\nu}$ be the smallest interval such that $\supp{\nu
}\subseteq I_{\nu}$. Then there exists $\varpi>0$ such that
%
\begin{eqnarray}
\label{eqassumptionmuV} \inf_{x\in I_{\nu}}\int\frac{\dd\nu
(v)}{(v-x)^2}\ge1+\varpi.
\end{eqnarray}
Similarly, let $I_{\widehat\nu}$ be the smallest interval such that
$\supp{\widehat\nu}\subseteq I_{\widehat\nu}$. Then:
\begin{longlist}[(2)]
\item[(1)] for \textit{random} $(v_i)$, there is a constant $\mathfrak
{t}>0$, such that
%
\begin{eqnarray}
\label{eqassumptionmuV2} \mathbb{P} \biggl(\inf_{x\in I_{\widehat\nu
}}\int
\frac{\dd
\widehat\nu(v)}{(v-x)^2}\ge1+{\varpi} \biggr)\ge1-N^{-\mathfrak
{t}},
\end{eqnarray}
for $N$ sufficiently large;

\item[(2)] for \textit{deterministic} $(v_i)$,
%
\begin{eqnarray}
\label{eqassumptionmuV3} \inf_{x\in I_{\widehat\nu}}\int\frac{\dd
\widehat\nu
(v)}{(v-x)^2}\ge1+{\varpi},
\end{eqnarray}
for $N$ sufficiently large.
\end{longlist}
\end{assumption}

We give two examples for which~(\ref{eqassumptionmuV}) is
satisfied:
\begin{longlist}[(2)]
\item[(1)] Choosing $\nu=\frac{1}{2}(\mathrm{\delta
}_{-a}+\mathrm{\delta}_{a})$, $a\ge0$, we have $I_{\nu}=[-a,a]$.
For $a<1$, one checks that there is a $\varpi=\varpi(a)$ such
that~(\ref{eqassumptionmuV}) is satisfied and that the deformed
semicircle law is supported on a single interval with a square root
type behavior at the edges. However, for $a>1$, the deformed semicircle
law is supported on two disjoint intervals; for further details,
see~\cite{BK1,BK2,BK3}.

\item[(2)] Let $\nu$ be a centered Jacobi measure of the form
%
\begin{eqnarray}
\label{Jacobimeasure} \nu(v)=Z^{-1}(v-1)^{\mathit{a}}(1-v)^{\mathit{b}}
d(v)\lone_{[-1,1]}(v),
\end{eqnarray}
where $d\in C^{1}([-1,1])$, $d(v)>0$, $-1<\mathit{a},\mathit{b}<\infty$ and $Z$, a normalization constant. Then for $\mathit{a},\mathit{b}<1$, there is $\varpi>0$ such that~(\ref
{eqassumptionmuV}) is satisfied with $I_{\nu}=[-1,1]$. However, if
$\mathit{a}>1$ or $\mathit{b}>1$, then~(\ref{assumptionmuV}) may
not be satisfied. In this setting the deformed semicircle law is still
supported on a single interval; however, the square root behavior at
the edge may fail. We refer to~\cite{LS,LS2} for a detailed discussion.
\end{longlist}

%
\begin{lemma}\label{lemmavorbereitung}
Let $\nu$ satisfy~(\ref{eqassumptionmuV}) for some $\varpi>0$.
Then there are $L_-,L_+$, with $L_-\le-2$, $2\le L_+$, such that
$\supp\rho_{\mathrm{fc}}=[L_-,L_+]$. Moreover,~$\rho_{\mathrm{fc}}$ has a strictly
positive density in $(L_-,L_+)$.
\end{lemma}

Lemma~\ref{lemmavorbereitung} follows directly from Lemma~\ref
{lemmamfc} below.

\subsection{Results on bulk universality}

Recall that we denote by $\varrho_{H,n}^N$ the $n$-point correlation
function of $H=V+W$, where~$V$ is either a real deterministic or real
random diagonal matrix. We denote by $\varrho_{G,n}^N$ the $n$-point
correlation function of the GUE, respectively, the GOE.\vadjust{\goodbreak}

A function $O\dvtx \R^{n}\to\R$ is called an $n$-particle \emph
{observable} if $O$ is symmetric, smooth and compactly supported.
Recall from Lemma~\ref{lemmavorbereitung} that we denote by $L_\pm$
the endpoints of the support of the measure $\rho_{\mathrm{fc}}$. For
deterministic~$V$ we have the following result.

\begin{theorem}\label{theorem1}
Let~$W$ be a complex Hermitian or a real symmetric Wigner matrix
satisfying the assumptions in Definition~\ref{assumptionwigner}.
Let~$V$ be a deterministic real diagonal matrix satisfying
Assumptions \ref{assumptionmuVconvergence}~and~\ref
{assumptionmuV}. Set $H=V+W$. Let $E,E'$ be two energies satisfying
$E\in(L_-,L_+)$, $E'\in(-2,2)$. Fix $n\in\N$, and let $O$ be an
$n$-particle observable. Let $\delta>0$ be arbitrary, and choose
$b\equiv b_N$ such that $N^{-\delta}\ge b_{N}\ge N^{-1+\delta}$. Then
%
\begin{eqnarray}
\label{theorem1equation1}
&& \lim_{N\to\infty}\int_{\bbR^n} \dd \alpha_1\cdots\,\dd\alpha_n O(\alpha_1,\ldots,
\alpha_n)\nonumber
\\
&&\hspace*{41pt}{}\times \biggl[\frac
{1}{2b}\int_{E-b}^{E+b}
\frac{\dd x}{[\rho_{\mathrm{fc}}(E)]^n}\varrho^N_{H,n} \biggl(x+
\frac{\alpha_1}{\rho_{\mathrm{fc}}(E)N},\ldots,x+\frac
{\alpha_n}{\rho_{\mathrm{fc}}(E)N} \biggr)
\nonumber\\[-8pt]\\[-8pt]\nonumber
&&\hspace*{72pt}{} -\frac{1}{[\rho_{\mathrm{sc}}(E')]^n}\varrho^N_{G,n} \biggl(E'+
\frac{\alpha
_1}{\rho_{\mathrm{sc}}(E')N},\ldots,E'+\frac{\alpha_n}{\rho_{\mathrm{sc}}(E')N} \biggr)
\biggr]\hspace*{-10pt}
\\
&&\qquad{} = 0,\nonumber
\end{eqnarray}
where $\rho_{\mathrm{fc}}$ denotes the density of the deformed semicircle law
and $\rho_{\mathrm{sc}}$ denotes the density of the standard semicircle law.
Here, $\varrho_{G,n}^N$ denotes the $n$-point correlation function of
the GUE in case~$W$ is a complex Hermitian Wigner matrix, respectively,
the $n$-point correlation function of the GOE in case~$W$ is a real
symmetric Wigner matrix.
\end{theorem}

For random~$V$ we have the following result.

%
\begin{theorem}\label{theorem2}
Let~$W$ be a complex Hermitian or a real symmetric Wigner matrix
satisfying the assumptions in Definition~\ref{assumptionwigner}.
Let~$V$ be a random real diagonal matrix whose entries are i.i.d.
random variables that are independent of~$W$ and satisfy
Assumptions \ref{assumptionmuVconvergence}~and~\ref
{assumptionmuV}. Set $H=V+W$. Let $E,E'$ be two energies satisfying
$E\in(L_-,L_+)$, $E'\in(-2,2)$. Fix $n\in\N$, and let $O$ be an
$n$-particle observable. Let $\delta>0$ be arbitrary, and choose
$b\equiv b_N$ such that $N^{-\delta}\ge b_{N}\ge N^{-1/2+\delta}$. Then
%
\begin{eqnarray}
\label{theorem2equation1}
&& \lim_{N\to\infty}\int_{\bbR^n} \dd \alpha_1\cdots\,\dd\alpha_n O(\alpha_1,\ldots,
\alpha_n)\nonumber
\\
&&\hspace*{41pt}{}\times \biggl[\frac
{1}{2b}\int_{E-b}^{E+b}
\frac{\dd x}{[\rho_{\mathrm{fc}}(E)]^n}\varrho^N_{H,n} \biggl(x+
\frac{\alpha_1}{\rho_{\mathrm{fc}}(E)N},\ldots,x+\frac
{\alpha_n}{\rho_{\mathrm{fc}}(E)N} \biggr)
\nonumber\\[-8pt]\\[-8pt]\nonumber
&&\hspace*{72pt}{} -\frac{1}{[\rho_{\mathrm{sc}}(E')]^n}\varrho^N_{G,n} \biggl(E'+
\frac{\alpha
_1}{\rho_{\mathrm{sc}}(E')N},\ldots,E'+\frac{\alpha_n}{\rho_{\mathrm{sc}}(E')N} \biggr)
\biggr]\hspace*{-10pt}
\\
&&\qquad = 0,\nonumber
\end{eqnarray}
where $\rho_{\mathrm{fc}}$ denotes the density of the deformed semicircle law
and $\rho_{\mathrm{sc}}$ denotes the density of the standard semicircle law.
Here, $\varrho_{G,n}^N$ denotes the $n$-point correlation function of
the GUE in case~$W$ is a complex Hermitian Wigner matrix, respectively,
the $n$-point correlation function of the GOE in case~$W$ is a real
symmetric Wigner matrix.
\end{theorem}
%

\begin{remark}
Theorem~\ref{theorem1} and Theorem~\ref{theorem2} show that the
averaged local correlation functions of $H=V+W$ are universal in the
limit of large~$N$ in the sense that they are independent of the
diagonal matrix~$V$ and also independent of the precise distribution of
the entries of~$W$. Both theorems hold for real symmetric and complex
Hermitian matrices. For the former choice, $\varrho_{G,n}^N$ stands
for the \mbox{$n$-}point correlation functions of the GOE. For the latter
choice, $\varrho_{G,n}^N$ stands for the $n$-point correlation
functions of the GUE.

Note that we can choose $b_N$ of order $N^{-1+\delta}$, $\delta>0$,
for deterministic~$V$ in Theorem~\ref{theorem1}, while we have to
choose $b_N$ of order $N^{-1/2+\delta}$, $\delta>0$, for random~$V$
in Theorem~\ref{theorem2}. The latter condition is technical and not
optimal. It is related to our next comment.

For random~$V$ with $(v_i)$ i.i.d. bounded random variables, the
eigenvalues of $H$ fluctuate on scale $N^{-1/2}$ in the bulk~\cite
{LS}. Yet, under the assumptions of Theorem~\ref{theorem2}, the
eigenvalue gaps remain rigid over small scales so that the universality
of local correlation functions, a statement about the eigenvalue gaps,
is unaffected by these mesoscopic fluctuations. We thus expect
Theorem~\ref{theorem2} to hold with $b_N\gg N^{-1}$. Relying on
explicit integration formulas in the complex Hermitian setting, we
suppose that the averaging over an energy window can be dropped; cf.
the results for the deformed GUE in~\cite{S1}.
\end{remark}

%
\begin{remark} \label{remarkforgapuniversality}
The main ingredient of our proofs of Theorem~\ref{theorem1} and
Theorem~\ref{theorem2} is an entropy
estimate; see Proposition~\ref{thesuperproposition}. Once such an
estimate is obtained, the
method in~\cite{singlegap} also implies the single gap universality
in the sense that the distribution of any single
gap in the bulk is the same (up to a scaling) as the one from the
corresponding Gaussian case.
More precisely, fix $\alpha>0$, and let $k\in\N$ be such that
$\alpha N\le k\le(1-\alpha)N$. Let $O$ be an $n$-particle observable.
Then there are $\chi>0$ and $C$ such that
\begin{eqnarray*}
&& \bigl\llvert\E^{H} O \bigl( (N\rho_{\mathrm{fc},k}) (
\lambda_k-\lambda_{k+1}),(N\rho_{\mathrm{fc},k}) (
\lambda_k-\lambda_{k+2}),\ldots,(N\rho_{\mathrm{fc},k}) (
\lambda_k-\lambda_{k+n}) \bigr)
\\
&&\quad{} - \E^{\mu_{G}} O
\bigl( (N\rho_{\mathrm{sc},k}) (\lambda_k-\lambda
_{k+1}),(N\rho_{\mathrm{sc},k}) (\lambda_k-
\lambda_{k+2}),\ldots,
\\
&&\hspace*{241pt} (N\rho_{\mathrm{sc},k}) (\lambda_k-
\lambda_{k+n}) \bigr) \bigr\rrvert
\nonumber
\\
&&\qquad\leq C N^{-\chi},
\end{eqnarray*}
\normalsize
for $N$ sufficiently large, where $\mu_G$ is the standard GOE or GUE
ensemble, depending
on the symmetry class of $H$. Here $\rho_{\mathrm{fc},k}$ stands for the
density of the measure $\rho_{\mathrm{fc}}$ at the classical location, $\gamma
_k$, of the $k$th eigenvalue defined through
%
\begin{eqnarray}
\int_{-\infty}^{\gamma_k}\rho_{\mathrm{fc}}(x)\,\dd x=
\frac{k-1/2}{N}.
\end{eqnarray}
Similarly, $\rho_{\mathrm{sc},k}$ stands for the density of the standard
semicircle law $\rho_{\mathrm{sc}}$ at the classical location of the $k$th
eigenvalue of the Gaussian ensembles.
\end{remark}
%
%
\begin{remark}
To conclude, we mention two extensions of the above results. In
Theorem~\ref{theorem2} we may relax the assumption that $(v_i)$ are
independent among themselves: our results can be extended to dependent
random variables provided that $(v_i)$ satisfy~(\ref
{equationassumptionmuVconvergence}),~(\ref{eqassumptionmuV}) and~(\ref
{eqassumptionmuV2}) for some constants $\alpha_0,\varpi,\mathfrak
{t}>0$, and provided that $(v_i)$ are independent of $(w_{ij})$. In
such a setting the required lower bound on~$b_N$ depends on $\alpha_0$.

The assumption that $V$ is diagonal can be relaxed by assuming in turn
that $W$ belongs to the GUE/GOE. Then using the invariance of $W$, we
can diagonalize $V$ and apply our approach for diagonal potentials. For
$W$ a Wigner matrix and $V$ a~nondiagonal matrix, we expect that
similar results hold by slowly changing W to a GUE/GOE. This, however,
involves many more technical steps.
\end{remark}

\subsection{Results on edge universality}\label{sectionresultsonedgeuniversality}

In this subsection, we show that our model also satisfies the edge
universality. Edge universality states that the statistics of the
extremal eigenvalues of many random matrix ensembles are universal: let
$\lambda_N$ denote the largest eigenvalue of a Wigner matrix $W$. The
limiting distribution of $\lambda_N$ was identified for the Gaussian
ensembles by Tracy and Widom~\cite{TW1,TW2}. They proved that
%
\begin{eqnarray}
\label{prototype} \lim_{N\to\infty} \mathbb{P}\bigl(N^{2/3}(
\lambda_N -2) \le s\bigr) = F_\beta(s)\qquad\qquad\bigl(
\beta\in\{1,2,4\}\bigr),
\end{eqnarray}
$s\in\R$, where the Tracy--Widom distribution functions $F_{\beta}$
are described by Painlev\'{e} equations. The edge universality can also
be extended to the $k$ largest eigenvalues, where the joint
distribution of the $k$ largest eigenvalues can be written in terms of
the Airy kernel, as first shown for the GUE/GOE in~\cite{F}. These
results also hold for the $k$ smallest eigenvalues.

Edge universality for Wigner matrices was first proved in~\cite{So1}
(see also~\cite{SiSo1}) for real symmetric and complex Hermitian
ensembles with symmetric distributions. The symmetry assumption on the
entries' distribution was partially removed in~\cite{PeSo1,PeSo2}.
Edge universality was proved in~\cite{TV2} under the condition that
the distribution of the matrix elements has subexponential decay, and
its first three
moments match those of the Gaussian distribution; that is, the third
moment of the entries vanish. The vanishing third moment condition was
removed in~\cite{EYY}. Finally, edge universality for generalized
Wigner matrices was proved only recently in~\cite{BEY}.

Edge universality for the deformed GUE was obtained for the special
case when~$V$ has two eigenvalues~$\pm a$, each with equal
multiplicity, via a Riemann--Hilbert approach in~\cite{BK1,BK2}. For
general $V$, the joint distribution of the eigenvalues of the deformed
GUE can be expressed explicitly by the Brezin--Hikami/Johansson formula
that may be used to prove the edge universality various choices and
ranges of $V$; see~\cite{J2,S1,CP}.

Our result on the edge universality for real symmetric and complex
Hermitian deformed Wigner matrices is as follows.

%
\begin{theorem}\label{theedgeuniversalitytheorem}
Let~$W$ be a complex Hermitian or a real symmetric Wigner matrix
satisfying the assumptions in Definition~\ref{assumptionwigner}.
Let~$V$ be either a random real diagonal matrix whose entries are
i.i.d. random variables that are independent of~$W$, or a
deterministic real diagonal matrix. Assume that $V$ satisfies
Assumptions \ref{assumptionmuVconvergence}~and~\ref{assumptionmuV}. Set $H=V+W$.

Then there are
$\varkappa> 0, \chi>0, c_0>0$ such that the following result holds
for any fixed $n\in\N$.
For any $n$-particle observable $O$ and for $\Lambda\subset\llbracket
1,N^\varkappa\rrbracket$,
respectively, $\Lambda\subset\llbracket N-N^{\varkappa},N\rrbracket$,
with $\llvert \Lambda\rrvert = n$, we have
%
\begin{eqnarray}\label{theedgeuniversalityequation}
&& \bigl\llvert\E^{H} O \bigl( \bigl
(c_0N^{ 2/3}
j^{1/3}(\lambda_j-\widehat\gamma_j)
\bigr)_{j\in\Lambda} \bigr) - \E^{\mu_{G}} O \bigl( \bigl(N^{ 2/3}
j^{1/3}(\lambda_j-\gamma_j)
\bigr)_{j\in\Lambda} \bigr) \bigr\rrvert\hspace*{-20pt}
\nonumber\\[-8pt]\\[-8pt]\nonumber
&&\qquad \leq C_O N^{-\chi},
\end{eqnarray}
for $N$ sufficiently large, for some constant $C_O$ (depending on $O$),
where~$\mu_G$ is the standard GUE/GOE, depending
on the symmetry class of~$W$. Here, the constant~$c_0$ is a scaling
factor so that
the eigenvalue density at the edge of~$H$ can be compared with the
Gaussian case. It only depends on~$\nu$. Further,~$\widehat\gamma
_j$,~$\gamma_j$ denote here the classical locations of the $j$th
eigenvalue with respect to the measure~$\widehat\rho_{\mathrm{fc}}$ introduced
in~(\ref{stieltjesinversionformula3}) below, respectively, with
respect to the standard semicircle law $\rho_{\mathrm{sc}}$.
\end{theorem}

Theorem~\ref{theedgeuniversalitytheorem} shows that the local
statistics of the $k$ largest, respectively, smallest, eigenvalues of
our model are given by the Tracy--Widom--Airy statistics.

The measure $\widehat\varrho_{\mathrm{fc}}$ depends solely on the empirical
eigenvalue distribution,~$\widehat\nu$, of~$V$, and so do the
classical locations $(\widehat\gamma_k)$. The scaling factor $c_0$
in~(\ref{theedgeuniversalityequation}) may be computed
explicitly~\cite{S1}.

Theorem~\ref{theedgeuniversalitytheorem} is proved in a similar way
to Theorems~\ref{theorem1} and~\ref{theorem2}. Using the Dirichlet
form bound obtained in Proposition~\ref{thesuperproposition} below,
we invoke the edge universality result for localized $\beta
$-ensembles, Theorem~3.3 of~\cite{BEY}, and follow the same strategy
as for the bulk universality. The proof of Theorem~\ref
{theedgeuniversalitytheorem} is given in Section~\ref{sectionedgeuniversality}.

To conclude, we mention that Theorem~\ref{theedgeuniversalitytheorem}
has recently been proved in~\cite{LS3}
using a completely different approach based on the Green function
comparison theorem; see, for example, \cite{EYY1} for earlier
ideas of using the Green function comparison for edge universality.

\subsection{Notation and conventions}
In this subsection, we introduce some more notation and conventions
used throughout the paper. For high probability estimates we use two
parameters $\xi\equiv\xi_N$ and $\varphi\equiv\varphi_N$: we let
%
\begin{eqnarray}
\label{eqxi} a_0<\xi\le A_0\log\log N,\qquad
\varphi=(\log N)^{C_1},
\end{eqnarray}
for some constants $a_0>2$, $A_0\ge10$, $C_1> 1$.

%
\begin{definition}
We say an event $\Xi$ has $(\xi,\upsilon)$-high probability if
\begin{eqnarray*}
\mathbb{P}\bigl(\Xi^c\bigr)\le\mathrm{e}^{-\upsilon(\log N)^{\xi}}
\qquad(
\upsilon>0),
\end{eqnarray*}
for $N$ sufficiently large. We say an event $\Xi$ has $\varsigma
$-exponentially high probability if
\begin{eqnarray*}
\mathbb{P}\bigl(\Xi^c\bigr)\le\mathrm{e}^{-N^{\varsigma}} \qquad(
\varsigma>0),
\end{eqnarray*}
for $N$ sufficiently large.
Similarly, for a given event $\Xi_0$ we say an event $\Xi$ holds with
$(\xi,\upsilon)$-high probability, respectively, $\varsigma
$-exponentially high probability, on $\Xi_0$, if
\begin{eqnarray*}
\mathbb{P}\bigl(\Xi^c\cap\Xi_0\bigr)\le
\mathrm{e}^{-\upsilon(\log N)^{\xi
}} \qquad(\upsilon>0),\qquad\mathbb{P}\bigl(
\Xi^c\cap\Xi_0\bigr)\le\mathrm{e}^{-N^{\varsigma}} \qquad
(\varsigma>0),
\end{eqnarray*}
respectively, for $N$ sufficiently large.
\end{definition}

For brevity, we occasionally say an event holds with exponentially high
probability, when we mean $\varsigma$-exponentially high probability.
We do not keep track of the explicit value of $\upsilon$ or $\varsigma
$ in the following, allowing $\upsilon$ and $\varsigma$ to decrease
from line to line such that $\upsilon,\varsigma>0$.

We use the symbols $\caO( \cdot)$ and $o( \cdot)$ for the
standard big-O and little-o notation. The notation $\caO, o$, $\ll
$, $\gg$, refers to the limit $N\to\infty$, if not indicated
otherwise. Here $a\ll b$ means $a=o(b)$. We use~$c$ and~$C$ to denote
positive constants that do not depend on $N$. Their value may change
from line to line. We write $a\sim b$ if there is $C\ge1$ such that
$C^{-1}\llvert b\rrvert \le\llvert a\rrvert \le C \llvert
b\rrvert $, and occasionally we write for
$N$-dependent quantities $a_N\lesssim b_N$ if there exist constants
$C,c>0$ such that $\llvert a_N\rrvert \le C(\varphi_N)^{c\xi}\llvert
b_N\rrvert $.

Finally, we abbreviate
\begin{eqnarray*}
\sum_{j}^{(i)}( \cdot)\equiv\mathop{\sum
_{j=1}}_{j\neq i}^N( \cdot),
\end{eqnarray*}
and we use double brackets to denote index sets, that is,
\[
\llbracket n_1, n_2 \rrbracket: =[n_1, n_2] \cap\Z, %
\]
for $n_1, n_2 \in\R$.

\section{Local law and rigidity estimates}\label{sectionlocallaw}
Recall the constant $\varpi>0$ in Assumption~\ref{assumptionmuV}.
Set $\varpi': =\varpi/10$. In this section we
consider the family
of interpolating random matrices
%
\begin{eqnarray}
\label{definitionTheta} H^{\vartheta}: =\vartheta
V+W,\qquad\vartheta
\in\Theta_{\varpi
}: =\bigl[0,1+\varpi'\bigr],
\end{eqnarray}
where~$V$ and~$W$ are chosen to satisfy Assumptions~\ref
{assumptionmuVconvergence} and~\ref{assumptionmuV}, respectively, the
assumptions in Definition~\ref{assumptionwigner}. Here~$\vartheta$
has the interpretation of a possibly \mbox{$N$-}dependent positive ``coupling
parameter.''

We define the resolvent or \emph{Green function}, $G^{\vartheta}(z)$,
and the averaged Green function, $m^{\vartheta}(z)$, of $H^{\vartheta
}$ by
%
\begin{eqnarray}
\label{definitionofthem} G^{\vartheta}(z)=\bigl(G^{\vartheta}_{ij}(z)
\bigr): =\frac
{1}{\vartheta
V+W-z}, \qquad m_N^{\vartheta}(z)
: =\frac{1}{N}\Tr G^{\vartheta
}(z),
\end{eqnarray}
$z\in\mathbb{C}^{+}$. Frequently, we abbreviate $G^{\vartheta}\equiv
G^{\vartheta}(z)$, $m_N^{\vartheta}\equiv m_N^{\vartheta}(z)$, etc.

To conveniently cope with the cases when $(v_i)$ are random,
respectively, deterministic, we introduce an event~$\Omega$ on which
the random variables $(v_i)$ exhibit ``typical'' behavior. Recall that
we denote by $m_{\widehat\nu}$ and $m_{\nu}$ the Stieltjes
transforms of~$\widehat\nu$, respectively, $\nu$.

\begin{definition}\label{definitionoftheeventomega}
Let $\Omega\equiv\Omega(N)$ be an event on which the following holds:
\begin{longlist}[(2)]
\item[(1)]There is a constant $\alpha_0>0$ such that, for any fixed
compact set $\caD\subset\C^+$ (independent of $N$) with $\operatorname
{dist}(\caD,\supp\nu)>0$, there is $C$ such that
%
\begin{eqnarray}
\label{definitionofomegaV} \bigl\llvert m_{\widehat\nu}(z)-m_{\nu}(z)
\bigr
\rrvert\le CN^{-\alpha_0},
\end{eqnarray}
for $N$ sufficiently large.

\item[(2)] Recall the constant $\varpi>0$ in Assumption~\ref
{assumptionmuV}. We have
%
\begin{eqnarray}
\label{definitionofomegaVIII} \inf_{x\in I_{\widehat\nu}}\int\frac{\dd
\widehat\nu
(v)}{(v-x)^2}\ge1+\varpi,
\qquad\inf_{x\in I_{\nu}}\int\frac
{\dd\nu}{(v-x)^2}\ge1+\varpi,
\end{eqnarray}
for $N$ sufficiently large.
\end{longlist}
\end{definition}

In case $(v_i)$ are deterministic, $\Omega$ has full probability for
$N$ sufficiently large by the Assumptions in~\ref{assumptionmuVconvergence}.

Similar to the definition of $m_{\mathrm{fc}}$, we define $m_{\mathrm{fc}}^{\vartheta
}$ and $\widehat{m}^{\vartheta}_{\mathrm{fc}}$ as the solutions to the equations
%
\begin{eqnarray}
\label{lambdamfc} m^{\vartheta}_{\mathrm{fc}}(z)=\int\frac{\dd\nu(v)}{\vartheta
v-z-m_{\mathrm{fc}}
^{\vartheta}(z)},\qquad
\im m^{\vartheta}_{\mathrm{fc}}(z)\ge0\qquad\bigl(z\in\C^+\bigr)
\end{eqnarray}
and
%
\begin{eqnarray}
\label{hatmfc} \widehat{m}^{\vartheta}_{\mathrm{fc}}(z)=\int
\frac{\dd\widehat\nu
(v)}{\vartheta v-z-\widehat{m}^{\vartheta}_{\mathrm{fc}}(z)},\qquad\im\widehat
m^{\vartheta}_{\mathrm{fc}}(z)\ge0,
\qquad\bigl(z\in\C^+\bigr),
\end{eqnarray}
respectively. Following the discussion of Section~\ref
{subsectiondeforemdsemicirlelaw}, $m^{\vartheta}_{\mathrm{fc}}$ and $\widehat
{m}^{\vartheta}_{\mathrm{fc}}$ define two probability measures $\rho
^{\vartheta}_{\mathrm{fc}}$ and $\widehat\rho^{\vartheta}_{\mathrm{fc}}$ through the densities
%
\begin{eqnarray}
\label{stieltjesinversionformula2} \rho^{\vartheta}_{\mathrm{fc}}(E):=\lim
_{\eta\searrow
0}\frac{1}{\pi
}\im m^{\vartheta}_{\mathrm{fc}}(E+
\ii\eta)\qquad(E\in\R)
\end{eqnarray}
and
%
\begin{eqnarray}
\label{stieltjesinversionformula3} \widehat\rho^{\vartheta
}_{\mathrm{fc}}(E): =
\lim_{\eta
\searrow0}\frac
{1}{\pi}\im\widehat m^{\vartheta}_{\mathrm{fc}}(E+
\ii\eta) \qquad(E\in\R);
\end{eqnarray}
cf.~(\ref{stieltjesinversionformula}). More precisely, we have the
following result which follows directly from the proofs of Lemmas~\ref
{lemmamfc} and~\ref{lemmahatmfc} below. Recall the definition
of $\Theta_{\varpi}$ in~(\ref{definitionTheta}).

%
\begin{lemma}\label{definitionhatmufc}
Let $\widehat\nu$ and $\nu$ satisfy the Assumptions~\ref
{assumptionmuVconvergence} and~\ref{assumptionmuV}. Then, for any
$\vartheta
\in\Theta_{\varpi}$ and $N\in\N$, equations~(\ref{lambdamfc})
and~(\ref{hatmfc}) define, through the inversion formulas in~(\ref
{stieltjesinversionformula2}) and~(\ref{stieltjesinversionformula3}),
absolutely continuous measures $\rho^{\vartheta}_{\mathrm{fc}}$ and
$\widehat\rho^{\vartheta}_{\mathrm{fc}}$. Moreover, the measure~$\rho
_{\mathrm{fc}}^{\vartheta}$ is supported on a single interval with strictly
positive density inside this interval. The same holds true on~$\Omega$
for the measures $\widehat\rho_{\mathrm{fc}}^{\vartheta}$, for $N$
sufficiently large.
\end{lemma}

Note that if $(v_i)$ are random, then so are $\widehat
{m}_{\mathrm{fc}}^{\vartheta}$, respectively, $\widehat{\rho}_{\mathrm{fc}}^{\vartheta
}$. As noted above, we use the symbol~$ \widehat{\phantom{m}}$~to denote
quantities that depend on the empirical distribution $\widehat\nu$ of
the $(v_i)$, while we drop this symbol for quantities depending on the
limiting distribution $\nu$ of $(v_i)$.\vspace*{1pt}

We denote by $\widehat L_{\pm}^{\vartheta}$, respectively, $L_{\pm
}^{\vartheta}$, the endpoints of the support of $\widehat\rho
_{\mathrm{fc}}^{\vartheta}$, respectively, $\rho_{\mathrm{fc}}^{\vartheta}$. Let
$E_0\ge1+\max\{\llvert L_-^1\rrvert,L_+^1\}$, and define the domain
%
\begin{equation}
\label{eqDL} \caD_L: =\bigl\{z=E+\ii\eta\in\C:
\llvert E\rrvert\le E_0, (\varphi_N)^{L}
\le N\eta\le3N\bigr\},
\end{equation}
with $L\equiv L(N)$, such that $L\ge12\xi$; see~(\ref{eqxi}).

The following theorem is the main result of this section.

\begin{theorem}[(Strong local deformed semicircle law)]\label{thmstrong}
Let $H^\vartheta=\vartheta V+W$, $\vartheta\in\Theta_{\varpi}$
[see~(\ref{definitionTheta})], where~$W$ is a real symmetric or
complex Hermitian Wigner matrix satisfying the assumptions in
Definition~\ref{assumptionwigner} and~$V$ is a deterministic or
random real diagonal matrix satisfying Assumptions~\ref
{assumptionmuVconvergence} and~\ref{assumptionmuV}. Let
%
\begin{equation}
\label{eqstrongxi} \xi=\frac{A_0+o(1)}{2}\log\log N.
\end{equation}
Then there are constants $\upsilon>0$ and $c_1$, depending on the
constants $E_0$ in~(\ref{eqDL}), $\alpha_0$ in~(\ref
{definitionofomegaV}), $A_0$, $a_0$, $C_1$ in~(\ref{eqxi}), $\theta$, $C_0$
in~(\ref{eqC0}) and the measure~$\widehat\nu$ such that the
following holds for $L\ge40\xi$. For any $z\in\caD_L$ and any
$\vartheta\in\Theta_{\varpi}$, we have
%
\begin{eqnarray}
\label{strong1} \bigl\llvert m_N^{\vartheta}(z)-
\widehat{m}^{\vartheta}_{\mathrm{fc}}(z) \bigr\rrvert&\le&(
\varphi_N)^{c_1\xi}\frac{1}{N\eta},
\end{eqnarray}
with $(\xi,\upsilon)$-high probability on $\Omega$.

Moreover, we have, for any $z\in\caD_L$, any $\vartheta\in\Theta
_{\varpi}$ and any $i,j\in\llbracket1,N\rrbracket$,
%
\begin{equation}
\label{strong2} \bigl\llvert G_{ij}^{\vartheta}(z)-
\delta_{ij} \widehat g_{i}^{\vartheta
}(z) \bigr\rrvert
\le(\varphi_N)^{c_1\xi} \biggl(\sqrt{
\frac{\im
\widehat m^{\vartheta}_{\mathrm{fc}}(z)}{N\eta}}+\frac{1}{N\eta} \biggr),
\end{equation}
with $(\xi,\upsilon)$-high probability on $\Omega$, where we have set
%
\begin{equation}
\widehat{g}^{\vartheta}_i(z): =
\frac{1}{\vartheta
v_i-z-\widehat
{m}^{\vartheta}_{\mathrm{fc}}(z)}.
\end{equation}
\end{theorem}

The study of local laws for Wigner matrices was initiated in~\cite
{ESY1,ESY2,ESY3}. For more recent results, we refer to~\cite{EKYY4}.
For deformed Wigner matrices with random potential, a local law was
obtained in~\cite{LS}.

Denote by $\bflambda^{\vartheta}=(\lambda_1^{\vartheta},\lambda
_2^{\vartheta},\ldots,\lambda_N^{\vartheta})$ the eigenvalues of
the random matrix $H^{\vartheta}=\vartheta V+W$ arranged in ascending
order. We define the classical location,~$\widehat\gamma
_{i}^{\vartheta}$, of the eigenvalue $\lambda_i^\vartheta$ by
%
\begin{equation}
\label{classicallocation} \int_{-\infty}^{\widehat\gamma_i^{\vartheta
}}\widehat{\rho
}_{\mathrm{fc}}^{\vartheta}(x)\,\dd x=\frac{ i-\sklfrac{1}{2}}{N}\qquad(1\le i\le N).
\end{equation}
Note that $(\widehat\gamma_i^{\vartheta})$ are random in case
$(v_i)$ are too. We have the following rigidity result on the
eigenvalue locations of $H^{\vartheta}$:

\begin{corollary}\label{rigidityofeigenvalues}
Let $H^{\vartheta}=\vartheta V+W$, $\vartheta\in\Theta_{\varpi}$,
where~$W$ is a real symmetric or complex Hermitian Wigner matrix
satisfying the assumptions in Definition~\ref{assumptionwigner},
and~$V$ is a deterministic or random real diagonal matrix satisfying
Assumptions~\ref{assumptionmuVconvergence} and~\ref{assumptionmuV}. Let
$\xi$ satisfy~(\ref{eqstrongxi}). Then there are
constants $\upsilon>0$ and $c_1,c_2$, depending on the constants $E_0$
in~(\ref{eqDL}), $\alpha_0$ in~(\ref{definitionofomegaV}),
$A_0$, $a_0$, $C_1$ in~(\ref{eqxi}), $\theta$, $C_0$ in~(\ref
{eqC0}) and the measure~$\widehat\nu$, such that
%
%
\begin{eqnarray}\label{rigidityofeigenvaluesequation}
\bigl\llvert\lambda_{i}^{\vartheta}-\widehat
\gamma_{i}^{\vartheta}\bigr\rrvert&\le& (\varphi_N)^{c_1\xi}
\frac{1}{N^{2/3}\check\alpha_i^{1/3} }\qquad(1\le i\le N),
\\
\sum_{i=1}^N\bigl\llvert
\lambda_{i}^{\vartheta}-\widehat\gamma_{i}^{\vartheta
}
\bigr\rrvert^2&\le& (\varphi_N)^{c_2\xi}
\frac{1}{N},
\end{eqnarray}
with $(\xi,\upsilon)$-high probability on $\Omega$, for all
$\vartheta\in\Theta_{\varpi}$, where we have abbreviated $\check
\alpha_i: =\min\{i,N-i+1\}$.
\end{corollary}

In the rest of this section we sum up the proofs of Theorem~\ref
{thmstrong} and Corollary~\ref{rigidityofeigenvalues}.

\subsection{Properties of \texorpdfstring{$m^{\vartheta}_{\mathrm{fc}}$}{mvarthetafc} and \texorpdfstring{$\widehat{m}^{\vartheta}_{\mathrm{fc}}$}{widehatmvarthetafc}}
In this subsection, we discuss properties of the Stieltjes transforms
$m_{\mathrm{fc}}^{\vartheta}$ and $\widehat{m}^{\vartheta}_{\mathrm{fc}}$. We first
derive the desired properties for $m_{\mathrm{fc}}^{\vartheta}$ (Lemma~\ref
{lemmamfc} and Corollary~\ref{remarkaboutrealandimaginarypartofmfc} in
the \hyperref[app]{Appendix}) and then show in a second step that
${m}_{\mathrm{fc}}^{\vartheta}$ is a good approximation to $\widehat
m_{\mathrm{fc}}^{\vartheta}$ so that $\widehat{m}_{\mathrm{fc}}^{\vartheta}$ also
shares these properties; see Lemma~\ref{lemmahatmfc}.

For $E_0$ as in~(\ref{definitionofD}), we define the domain, $\caD
'$, of the spectral parameter~$z$ by
%
\begin{eqnarray}
\label{definitionofD} \caD': =\bigl\{z=E+\ii\eta\dvtx E
\in[-E_0,E_0], \eta\in(0,3]\bigr\}.
\end{eqnarray}

The next lemma, whose proof is postponed to the \hyperref[app]{Appendix}, gives a
qualitative description of the deformed semicircle law $\rho
^{\vartheta}_{\mathrm{fc}}$ and its Stieltjes transform $m^{\vartheta}_{\mathrm{fc}}$.

%
\begin{lemma}\label{lemmamfc}
Let\vspace*{1pt} $\nu$ satisfy Assumption~\ref{assumptionmuV}, for some $\varpi
>0$. Then the following holds true for any $\vartheta\in\Theta
_{\varpi}$. There are $L_-^{\vartheta},L_+^{\vartheta}\in\R$, with
$L_-^{\vartheta}<0<L_+^{\vartheta}$, such that $\supp\rho
_{\mathrm{fc}}^{\vartheta}=[L_-^{\vartheta},L_+^{\vartheta}]$, and there
exists a constant $C>1$ such that, for all $\vartheta\in\Theta
_{\varpi}$,
%
\begin{eqnarray}
\label{thesquareroot} C^{-1}\sqrt{\kappa_E}\le
\rho_{\mathrm{fc}}^{\vartheta}(E)\le C \sqrt{\kappa_E} \qquad
\bigl(E\in\bigl[L_-^{\vartheta},L_+^{\vartheta}\bigr]\bigr),
\end{eqnarray}
where $\kappa_E$ denotes the distance of $E$ to the endpoints of the
support of $\rho_{\mathrm{fc}}^{\vartheta}$, that is,
%
\begin{eqnarray}
\kappa_E: =\min\bigl\{\bigl\llvert
E-L_-^{\vartheta}\bigr\rrvert, \bigl\llvert E-L_+^{\vartheta}\bigr
\rrvert
\bigr\}.
\end{eqnarray}
The Stieltjes transform, $m^{\vartheta}_{\mathrm{fc}}$, of $\rho^{\vartheta
}_{\mathrm{fc}}$ has the following properties:
\begin{longlist}[(2)]
\item[(1)] for all $z=E+\ii\eta\in\caD'$,
%
\begin{eqnarray}
\label{behaviorofmfc} \im m^{\vartheta}_{\mathrm{fc}}(z)\sim\cases{ \sqrt{\kappa+
\eta}, &\quad$\displaystyle E\in\bigl[L_-^{\vartheta},L_+^{\vartheta}
\bigr]$,
\vspace*{3pt}\cr
\displaystyle\frac{\eta}{\sqrt{\kappa+\eta}}, &\quad $E\in
\bigl[L_-^{\vartheta},L_+^{\vartheta}\bigr]^c$;}
\end{eqnarray}

\item[(2)] there exists a constant $C>1$ such that for all $z\in
\caD'$ and all $x\in I_{\nu}$,
%
\begin{eqnarray}
\label{stabilitybound} C^{-1}\le\bigl\llvert\vartheta x -z-m^{\vartheta
}_{\mathrm{fc}}(z)
\bigr\rrvert\le C.
\end{eqnarray}
\end{longlist}
Moreover, the constants in~(\ref{thesquareroot}),~(\ref{behaviorofmfc})
and~(\ref{stabilitybound}) can be chosen uniformly in
$\vartheta\in\Theta_{\varpi}$.
\end{lemma}

Next, we argue that $\widehat m_{\mathrm{fc}}^{\vartheta}$ behaves
qualitatively in the same way as~$m_{\mathrm{fc}}^{\vartheta}$ on~$\Omega$
for~$N$ sufficiently large. Lemma~\ref{lemmahatmfc} below is proven
in the \hyperref[app]{Appendix}.

%
\begin{lemma}\label{lemmahatmfc}
Let $\widehat\nu$ satisfy Assumptions~\ref{assumptionmuVconvergence}
and~\ref{assumptionmuV}, for some $\varpi>0$. Then the
following holds\vspace*{1pt} for all $\vartheta\in\Theta_{\varpi}$ and all
sufficiently large $N$ on~$\Omega$. There are $\widehat L_-^{\vartheta
}, \widehat L_+^{\vartheta}\in\R$, with $\widehat L_-^{\vartheta
}<0<\widehat L_+^{\vartheta}$, such\vspace*{2pt} that $\supp\widehat\rho
_{\mathrm{fc}}^{\vartheta}=[\widehat L_-^{\vartheta},\widehat L_+^\vartheta
]$. Let $\widehat\kappa_E: =\min\{\llvert E-\widehat
{L}_-^{\vartheta}\rrvert,
\llvert E-\widehat{L}_+^{\vartheta}\rrvert \}$. Then~(\ref
{thesquareroot}),~(\ref{behaviorofmfc}) and~(\ref{stabilitybound}) of
Lemma~\ref{lemmamfc}, hold true on $\Omega$, for $N$ sufficiently
large, with~$m_{\mathrm{fc}}^{\vartheta}$ replaced by~$\widehat{m}^{\vartheta
}_{\mathrm{fc}}$,~$\rho_{\mathrm{fc}}^\vartheta$ replaced by~$\widehat\rho
_{\mathrm{fc}}^{\vartheta}$, etc. Moreover, the constants in these inequalities
can be chosen uniformly in $\vartheta\in\Theta_{\varpi}$ and $N$,
for $N$
sufficiently large.

Further, there is $c>0$ such that for all $z\in\caD'$ we have
%
\begin{eqnarray}
\label{minilocallaw} \bigl\llvert\widehat{m}^{\vartheta
}_{\mathrm{fc}}(z)-m_{\mathrm{fc}}^{\vartheta}(z)
\bigr\rrvert\le N^{-c\alpha_0/2} \qquad\bigl\llvert\widehat L_\pm
^{\vartheta}-{L}_\pm
^{\vartheta}\bigr\rrvert\le N^{-c\alpha_0},
\end{eqnarray}
on $\Omega$ for $N$ sufficiently large and all $\vartheta\in\Theta
_{\varpi}$.
\end{lemma}

\subsection{Proof of Theorem~\texorpdfstring{\protect\ref{thmstrong}}{3.3} and Corollary~\texorpdfstring{\protect\ref{rigidityofeigenvalues}}{3.4}}

The proof of Theorem~\ref{thmstrong} follows closely the proof of
Theorem~2.10 in~\cite{LS}. The difference between Theorem~\ref
{thmstrong} of the present paper and Theorem~2.10 in~\cite{LS} is
that we presently condition on the diagonal entries $(v_i)$; that is,
we consider the entries of~$V$ as fixed. Accordingly,\vspace*{1pt} we compare [on
the event $\Omega$ of typical $(v_i)$] the averaged Green function
$m^{\vartheta}$ with $\widehat m^{\vartheta}_{\mathrm{fc}}$ [see~(\ref{hatmfc})]
instead of $m_{\mathrm{fc}}^{\vartheta}$; see~(\ref{lambdamfc}). For
consistency, we momentarily drop the $\vartheta$ dependence form our
notation. To establish Theorem~\ref{thmstrong}, we first derive a
weak local deformed semicircle law (see Theorem~4.1 in~\cite{LS}) by
following the proof in~\cite{LS}. Using the Lemma~\ref{lemmamfc},
Lemma~\ref{lemmahatmfc} and the results in the \hyperref[app]{Appendix}, it is then
straightforward to obtain the following result.

\begin{lemma}\label{weaklocallaw}
Under the assumption of Theorem~\ref{thmstrong}, there are $c_1$ and
$\upsilon>0$ such that
\begin{eqnarray*}
\bigl\llvert m_N(z)-\widehat m_{\mathrm{fc}}(z)\bigr\rrvert\le(
\varphi_N)^{c_1\xi}\frac{1}{(N\eta
)^{1/3}},\qquad\bigl\llvert
G_{ij}(z)\bigr\rrvert\le(\varphi_N)^{c_1\xi}
\frac{1}{\sqrt
{N\eta}},
\end{eqnarray*}
with $(\xi,\upsilon)$-high probability on $\Omega$, uniformly in
$z\in\caD_L$ and $\vartheta\in\Theta_{\varpi}$.
\end{lemma}

To prove~Theorem~\ref{thmstrong} we follow mutatis mutandis the proof
of Theorem~4.1 in~\cite{LS}. But we note that in the corresponding
equation to~(5.25) in~\cite{LS}, we may set $\lambda=0$ in the error
term, at the cost of replacing~$m_{\mathrm{fc}}$ by~$\widehat m_{\mathrm{fc}}$. In the
subsequent analysis, we can simply set $\lambda=0$ in the error terms.
In this way, one establishes the proof of Theorem~\ref{thmstrong}.
Similarly, Corollary~\ref{rigidityofeigenvalues} can be proven in
the same way as is Theorem~2.21 in~\cite{LS}. It suffices to set
$\lambda=0$ in the analysis in~\cite{LS}. We leave the details aside.

\section{Reference \texorpdfstring{$\beta$}{beta}-ensemble}\label{sectionbetaensemble}

\subsection{Definition of \texorpdfstring{$\beta$}{beta}-ensemble and known results}\label{subsectionbetaensembleintro}
We first recall the notion of $\beta$-ensembles. Let $N\in\N$, and
let ${\digamma}^{(N)}\subset\R^N$ denote the set
%
\begin{eqnarray}
\label{Weylchamber} \digamma^{(N )}: =\bigl\{ \bfx=
(x_1,x_2,\ldots, x_N )\dvtx x_1\le
x_2\le\cdots\le x_N\bigr\}.
\end{eqnarray}
Consider the probability distribution, $\mu_U\equiv\mu_U^N$, on
$\digamma^{(N)}$ given by
%
\begin{eqnarray}
\label{eqnmeasure} \qquad\mu_{U}^N (\dd\bfx): =
\frac{1}{Z_{U}^{N}}\mathrm{e}^{-\beta N
\caH(\bfx) }\,\dd\bfx,\qquad\dd\bfx: =\lone\bigl(\bfx\in\digamma^{(N)}\bigr)\,\dd x_1\,\dd
x_2\cdots\,\dd x_N,
\end{eqnarray}
where $\beta>0$,
%
\begin{eqnarray}
\label{lehamiltonianwithU} \caH(\bfx) : =\sum_{i=1}^N
\frac{1}{2} \biggl(U(x_i)+\frac
{x_i^2}{2} \biggr)-
\frac{1}{N}\sum_{1\le i<j\le N}\log(x_j-x_i)
\end{eqnarray}
and $Z_{U}^{N}\equiv Z_{U}^N(\beta)$ is a normalization. Here $U$ is a
potential, that is, a real-valued, sufficiently regular function on $\R
$. In the following, we often omit the parameters $N$ and $\beta$ from
the notation. We use $\mathbb{P}^{\mu_U}$ and $\E^{\mu_U}$ to
denote the probability and the expectation with respect to $\mu_U$.
We view $\mu_U$ as a Gibbs measure of~$N$ particles on~$\R$ with a
logarithmic interaction, where the parameter
$\beta> 0$ may be interpreted as the inverse temperature. (For the
results in the present paper, we choose $\beta=2$ in case~$W$ is
complex Hermitian Wigner matrix and $\beta=1$ in case~$W$ is a real
symmetric Wigner matrix.) We refer to the variables $(x_i)$ as
particles or points, and we call the system a log-gas or a $\beta
$-ensemble. We assume that the potential $U$ is a $C^4$
function on $\R$ such that its second derivative is bounded below;
that is, we have
%
\begin{eqnarray}
\label{assumption1forbetauniversality} \inf_{x\in\R} U'' (x)
\ge-2C_U,
\end{eqnarray}
for some constant $C_U\ge0$, and we further assume that
%
\begin{eqnarray}
\label{assumption2forbetauniversality} U (x)+\frac{x^2}{2} > (2 +
\varepsilon) \log\bigl(1 + \llvert x
\rrvert\bigr)\qquad(x\in\R),
\end{eqnarray}
for some $\varepsilon> 0$, for large enough $\llvert x\rrvert $. It
is well known (see,
e.g.,~\cite{BPS}) that under these conditions
the measure is normalizable, $Z_{U}^{N} < \infty$. Moreover, the
averaged density of the empirical spectral measure, $\rho_U^N$,
defined as
%
\begin{eqnarray}
\label{referencesection71} \rho_{U}^{N}: =
\E^{\mu_{U}}\frac{1}{N}\sum_{i=1}^N
\delta_{x_i},
\end{eqnarray}
converges weakly in the limit $N\to\infty$ to a continuous function
$\rho_U$, the equilibrium density, which is of compact support. It is
well known that~$\rho_U$ can be obtained as the unique solution to the
variational problem
%
\begin{eqnarray}
\label{minimizationproblem}
&&\inf\biggl\{\int_{\R} \biggl(
\frac{x^2}{2}+U(x) \biggr)\,\dd\rho(x)-\int_{\R}\log
\llvert x-y\rrvert\,\dd\rho(x)\,\dd\rho(y)\dvtx
\nonumber\\[-8pt]\\[-8pt]\nonumber
&&\hspace*{141pt}\rho\mbox{ is a probability measure} \biggr\}
\end{eqnarray}
and that the equilibrium density $\rho=\rho_U$ satisfies
%
\begin{eqnarray}
\label{regularpotentialU} U'(x)+x=-2\dashint_\R\frac{\rho(y)\,\dd y}{y-x}
\qquad(x\in\supp\rho_U).
\end{eqnarray}
In fact,~(\ref{regularpotentialU}) holds if and only if $x\in\supp
\rho_U$. We will assume in addition that the minimizer $\rho_U$ is
supported on a single interval $[A_-,A_+]$ and that $U$ is ``regular''
in the sense of~\cite{KML}; that is, the equilibrium density of $U$ is
positive on $(A_-,A_+)$ and vanishes
like a square root at each of the endpoints of $[A_-, A_+]$. Viewing
the points $\bfx=(x_i)$ as points or particles on $\R$, we define the
classical location of the $k$th particle, $\gamma_k$, under the $\beta
$-ensemble $\mu_U$ by
%
\begin{equation}
\label{generalbetaensembleclassicallocation} \int_{-\infty}^{\gamma
_k}\rho_U(x)\,\dd x=\frac{k-\sklfrac{1}{2}}{N}.
\end{equation}
For a detailed discussion of general $\beta$-ensemble we refer, for
example, to~\cite{AGZ,BEY}.

For $U\equiv0$, we write $\mu_G\equiv\mu_G^N$ instead of $\mu_0$,
since $\mu_0$ is the equilibrium measure for the GUE ($\beta=2$),
respectively, the GOE ($\beta=1$). More precisely, setting
%
\begin{equation}
\label{gaussianhamiltonian} \caH_G(\bfx) : =\sum
_{i=1}^N\frac
{1}{4}x_i^2-
\frac{1}{N}\sum_{1\le i<j\le N}\log(x_j-x_i),
\end{equation}
the GUE, respectively, GOE, distribution on $\digamma^{(N)}$ are given by
%
\begin{equation}
\label{gaussianmeasures} \mu_{G}^N(\dd\bfx)=\frac{1}{Z_{G}^N}
\mathrm{e}^{-\beta N \caH
_G(\bfx)}\,\dd\bfx,
\end{equation}
where $Z_{G}^N\equiv Z_G^N(\beta)$ is a normalization, and we either
choose $\beta=2$ or $\beta=1$.

We are interested in the $n$-point correlation functions defined by
%
\begin{eqnarray}
\label{eqncorrFunct} \varrho^{N}_{U,n}(x_1,
\dots,x_n)= \int_{\R^{N-n}} \mu_U^\#(
\bfx)\,\dd x_{n+1}\cdots\,\dd x_{N},
\end{eqnarray}
where $ \mu_U^\#$ is the symmetrized version of $\mu_U$ given
in~(\ref{eqnmeasure})
but defined on $\R^N$ instead of the simplex~$\digamma^{(N)}$,
%
\begin{eqnarray}
\label{symmetrization} \mu_U^{\# }(\dd x)=\frac{1}{N!}
\mu_U\bigl(\dd\bfx^{(\sigma)}\bigr),\qquad\dd x= \dd
x_1\cdots\,\dd x_N,
\end{eqnarray}
where
$\bfx^{(\sigma)}=(x_{\sigma(1)},\dots,x_{\sigma(N)})$, with
$x_{\sigma(1)}<\cdots<x_{\sigma(N)}$. The following universality
result is proven in~\cite{BEY}.

%
\begin{theorem}[(Bulk universality for $\beta$-ensembles, Theorem~2.1
in~\cite{BEY})]\label{thmbulkuniversalitybeta}
Let $U$ be a ${C}^4$ regular potential with equilibrium density
supported on a single interval $[A_-,A_+]$ that satisfies~(\ref
{assumption1forbetauniversality}) and~(\ref
{assumption2forbetauniversality}). Then the following result holds.
For any fixed $\beta>0$, $E\in(A_-,A_+)$, $\llvert E'\rrvert <2$, $
n\in\N$,
$0<\delta\le\frac{1}{2}$ and any $n$-particle observable $O$,
we have with $b: =N^{-1+\delta}$,
\begin{eqnarray*}
&& \lim_{N\to\infty}\int_{\R^n} \dd
\alpha_1 \cdots\,\dd\alpha_n O(\alpha_1, \dots,
\alpha_n)
\nonumber
\\
&&\hspace*{41pt}{}\times\biggl[ \int_{E - b}^{E + b} \frac{\dd x}{2 b}
\frac{1}{[ \rho_U (E)]^n } \varrho_{U,n}^{N} \biggl( x +
\frac{\alpha_1}{N\rho_U(E)}, \dots, x + \frac{\alpha_n}{N\rho
_U(E)} \biggr)
\\
&&\hspace*{71pt}{}- \frac{1}{[\rho_{\mathrm{sc}}(E')]^n} \varrho_{{G}, n}^{N} \biggl(
E' + \frac{\alpha_1}{N\rho_{\mathrm{sc}}(E')}, \dots, E' +
\frac{\alpha
_n}{N\rho_{\mathrm{sc}}(E')} \biggr) \biggr]=0.
\end{eqnarray*}
Here, $\rho_{\mathrm{sc}}$ denotes the density of the semicircle law,
and $\varrho_{{G}, n}^{N}$ is the $n$-point the correlation function
of the Gaussian $\beta$-ensemble, that is,
with $U\equiv0$.
\end{theorem}

Theorem~\ref{thmbulkuniversalitybeta} was first proved in~\cite
{BEYI} under the assumption that $U$ is analytic,
a hypothesis that was only required for proving rigidity. The
analyticity assumption has been removed in~\cite{BEY}. Recently,
alternative proofs of bulk universality for \mbox{$\beta$-}ensembles with
general $\beta>0$, that is, results similar to Theorem~\ref
{thmbulkuniversalitybeta}, have been obtained in~\cite{ShMa} and~\cite{BAG}.
In the present paper, we will not use Theorem~\ref
{thmbulkuniversalitybeta}; it is stated here for completeness.

To conclude this subsection, we recall an important tool in the study
of \mbox{$\beta$-}ensembles, the ``first order loop'' equation. In the
notation above it reads (in the limit $N\to\infty$)
%
\begin{eqnarray}
\label{loopequation} m_U(z)^2=\int\frac{x+ U'(x)}{x-z}
\rho_U(x)\,\dd x \qquad\bigl(z\in\C^+\bigr),
\end{eqnarray}
where $m_U$ denotes the Stieltjes transform of the equilibrium measure
$\rho_U$, that is,
\begin{eqnarray*}
m_U(z)\equiv m_{\rho_U}(z)=\int\frac{\rho_U(x)}{x-z} \,\dd x
\qquad\bigl(z\in\C^+\bigr).
\end{eqnarray*}
The loop equation~(\ref{loopequation}) can be obtained by a change of
variables in~(\ref{eqnmeasure}) (see~\cite{J3}) or
by integration by parts; see~\cite{ShMa1}.

\subsection{Time-dependent modified \texorpdfstring{$\beta$}{beta}-ensemble}
In this subsection, we introduce a modified $\beta$-ensemble by
specifying potentials $\widehat U$ and $ U$ that depend, among other
things, on a parameter $t\ge0$ which has the interpretation of a time.
The potential $\widehat U$ also depends on $N$, the size of our
original matrix $H=V+W$, yet the $N$ dependence is only through the
fixed random variables $(v_i)$. Recall that we have defined $\widehat
m_{\mathrm{fc}}^{\vartheta}$, respectively, $ m_{\mathrm{fc}}^{\vartheta}$, as the
solutions to the equations
%
\begin{eqnarray}
\label{nochmalsmfc} \qquad\widehat{m}^{\vartheta}_{\mathrm{fc}}(z)=\int
\frac{\dd\widehat\nu
(v)}{\vartheta v_i-z-\widehat{m}^{\vartheta}_{\mathrm{fc}}(z)},\qquad
m^{\vartheta}_{\mathrm{fc}}(z)=\int
\frac{\dd\nu(v)}{\vartheta v-z-m_{\mathrm{fc}}
^{\vartheta}(z)},
\end{eqnarray}
$z\in\C^+$, subject to the conditions $\im\widehat m^{\vartheta
}_{\mathrm{fc}}(z), \im m^{\vartheta}_{\mathrm{fc}}(z)\ge0$, for $\im z>0$. Recall
from~(\ref{definitionTheta}) that we denote $\Theta_{\varpi
}=[0,1+\varpi']$, $\varpi'=\varpi/10$. We then fix some $ t_0\ge0$
such that $\mathrm{e}^{t_0/2}\in\Theta_{\varpi}$ and let
%
\begin{eqnarray}
\label{equationforvartheta} \vartheta\equiv\vartheta(t):=\mathrm
{e}^{-(t-t_0)/2}\qquad(t\ge0).
\end{eqnarray}
In the following we consider $t\ge0$ as time, and we henceforth
abbreviate $m^{\vartheta(t)}_{\mathrm{fc}}(z)\equiv m_{\mathrm{fc}}(t,z)$, etc.
Equation~(\ref{nochmalsmfc}) defines time dependent measures
$\widehat\rho_{\mathrm{fc}}(t)$, $\rho_{\mathrm{fc}}(t)$, respectively, whose densities
at the point $x\in\R$ are denoted by $\widehat\rho_{\mathrm{fc}}(t,x)$,
respectively, $\rho_{\mathrm{fc}}(t,x)$.

We denote by $\widehat U'(t,x)$, $\widehat U^{(n)}(t,x)$ the first,
respectively, the $n$th derivative of $\widehat U(t,x)$ with respect to
$x$, and we use the same notation for $U$. We define~$\widehat U$
and~$U$ (up to finite additive constants that enter the formalism only
in normalizations) through their derivatives $\widehat U'$ and $U'$.
For $t\ge0$, we set
%
\begin{eqnarray}
\label{definitionwidehatU} \widehat U'(t,x)+x: =-2
\dashint_\R\frac{\widehat
\rho
_{\mathrm{fc}}(t,y)}{y-x} \,\dd y,
\end{eqnarray}
for $x\in\supp\widehat\rho_{\mathrm{fc}}(t)$, respectively,
%
\begin{eqnarray}
\label{definitionwidehatU2} U'(t,x)+x: =-2\dashint_\R
\frac{\rho
_{\mathrm{fc}}(t,y)}{y-x} \,\dd y,
\end{eqnarray}
for $x\in\supp\rho_{\mathrm{fc}}(t)$. Outside the support of the measures
$\widehat\rho_{\mathrm{fc}}(t)$ and $\rho_{\mathrm{fc}}(t)$, we define $\widehat U'$
and $U'$ as $C^3$ extensions such that they are ``regular'' potentials
satisfying~(\ref{assumption1forbetauniversality}) and~(\ref
{assumption2forbetauniversality}) for all $t\ge0$. The definitions
of such potentials are obviously not unique. One possible construction
is outlined in the \hyperref[app]{Appendix} in the form of the proof of the next lemma.

%
\begin{lemma}\label{superlemma}
There exist potentials $\widehat U, U\dvtx \R^+\times\R\to\R$,
$(t,x)\mapsto\widehat U(t,x)$, $ U(t,x)$ such that for $n\in\llbracket
1,4\rrbracket$, $\widehat U^{(n)}(t,x)$, $U^{(n)}(t,x)$, $\partial_t
\widehat U^{(n)}(t,x)$, $\partial_t U^{(n)}(t,x)$ are continuous
functions of $x\in\R$ and $t\in\R^+$, which can be uniformly
bounded in~$ x$ on compact sets, uniformly in $t\in\R^+$ and
sufficiently large~$N$. Moreover the following holds for all $t\ge0$
on $\Omega$ for $N$ sufficiently large:
\begin{longlist}[(2)]
\item[(1)] $\widehat U'(t,x)$ and $ U'(t,x)$ satisfy~(\ref
{definitionwidehatU}) and~(\ref{definitionwidehatU2}) for $x\in
\supp\widehat\rho_{\mathrm{fc}}(t)$, respectively, $x\in\supp\rho_{\mathrm{fc}}(t)$.
For $x\notin\supp\widehat\rho_{\mathrm{fc}}(t)$, respectively, $x\notin\supp
\rho_{\mathrm{fc}}(t)$, we have
\begin{eqnarray*}
\bigl\llvert\widehat U'(t,x)+x\bigr\rrvert> 2\bigl\llvert\re
\widehat m_{\mathrm{fc}}(t,x)\bigr\rrvert,\qquad\bigl\llvert
U'(t,x)+x\bigr\rrvert> 2\bigl\llvert\re m_{\mathrm{fc}}(t,x)
\bigr\rrvert.
\end{eqnarray*}
\item[(2)] There is a constant $c>0$ such that for all $x\in\R$
and all $t\ge0$, we have
%
\begin{eqnarray}
\bigl\llvert\widehat U'(t,x)-U'(t,x)\bigr\rrvert
\le N^{-c\alpha_0/2},
\end{eqnarray}
where $\alpha_0>0$ is the constant in~(\ref{definitionofomegaV}).
\item[(3)] The potentials $\widehat U$ and $U$ satisfy~(\ref
{assumption1forbetauniversality}) and~(\ref
{assumption2forbetauniversality}). In particular, there is \mbox{$C_U\ge0$}
(independent of $N$),
such that
%
\begin{eqnarray}
\label{uniformconvexitybound} \inf_{x\in\R, t\in\R^+}\widehat U''(t,x)
\ge-2C_U, \qquad\inf_{x\in\R, t\in\R^+} U''(t,x)
\ge-2C_U.
\end{eqnarray}
Moreover, $\widehat U$ and $U$ are ``regular''; see the paragraph
below~(\ref{regularpotentialU}) for the definition of ``regular'' potential.
\end{longlist}
\end{lemma}

Below, we are mainly interested in $\beta$-ensembles determined by the
potential $\widehat U$. For ease of notation, we thus limit the
discussion to $\widehat U$.

For $N\in\N$ we define a measure on $\digamma^{(N)}$ by setting
%
\begin{eqnarray}
\label{definitionofbt} \qquad{\widehat\psi_t}(\bfx) \mu_G(\dd\bfx)
: =\frac
{1}{Z_{{\widehat
\psi_t}}}\mathrm{e}^{-\sklvfrac{\beta N}{2}{\sum_{i=1}^N\widehat
U(t,x_i)}}\mu_G(
\dd\bfx) \qquad\bigl(\bfx\in\digamma^{(N)}\bigr),
\end{eqnarray}
where $Z_{{\widehat\psi_t}}\equiv Z_{{\widehat\psi_t}}(\beta)$ is
a normalization, and we usually choose $\beta=1,2$. By Lemma~\ref
{superlemma}, ${\widehat\psi_t} \mu_G$ is a well-defined $\beta
$-ensemble, and from the discussion in Section~\ref
{subsectionbetaensembleintro} we further infer that the equilibrium
density of
$\widehat\psi_t \mu_G$, that is, the unique measure solving the
minimization problem in~(\ref{minimizationproblem}), is for any $t\ge
0$, $\widehat\rho_{\mathrm{fc}}(t)$. Viewing~${\widehat\psi_t} \mu_G$ as a
Gibbs measure of~$N$ (ordered) particles~$(x_i)$ on the real line, we
define the classical location of the $i$th particles,~$\widehat\gamma
_i(t)$, as in~(\ref{generalbetaensembleclassicallocation}), that is,
%
\begin{eqnarray}
\label{nochmalsclassicallocations} \int_{-\infty}^{\widehat\gamma
_i(t)}\widehat
\rho_{\mathrm{fc}}(t,x)\,\dd x=\frac{i-\sklfrac{1}{2}}{N} \qquad\bigl(i\in\llbracket1,N
\rrbracket\bigr).
\end{eqnarray}
From~\cite{BEY} we have the following rigidity result.

\begin{proposition}\label{rigidityfortimedependentbetaone}
Let $\widehat U(t,\cdot)$, with $t\ge0$ and $N\in\N$, be given by
Lem\-ma~\ref{superlemma}. Then the following holds on $\Omega$. For any
$\delta>0$, there is $\varsigma>0$, such that for any $t\ge0$,
%
\begin{eqnarray}
\label{rigidityforbensembleone}
\qquad\mathbb{P}^{{\widehat\psi_t} \mu_G}
\bigl(\bigl\llvert x_i-
\widehat\gamma_i(t)\bigr\rrvert> N^{-\sklfrac{2}{3}+ \delta}\check
\alpha_i^{-\sfrac{1}{3}} \bigr)\leq\mathrm{e}^{- N^\varsigma} \qquad(1
\le i\le N),
\end{eqnarray}
for $N$ sufficiently large, where $\mathbb{P}^{{\widehat\psi_t} \mu
_G}$ stands for the probability under ${\widehat\psi_t} \mu_{G}$
conditioned on~$V$. Here, $\check\alpha_i: =\min\{
i,N-i+1\}$.
\end{proposition}

\begin{pf}
The rigidity estimate (\ref{rigidityforbensembleone}) is taken
from Theorem~2.4 of~\cite{BEY}. To achieve uniformity in $t\ge0$ and
$N$ sufficiently large, we note that estimate~(\ref
{rigidityforbensembleone}) depends on the potential mainly through the convexity
bounds\break (\ref{assumption1forbetauniversality}) and~(\ref
{assumption2forbetauniversality}). Starting from the uniform bounds
of Lemma~\ref{superlemma}, one checks that Proposition~\ref
{rigidityfortimedependentbetaone} holds
uniformly in $t$ and $N$ large enough.
\end{pf}

In the rest of this section, we derive equations of motion for the
potential $\widehat U(t,\cdot)$ and the classical locations $(\widehat
\gamma_i(t))$. To derive these equations we observe that the Stieltjes
transform $\widehat m_{\mathrm{fc}}(t,z)$ can be obtained from $\widehat
m_{\mathrm{fc}}(t=0,z)$ as the solution to the following complex Burgers
equation~\cite{P}:
%
\begin{eqnarray}
\label{burgerequation} \partial_t \widehat m_{\mathrm{fc}}(t,z)=
\tfrac{1}{2}\partial_z \bigl[ \widehat m_{\mathrm{fc}}(t,z)
\bigl(\widehat m_{\mathrm{fc}}(t,z)+z\bigr) \bigr] \qquad\bigl(z\in\C^+,t\ge0
\bigr).
\end{eqnarray}
This can be checked by differentiating~(\ref{nochmalsmfc}). Combining
the complex Burgers equation~(\ref{burgerequation}) and the loop
equation~(\ref{loopequation}) we obtain the following result.

%
\begin{lemma}
Let $N\in\N$. Assume that $\widehat\nu$ satisfies the
Assumptions~\ref{assumptionmuVconvergence} and~\ref{assumptionmuV}.
Then the following holds on $\Omega$ for~$N$ sufficiently
large. For $t\ge0$, we have
%
\begin{eqnarray}
\partial_t\widehat\gamma_i(t)=\tfrac{1}{2}
\widehat U'\bigl(t,\widehat\gamma_i(t)\bigr),\label{EoMclassicallocations}
\end{eqnarray}
respectively,
%
\begin{eqnarray}
\label{EoMzwei} \partial_t\widehat\gamma_i(t)=-
\dashint_\R\frac{\widehat\rho
_{\mathrm{fc}}(t,y)}{y-\widehat\gamma_i(t)} \,\dd y-\frac{1}{2}\widehat\gamma
_i(t)\qquad\bigl(i\in\llbracket1,N\rrbracket\bigr).
\end{eqnarray}
Further, the potential $\widehat U$ satisfies
%
\begin{eqnarray}
\label{EoMdrei} \partial_t \widehat U(t,x)&=&\dashint_\R
\frac{\widehat
U'(t,y)\widehat\rho_{\mathrm{fc}}(t,y)}{y-x} \,\dd y \qquad\bigl(x\in\supp\widehat
\rho_{\mathrm{fc}}(t)
\bigr).
\end{eqnarray}
Moreover, there exist constants $C,C'$ such that the following bounds
hold on~$\Omega$:
%
\begin{eqnarray}
\label{boundsEoM} \bigl\llvert\partial_{t}\widehat\gamma_i(t)
\bigr\rrvert\le C,\qquad\bigl\llvert\partial_t \widehat U(t,x)
\bigr\rrvert\le C',
\end{eqnarray}
for all $i\in\llbracket1,N\rrbracket$, uniformly in $t\ge0$, $x\in
\supp\widehat\rho_{\mathrm{fc}}(t)$ and $N$, for $N$ sufficiently large.

Finally, $U(t,\cdot)$ and $(\gamma_i(t))$, share the same properties.
\end{lemma}

\begin{pf}
Combining~(\ref{burgerequation}) and~(\ref{loopequation}), we find,
for $z\in\C^+$, $t\ge0$,
\begin{eqnarray*}
\partial_t\widehat m_{\mathrm{fc}}(t,z)&=&\frac{1}{2}
\partial_z \biggl(-\int\frac{v+\widehat U'(t,v)}{v-z}\widehat
\rho_{\mathrm{fc}}(t,v)\,\dd v+z\int\frac{\widehat\rho_{\mathrm{fc}}(t,v)}{v-z} \,\dd v
\biggr)
\\
&=&\frac{1}{2}\partial_z \biggl(-\int\frac{\widehat
U'(t,v)}{v-z}
\widehat\rho_{\mathrm{fc}}(t,v)\,\dd v-1 \biggr)
\\
&=&-\frac{1}{2}\partial_z\int\frac{\widehat U'(t,v)}{v-z}\widehat
\rho_{\mathrm{fc}}(t,v)\,\dd v.
\end{eqnarray*}
Hence, for $\im z>0$, we get
\begin{eqnarray*}
\partial_t\widehat m_{\mathrm{fc}}(t,z)&=&-\frac{1}{2}\int
\frac{\widehat
U'(t,v)}{(v-z)^2}\widehat\rho_{\mathrm{fc}}(t,v)\,\dd v=-\frac{1}{2}\int
\frac{ (\widehat U'(t,v)\widehat\rho_{\mathrm{fc}}(t,v)
)'}{(v-z)} \,\dd v.
\end{eqnarray*}
Clearly $\widehat U'(t,v)\widehat\rho_{\mathrm{fc}}(t,v)$ is a $C^3$ function
inside the support of $\widehat\rho_{\mathrm{fc}}(z)$ that has a square root
behavior at the endpoints. Thus we obtain from the Stieltjes inversion
formula that
%
\begin{eqnarray}
\label{tapir} \partial_t\widehat\rho_{\mathrm{fc}}(t,E)&=&
\frac{1}{\pi}\lim_{\eta
\searrow0}\im\partial_t
\widehat m_{\mathrm{fc}}(t,z)=-\frac{1}{2} \bigl(\widehat
U'(t,E)\widehat\rho_{\mathrm{fc}}(t,E) \bigr)',
\end{eqnarray}
for all $E\in(\widehat L_-(t),\widehat L_+(t))$, where $\widehat
L_{\pm}(t)$ denote the endpoints of the support of $\widehat\rho_{\mathrm{fc}}(t)$.

On the other hand, differentiating~(\ref{nochmalsclassicallocations})
with respect to time, we obtain
\begin{eqnarray*}
\int_{-\infty}^{\widehat\gamma_i(t)}\partial_t\widehat
\rho_{\mathrm{fc}}(t,v)\,\dd v=-\widehat\rho_{\mathrm{fc}}\bigl(t,\widehat
\gamma_i(t)\bigr)\partial_t\widehat
\gamma_i(t).
\end{eqnarray*}
Substituting from~(\ref{tapir}), we get
\begin{eqnarray*}
\partial_t\widehat\gamma_i(t)=\frac{1}{2}
\frac{1}{\widehat\rho
_{\mathrm{fc}}(t,\widehat\gamma_i(t))}\int_{-\infty}^{\widehat\gamma
_i(t)} \,\dd v \bigl(
\widehat{U}'_{\mathrm{fc}}(t,v)\widehat\rho_{\mathrm{fc}}(t,v)
\bigr)'.
\end{eqnarray*}
Hence
\begin{eqnarray*}
\partial_t\widehat\gamma_i(t)&=&\frac{1}{2}
\frac{1}{\widehat\rho
_{\mathrm{fc}}(t,\widehat\gamma_i(t))}\widehat U' \bigl(t,\widehat\gamma
_i(t) \bigr) \widehat\rho_{\mathrm{fc}} \bigl(t,\widehat
\gamma_i(t) \bigr),
\end{eqnarray*}
and~(\ref{EoMclassicallocations}) follows. Using that $\widehat U$
satisfies~(\ref{definitionwidehatU}), we can recast this last
equation as
\begin{eqnarray*}
\partial_t\widehat\gamma_i(t)=-\dashint_\R
\frac{\widehat\rho
_{\mathrm{fc}}(t,y)}{y-\widehat\gamma_i(t)} \,\dd y-\frac{1}{2}\widehat\gamma_i(t),
\end{eqnarray*}
and we find~(\ref{EoMzwei}). Equation~(\ref{EoMzwei}) follows in a
similar way by differentiating~(\ref{definitionwidehatU}) with
respect to time. By a similar computation we obtain~(\ref{EoMdrei}).
The bound in~(\ref{boundsEoM}) follows from Lemma~\ref{superlemma}.
\end{pf}
Starting from the relations in~(\ref{nochmalsmfc}), we derived via
the time dependent potential~$\widehat U$, an equation of motions for
the classical locations $(\widehat\gamma_i(t))$. The
points~$(\widehat\gamma_i(t))$ may also be viewed as the classical
locations of the eigenvalues of a family of random matrices which is
parametrized by the times $t_0$ and $t$. This is the subject of the
next section.\vadjust{\goodbreak}

\section{Dyson Brownian motion: Evolution of the entropy}\label{sectionDBM}

\subsection{Dyson Brownian motion}\label{subsectionledysonbrownianmotion}

Let $H_{0}=(h_{ij,0})$ be the matrix
\begin{eqnarray*}
H_{0}: =\mathrm{e}^{t_0/2} V+W,
\end{eqnarray*}
where~$V$ satisfies Assumptions~\ref{assumptionmuVconvergence}
and~\ref{assumptionmuV}, and~$W$ is real symmetric or complex
Hermitian satisfying the assumptions in Definition~\ref
{assumptionwigner}. Here, $t_0\ge0$ is chosen such that $\vartheta
=\mathrm
{e}^{t_0/2}\in\Theta_{\varpi}$ [see~(\ref{definitionTheta})], and
we consider $\vartheta$ as an a priori free ``coupling parameter''
that we fix in Section~\ref{Proofsofmainresults} below. Let
$B=(b_{ij})\equiv(b_{ij,t})$ be a real symmetric, respectively, a
complex Hermitian, matrix whose entries are a collection of
independent, up to the symmetry constraint, real (complex) Brownian
motions, independent of $(h_{ij,0})$. More precisely, in case~$W$ is a
complex Hermitian Wigner matrix, we choose the entries~$(b_{ij,t})$ to
have variance~$t$; in case~$W$ is a real symmetric Wigner matrix, we
choose the off-diagonal entries of $(b_{ij,t})$ to have variance $t$,
while the diagonal entries are chosen to have variance~$2t$. Let
$H_t=(h_{ij,t})$ satisfy the
stochastic differential equation
%
\begin{eqnarray}
\label{lesde} \dd h_{ij}=\frac{\dd b_{ij}}{\sqrt N}-\frac{1}{2}h_{ij}
\,\dd t \qquad(t\ge0).
\end{eqnarray}
It is then easy to check that the distribution of $H_{t}$ agrees with
the distribution of the matrix
%
\begin{eqnarray}
\label{thenewmatrix} \mathrm{e}^{-(t-t_0)/2}V+\mathrm{e}^{-t/2}W+\bigl(1-
\mathrm{e}^{-t}\bigr)^{1/2}W',
\end{eqnarray}
where~$W'$ is, in case~$W$ is a complex Hermitian, a GUE matrix,
independent of~$V$ and~$W$, respectively, a GOE matrix, independent
of~$V$ and~$W$, in case~$W$ is a real symmetric Wigner matrix. The law
of the eigenvalues of the matrix~$W'$ is explicitly given by~(\ref
{gaussianmeasures}) with $\beta=2$, respectively, $\beta=1$.

Denote by $\bflambda(t)=(\lambda_{1}(t),\lambda_{2}(t),\ldots,\lambda
_{N}(t) )$ the ordered eigenvalues of $H_t$. It is well known
that $\bflambda(t)$ satisfy the following stochastic differential equation:
%
\begin{eqnarray}
\label{dysonbrownianmotion} \dd\lambda_{i}=\frac{\sqrt{2}}{\sqrt{\beta
N}}\, {\dd
b_{i}}+ \Biggl(-\frac{\lambda_{i}}{2}+\frac{1}{N}\sum
_{j}^{(i)}\frac{1}{\lambda
_{i}-\lambda_{j}} \Biggr)\,\dd t\qquad\bigl(i
\in\llbracket1,N\rrbracket\bigr),
\end{eqnarray}
where $(b_i)$ is a collection of real-valued, independent standard
Brownian motions. If the matrix $(b_{ij})$ in~(\ref{lesde}) is real
symmetric, we have $\beta=1$ in~(\ref{dysonbrownianmotion}),
respectively, $\beta=2$, if $(b_{ij})$ is complex Hermitian. The
evolution of~$\bflambda(t)$ is the celebrated Dyson Brownian
motion~\cite{D}.

For $t\ge0$, we denote by $f_t \mu_G$ the distribution of $\bflambda
(t)$. In particular, $\int f_t \,\dd\mu_G\equiv\int
f_t(\bflambda) \mu_G(\dd\bflambda)=1$. Note that $f_t \mu_G$
depends on~$V$ through the initial condition $f_{0} \mu_G$. In the
following we always keep the $(v_i)$ fixed; that is, we condition
on~$V$. For simplicity, we omit this conditioning from our notation.
The density $f_t$ is the solution of the equation
\begin{eqnarray*}
\partial_tf_t=\caL f_t\qquad(t\ge0),
\end{eqnarray*}
where the generator $\caL$ is defined via the Dirichlet form
%
\begin{eqnarray}
\label{dirichletform} D_{\mu_G}(f) = -\int f \caL f \,\dd{\mu_G} =
\sum_{i=1}^N \frac
{1}{\beta N} \int(
\partial_i f)^2 \,\dd{\mu_G}\qquad(\partial_i \equiv\partial_{x_i}).
\end{eqnarray}
Formally, we have $\caL= \frac{1}{\beta N}\Delta- (\nabla\caH
_G)\cdot\nabla$, that is,
%
\begin{eqnarray}
\label{legenerator} \caL=\sum_{i=1}^N
\frac{1}{\beta N}\partial_i^2+\sum
_{i=1}^N \Biggl(-\frac{1}{2}
\lambda_i+\frac{1}{N}\sum_{j}^{(i)}
\frac{1}{\lambda
_i-\lambda_j} \Biggr) \partial_i.
\end{eqnarray}
We remark that we use a different normalization in the definition of
the Dirichlet form $D_{\mu_G}(f)$ in~(\ref{dirichletform}) (and the
generator $\caL$) than in earlier works, as in, for example,~\cite
{EY}, where the Dirichlet from is defined as $\sum_{i=1}^N\frac
{1}{2N}\times\break \int(\partial_i f)^2\,\dd\mu_G$.

%
\begin{lemma}[(Dyson Brownian motion)]
The equation $\partial_t f_t=\caL f_t$, with initial data
$f_t\mid_{t=0}=f_{0}$ has a unique solution on $\mathrm{L}^1(\mu
_G)\equiv\mathrm{L}^1(\R^N,\mu_G)$ for all $t\ge0$. Moreover, the
domain $\digamma^{(N)}$ is invariant under the dynamics; that is, if
$f_{0}$ is supported in $\digamma^{(N)}$, then is $f_t$ for all $t\ge0$.
\end{lemma}

(Strictly speaking, the eigenvalue distribution of $H_0$ may not allow
a density $f_0$, but for $t>0$, $H_t$ admits a density $f_t$. Our
proofs are not affected by this technicality.)

We refer, for example, to~\cite{AGZ} for more details and proofs. To
conclude, we record one of the technical tools used in the next sections.

\begin{lemma}\label{lemmamasterrigiditybound}
Denote by $f_t(\bflambda) \mu_G(\dd\bflambda)$ the distribution of
the eigenvalues of matrix~(\ref{thenewmatrix}) with $t\ge0$. Then,
for any $0<\mathfrak{a}<1/2$, we have
%
\begin{eqnarray}
\label{masterrigiditybound} \sup_{t\ge0}\int\frac{1}{N}\sum
_{i=1}^N\bigl(\lambda_i-\widehat
\gamma_i(t)\bigr)^2f_t(\bflambda)\,\dd\mu(
\bflambda)\le N^{-1-2\mathfrak{a}},
\end{eqnarray}
on $\Omega$ for $N$ sufficiently large, where $(\widehat\gamma
_i(t))$ denote the classical locations with respect to the measure
$\widehat\rho_{\mathrm{fc}}(t)$; that is, they are defined through the relation
%
\begin{eqnarray}
\label{classicallocationsoncemore} \int_{-\infty}^{\widehat\gamma
_i(t)}\widehat
\rho_{\mathrm{fc}}(t,x)\,\dd x=\frac{i-\sklfrac{1}{2}}{N}\qquad(1\le i\le N).
\end{eqnarray}
[They agree with the classical locations of~(\ref{nochmalsclassicallocations}).]
\end{lemma}

\begin{pf}
The random matrix $W_t\equiv(w_{ij,t}): =\mathrm
{e}^{-t/2}W+(1-\mathrm{e}^{-t})^{1/2}{W'}$ satisfies the assumptions
in Definition~\ref{assumptionwigner}: the entries are centered and
have variance $1/N$. Moreover, since the distributions of $(w_{ij,0})$,
satisfies~(\ref{eqC0}) and since $(w'_{ij})$ are real, respectively,
complex, centered Gaussian random variables with variance $1/N$,
respectively, $2/N$, the distributions of $(w_{ij,t})$ also
satisfy~(\ref{eqC0}). The claim now follows from~(\ref
{rigidityofeigenvaluesequation}) of Corollary~\ref
{rigidityofeigenvalues} and
the moment bounds $ \E\Tr W_t^{2p}\le C_p$ (see, e.g.,~\cite{AGZ}),
as well as the boundedness of $(v_i)$.
\end{pf}

\subsection{Entropy decay estimates}

Let $\omega$ and $\nu$ be two (probability) measures on $\R^N$ that
are absolutely continuous with respect to Lebesgue measure. We denote
the Radon--Nikodym derivative of $\nu$ with respect to $\omega$ by
$\frac{\dd\nu}{\dd\omega}$, define the relative entropy of $\nu$
with respect to $\omega$ by
%
\begin{eqnarray}
\label{definitionofrelativentropy} S(\nu\mid\omega) :=\int_{\R^N}
\frac{\dd\nu}{\dd
\omega} \log\frac{\dd\nu}{\dd\omega} \,\dd\omega,
\end{eqnarray}
and, in case $\nu=f\omega$, $f\in\mathrm{L}^1(\R^N)$, abbreviate
\begin{eqnarray*}
S_{\omega}(f)=S(f\omega\mid\omega).
\end{eqnarray*}
The entropy $S_{\omega}(f)$ controls the total variation norm of $f$
through the inequality
%
\begin{eqnarray}
\label{entropyinequality} \int\llvert f-1\rrvert\,\dd\omega\le\sqrt{2
S_{\omega}(f)},
\end{eqnarray}
a result we will use repeatedly in the next sections.

Besides the dynamics $(f_t)_{t\ge0}$ generated by $\caL$ introduced
in Section~\ref{subsectionledysonbrownianmotion}, we also consider
a (a priori undetermined) time dependent density, $(\tilde{\psi
}_t)_{t\ge0}$, with respect to $\mu_G$. We assume that $\tilde
{\psi}_t\neq0$, almost everywhere with respect to $\mu_G$ and
abbreviate $\tilde g_t: =\frac{f_t}{\tilde
{\psi}_t}$.
Setting $\tilde\omega_t: =\tilde{\psi}_t
\mu_G$, we can write
\begin{eqnarray*}
f_t(\bflambda) \mu_G(\dd\bflambda)=\tilde
g_t(\bflambda) \tilde\omega_t(\dd\bflambda).
\end{eqnarray*}
A natural choice for $\tilde{\psi}_t \mu_G$ is the time
dependent $\beta$-ensemble, ${\widehat\psi_t} \mu_G$, introduced
in~(\ref{definitionofbt}). Yet, following the arguments of Erd\H
{o}s et al.~\cite{ESYY} we make a slightly different choice for
$\tilde{\psi}_t$: for $\tau>0$, we define a measure $\tilde
{\psi}_t \mu_G$ on $\digamma^{(N)}$ by setting
%
\begin{eqnarray}
\label{definitionofpsit}
\tilde{\psi}_t(\bflambda) \mu_G(
\dd\bflambda): =\frac
{1}{Z_{\tilde{\psi}_t}'}\mathrm{e}^{-N\beta\sum_{i=1}^N\afrac
{(\lambda_i-\widehat\gamma_i(t))^2}{2\tau}}{
\widehat\psi_t}(\bflambda) \mu_G(\dd\bflambda),
\end{eqnarray}
where $Z_{\tilde{\psi}_t}'\equiv Z_{\tilde{\psi}_t}'(\beta
)$ is chosen such that $ \int\tilde{\psi}_t(\bflambda) \mu
_G(\dd\bflambda)=1$. In the following, we mostly choose $\tau$ to be
$N$-dependent with $1\gg\tau>0$.

We call the measure $\tilde{\psi}_t \mu_G$ the instantaneous
relaxation measure. The density $\tilde{\psi}_t$ depends on
$V=\operatorname{diag}(v_i)$ via the initial condition $\tilde\psi
_{0}$. As for the distribution~$f_t$, we condition on~$V$ and omit this
from the notation. We may write the measure $\tilde{\psi}_t \mu
_G$ in the Gibbs form
\begin{eqnarray*}
\tilde{\psi}_t(\bflambda)\mu_G(\dd\bflambda)=
\frac
{1}{Z_{\tilde{\psi}_t}}\mathrm{e}^{-\beta N\widetilde{\caH
}_t(\bflambda)}\,\dd\bflambda\qquad\bigl(\bflambda\in
\digamma^{(N)}\bigr),
\end{eqnarray*}
with
%
\begin{eqnarray}
\label{instantaneoushamiltonian} \widetilde{\caH}_t(\bflambda)=\caH_G(
\bflambda)+\sum_{i=1}^N \biggl(
\frac{(\lambda_i-\widehat\gamma_i(t))^2}{2\tau}+\frac{\widehat
U(t,\lambda_i)}{2} \biggr),
\end{eqnarray}
where $\caH_G$ is defined in~(\ref{gaussianhamiltonian}) and
$Z_{\tilde\psi_t}\equiv Z_{\tilde\psi_t}(\beta)$ is a
normalization. Then we compute
%
\begin{eqnarray}
\label{instantaneousconvexitybound} \nabla u\cdot\bigl(\nabla^2
\widetilde{\caH}_t\bigr)\cdot\nabla u&\ge&\sum_{i=1}^N(
\partial_i u)^2 \biggl(\frac{1}{\tau}+
\frac{\widehat
U''(t,\lambda_i)}{2}+\frac{1}{2} \biggr)
\nonumber
\\
&&{} +\frac
{1}{N}\sum_{i=1}^N
\sum_{j}^{(i)}\frac{1}{(x_i-x_j)^2}(
\partial_i u-\partial_j u)^2
\\
&\ge&\sum_{i=1}^N\frac{(\partial_i u)^2}{2\tau},\nonumber
\end{eqnarray}
for $u\in C^1(\R^N)$ and $\tau$ sufficiently small (independent of
$N$), where we use that $\widehat U''(t,\cdot)$ is uniformly bounded
below by Lemma~\ref{superlemma}. Then, by the Bakry--\'{E}mery
criterion~\cite{BE}, there is a constant $C$ such that the following
logarithmic Sobolev inequality holds for all sufficiently small $\tau>0$:
%
\begin{eqnarray}
\label{logsobolev} S_{\tilde\omega_t}(q)\le C \tau D_{\tilde\omega
_t}(\sqrt{q})\qquad(t\ge0),
\end{eqnarray}
where $q\in\mathrm{L}^\infty(\dd\tilde\omega_t)$ is such that
$\int q \,\dd\tilde\omega_t=1$. We refer, for example, to~\cite
{ESY4,ESYY,EY,EYY} for more details.

Recall the definition of ${\widehat\psi_t} \mu_G$ in~(\ref
{definitionofbt}). Let $\widehat\caL_t$ denote the generator defined
by the natural Dirichlet form with respect to $\widehat\omega_t$,
that is,
%
\begin{eqnarray}
\label{lethewidehatlt} D_{\widehat\omega_t}(q)=\frac{1}{\beta N}\sum
_{i=1}^N\int(\partial_i
q)^2\,\dd\widehat\omega_t=-\int q \widehat
\caL_t q \,\dd\widehat\omega_t\qquad(t> 0).
\end{eqnarray}
The main result of this section is the following proposition.

\begin{proposition}\label{thesuperproposition}
Let $\widehat g_t: =f_t/\widehat{\psi}_t$, and set
$\widehat\omega
_t: =\widehat{\psi}_t \mu_G$ such that
\begin{eqnarray*}
S(f_t \mu_G\mid\widehat{\psi}_t
\mu_G)= S_{\widehat\omega
_t}(\widehat g_t).
\end{eqnarray*}
Then there is a constant~$C$ (independent of $t$) such that, for all
$0<\mathfrak{a}<1/2$, we have
%
\begin{eqnarray}
\label{thesuperpropositionbound} \partial_t S_{\widehat\omega
_t}(\widehat
g_t)\le- 4 D_{\widehat
\omega_t}(\sqrt{\widehat g_t})+C{N^{1-2\mathfrak{a}}}\qquad(t> 0),
\end{eqnarray}
for $N$ sufficiently large on $\Omega$.
\end{proposition}

The results of Proposition~\ref{thesuperproposition} resemble the
relative entropy estimate of Theorem~2.5 in~\cite{EY} for Wigner
matrices. However, due to the fact that both distributions $f_t \mu
_G$ and $\widehat\psi_t \mu_G$ are not close\vspace*{1pt} to the global
equilibrium for the Dyson Brownian motion, $\mu_G$, the reference
ensemble $\widehat\psi_t \mu_G$ changes with time, too. Thus to
establish~(\ref{thesuperpropositionbound}), we need to include
additional factors coming from time derivatives of $\widehat\psi_t
\mu_G$. These can be controlled using the definition of the
potential~$\widehat U(t)$. The idea of choosing slowly varying time
dependent approximation states and controlling the entropy flow goes
back to the work~\cite{Y}.

The relative entropy $S_{\widehat\omega_t}$ and the Dirichlet form
$D_{\widehat\omega_t}$ do not satisfy the logarithmic Sobolev
inequality~(\ref{logsobolev}). However, we have for $t> 0$ the estimates
%
\begin{eqnarray}
D_{\widehat\omega_t}(\sqrt{\widehat g_t})&\le& 2D_{{\tilde\omega
_t}}(\sqrt{\rule{0pt}{7pt}\smash{\tilde g_t}})+C\frac{\beta N^2Q_t}{\tau^2} \label
{dirichletequivalence}
\end{eqnarray}
and
%
\begin{eqnarray}
D_{\tilde\omega_t}(\sqrt{\tilde g_t})&\le& 2D_{\widehat
\omega_t}(
\sqrt{\widehat g_t})+C\frac{\beta N^2Q_t}{\tau^2}\label
{dirichletequivalence2},
\end{eqnarray}
respectively,
%
\begin{eqnarray}
\label{entropyequivalence} S_{\tilde\omega_t}(\tilde g_t)=S_{\widehat
\omega
_t}(
\widehat g_t)+\caO\biggl(\frac{\beta N^2Q_t}{\tau} \biggr),
\end{eqnarray}
where we have set
%
\begin{eqnarray}
\label{leQt} Q_t: =\E^{f_t \mu_G}
\frac{1}{N}\sum_{i=1}^N\bigl(
\lambda_i-\widehat\gamma_i(t)\bigr)^2.
\end{eqnarray}
Estimates~(\ref{dirichletequivalence}),~(\ref{dirichletequivalence2})
and~(\ref{entropyequivalence}) can be checked by elementary
computations, which we omit here. In the following we always bound
$Q_t\le CN^{-1-2\mathfrak{a}}$ [$t\ge0$, $\mathfrak{a}\in(0,1/2)$];
see Lemma~\ref{masterrigiditybound}. Using~(\ref
{dirichletequivalence}),~(\ref{dirichletequivalence2}) and~(\ref
{entropyequivalence}) in combination with the logarithmic Sobolev
inequality~(\ref{logsobolev}) and with Proposition~\ref
{thesuperproposition}, we can follow~\cite{EY} to obtain a bound on the
Dirichlet form $D_{\widehat\omega_t}(\sqrt{\widehat g_t})$.

%
\begin{corollary}\label{lemmaboundonentroypandirichletform}
Under the assumptions of Proposition~\ref{thesuperproposition}, the
following holds on $\Omega$ for $N$ sufficiently large. For any
$\varepsilon'>0$ and $t \geq\tau N^{\varepsilon'}$ with $1\gg\tau\ge
N^{-2\mathfrak{a}}$, we have the entropy and Dirichlet form bounds
%
\begin{eqnarray}
\label{boundonentropyanddirichletform} S_{\widehat\omega_t}(\widehat
g_t)\le C
\frac{N^{1-2\mathfrak
{a}}}{\tau},\qquad D_{\widehat\omega_t}(\sqrt{\widehat g_t})
\leq C \frac{N^{1-2\mathfrak{a}}}{\tau^2},
\end{eqnarray}
where the constants depend on $\varepsilon'$.
\end{corollary}

Before we prove Proposition~\ref{thesuperproposition}, we obtain
rigidity estimates for the time dependent $\beta$-ensemble ${\widehat
\psi_t} \mu_G$. Recall that we denote by~$(\widehat\gamma_i(t))$
the classical locations with respect to the measure~$\widehat\rho
_{\mathrm{fc}}(t)$. Also recall the notation $\check\alpha_i=\min\{i,N-i+1\}$.

\begin{lemma}\label{rigidityfortimedependentbeta}
Let $\widehat U(t,\cdot)$, $t\ge0$ be as in Lemma~\ref{superlemma}.
Then the following holds on $\Omega$ for $N$ sufficiently large:

For any $\delta>0$, there is $\varsigma>0$ such that
%
\begin{eqnarray}
\label{rigidityforbensemble} \mathbb{P}^{{\widehat\psi_t}\mu_G} \bigl
(\bigl\llvert
\lambda_i-\widehat\gamma_i(t)\bigr\rrvert>
N^{-\sklfrac{2}{3}+ \delta}\check\alpha_i^{-\sfrac{1}{3}} \bigr)\leq
\mathrm{e}^{- N^\varsigma},
\end{eqnarray}
for all $t\ge0$, $1\le i\le N$, where $\mathbb{P}^{{\widehat\psi
_t}\mu_G}$, stands for the probability under ${\widehat\psi_t} \mu
_{G}$ conditioned on $\Omega$. Moreover, for any $0<\mathfrak
{a}<1/2$, we have
%
\begin{eqnarray}
\label{masterrigidityb} \sup_{t\ge0}\int\frac{1}{N}\sum
_{i=1}^N\bigl(\lambda_i-\widehat
\gamma_i(t)\bigr)^2 {\widehat\psi_t}(
\bflambda)\mu_G(\dd\bflambda)\le N^{-1-2\mathfrak{a}},
\end{eqnarray}
for $N$ sufficiently large.
\end{lemma}

\begin{pf}
The rigidity estimate (\ref{rigidityforbensemble}) follows from
Proposition~\ref{rigidityfortimedependentbetaone} by choosing
$N\in\N$ sufficiently large. Estimate~(\ref{masterrigidityb}) is a
direct consequence of~(\ref{rigidityforbensemble}) and the fast
decay of the distribution ${\widehat\psi_t}(\bflambda) \mu
_G(\bflambda)$.
\end{pf}
For brevity, we often drop the $t$-dependence of $\widehat\gamma
_i(t)$ from the notation.

\begin{pf*}{Proof of Proposition~\ref{thesuperproposition}}
Recall that we have set $\widehat g_t=f_t/\widehat{\psi}_t$ and
$\widehat\omega_t=\widehat{\psi}_t \mu_G$. The relative entropy
$S(f_t \mu_G\mid\widehat{\psi}_t \mu_G)=S_{\widehat\omega
_t}(\widehat g_t)$ satisfies~\cite{Y},
%
\begin{eqnarray}
\label{landauginzburgformula} \qquad\partial_{t}S(f_t \mu_G
\mid\widehat{\psi}_t \mu_G) = -\frac
{1}{\beta N}\int
\frac{\llvert \nabla\widehat g_t\rrvert ^2}{g_t} \widehat{\psi}_t \,\dd
\mu_G +\int
\frac{(\ls- \partial_t)\widehat{\psi}_t}{\widehat{\psi}_t} f_t \,\dd\mu_G.
\end{eqnarray}
We note that the first term on the right-hand side of~(\ref
{landauginzburgformula}) equals
%
\begin{eqnarray}
\label{landauginzburglemma0} -\frac{1}{\beta N}\int\frac{\llvert \nabla
\widehat g_t\rrvert ^2}{\widehat g_t} \widehat{
\psi}_t \,\dd\mu_G=-4D_{\widehat\omega_t}(\sqrt{\widehat
g_t}).
\end{eqnarray}
To bound the second term on the right-hand side of~(\ref
{landauginzburgformula}), we write
%
\begin{eqnarray}
\label{entropyestimate1}
&& \int\frac{(\ls- \partial_t)\widehat{\psi
}_t}{\widehat{\psi}_t} f_t \,\d\mu_G
\nonumber\\[-8pt]\\[-8pt]\nonumber
&&\qquad =\int(\widehat{\caL}_t \widehat g_t)\,\dd\widehat
\omega_t
 +\frac{1}{2}\int\sum_{i=1}^N
\widehat U'(t,\lambda_i) \bigl(\partial_i
\widehat g_t(\bflambda)\bigr)\,\dd\widehat\omega_t(
\bflambda)-\int\widehat g_t\partial_t\widehat{
\psi}_t\,\d\mu_G,\hspace*{-20pt}
\end{eqnarray}
with $\widehat\caL_t$ defined in~(\ref{lethewidehatlt}).

Note that the first term on the right-hand side of~(\ref
{entropyestimate1}) vanishes since, by construction, $\widehat\omega_t$ is
the reversible measure for the instantaneous flow generated by
$\widehat{\caL}_t$. The last term on the right-hand side of~(\ref
{entropyestimate1}) can be computed explicitly as (recall that the
normalization $Z_{\widehat{\psi}_t}$ in the definition of $\widehat
{\psi}_t \mu_G$ also depends on $t$),
%
\begin{eqnarray}
\label{neue1} -\int\widehat g_t\partial_t\widehat{
\psi}_t\,\d\mu_G &=& \bigl[\E^{f_t \mu_G}-
\E^{\widehat\psi_t \mu_G} \bigr] \Biggl[\frac
{\beta N}{2}\sum
_{i=1}^N\partial_t \widehat U(t,
\lambda_i) \Biggr].
\end{eqnarray}
To deal with the second term on the right-hand side of~(\ref
{entropyestimate1}), we integrate by parts to find
%
\begin{eqnarray}
\label{freude}
&& \frac{1}{2}\int\sum_{i=1}^N
\widehat U'(t,\lambda_i) \bigl(\partial
_i\widehat g_t(\bflambda)\bigr)\,\dd\widehat
\omega_t(\bflambda)
\nonumber
\\
&& \qquad
=\E^{f_t
\mu_G} \Biggl[-\frac{1}{2}\sum
_{i=1}^N \widehat U''(t,
\lambda_i) \Biggr]
\\
&&\quad\qquad{}
+\E^{f_t \mu_G} \Biggl[\frac{\beta
N}{4}\sum_{i=1}^N
\widehat U'(t,\lambda_i) \Biggl(\widehat
U'(t,\lambda_i)+\lambda_i-
\frac{2}{N}\sum_{j}^{(i)}
\frac{1}{\lambda_i-\lambda_j} \Biggr) \Biggr].\hspace*{-20pt}\nonumber
\end{eqnarray}
Setting $\widehat g_t\equiv1$ in the above computation, we also obtain
the identity
%
\begin{eqnarray}
\label{loopeq2} 0&=&\E^{\widehat\psi_t \mu_G} \Biggl[-\frac{1}{2}\sum
_{i=1}^N \widehat U''(t,
\lambda_i) \Biggr]
\nonumber\\[-8pt]\\[-8pt]\nonumber
&&{} +\E^{\widehat
\psi_t \mu_G} \Biggl[\frac{\beta N}{4}\sum
_{i=1}^N \widehat U'(t,
\lambda_i) \Biggl(\widehat U'(t,\lambda_i)+
\lambda_i-\frac
{2}{N}\sum_{j}^{(i)}
\frac{1}{\lambda_i-\lambda_j} \Biggr) \Biggr].
\end{eqnarray}
Equation~(\ref{loopeq2}) may alternatively be derived from the
``first order loop equation'' for the $\beta$-ensemble $\widehat\psi
_t \mu_G$. Equation~(\ref{freude}) can thus be rewritten as
%
\begin{eqnarray}
\label{neue35}
&& \frac{1}{2}\int\sum_{i=1}^N
\widehat U'(t,\lambda_i) \bigl(\partial
_i\widehat g_t(\bflambda)\bigr)\,\dd\widehat
\omega_t(\bflambda)\nonumber
\\
&&\qquad = \bigl[\E^{f_t \mu_G}-\E^{\widehat\psi_t \mu_G}
\bigr] \Biggl[-\frac
{1}{2}\sum_{i=1}^N
\widehat U''(t,\lambda_i) \Biggr]
\nonumber\\[-8pt]\\[-8pt]\nonumber
&&\quad\qquad{}+ \bigl[\E^{f_t \mu_G}-\E^{\widehat\psi_t \mu_G} \bigr]
\\
&&\qquad\qquad{}\times \Biggl[\frac
{\beta N}{4}
\sum_{i=1}^N \widehat U'(t,
\lambda_i) \Biggl(\widehat U'(t,\lambda_i)+
\lambda_i-\frac{2}{N}\sum_{j}^{(i)}
\frac
{1}{\lambda_i-\lambda_j} \Biggr) \Biggr].
\nonumber
\end{eqnarray}

Next, to control the second and third terms on the right-hand side
of~(\ref{entropyestimate1}), respectively, the right-hand side
of~(\ref{neue35}), we\vspace*{1pt} proceed as follows. We expand the potential
terms $\widehat U'(t,\lambda_i)$, respectively, $\widehat U''(t,\lambda
_i)$, in Taylor series in $\lambda_i$ to second order around the
classical location $\widehat\gamma_i$. The resulting zero order terms
cancel exactly since the classical locations of the ensembles $f_t \mu
_G$, and $\widehat\psi_t \mu_G$ agree by construction. The first
order terms in the Taylor expansion can~(1) either be bounded in terms
of the expectations of $\sum_{i=1}^N{(\lambda_i-\widehat\gamma
_i)^2}$ (which can be controlled with the rigidity estimates in
Lemmas~\ref{rigidityfortimedependentbeta} and~\ref
{lemmamasterrigiditybound}); or~(2) they cancel exactly due to the~definition of the potential $\widehat U(t,\cdot)$ and its equation of
motion in~(\ref{EoMdrei}). Finally, the second order
terms in the Taylor expansion can be bounded by the rigidity estimates
in Lemmas~\ref{rigidityfortimedependentbeta} and~\ref
{lemmamasterrigiditybound}. The details are as follows.

Expanding $\partial_t \widehat U(t,\lambda_i)$ to second order around
$\widehat\gamma_i$, we obtain from~(\ref{neue1}) that
%
\begin{eqnarray}
\label{neue2}
-\int\widehat g_t\partial_t\widehat{
\psi}_t\,\d\mu_G
&=& \bigl[\E^{f_t \mu_G}-
\E^{\widehat\psi_t \mu_G} \bigr] \nonumber
\\
&&{}\times
\Biggl[\frac{\beta
N}{2}\sum
_{i=1}^N\partial_t \widehat U(t,
\widehat\gamma_i)+\frac{\beta N}{2}\sum_{i=1}^N
\partial_t \widehat U'(t,\gamma_i) (
\lambda_i-\widehat\gamma_i) \Biggr]
\\
&&{} +\caO
\bigl(N^{1-2\mathfrak{a}}\bigr),\nonumber
\end{eqnarray}
on $\Omega$, where we use the rigidity estimates in Lemmas~\ref
{rigidityfortimedependentbeta} and~\ref
{lemmamasterrigiditybound}, and that $\partial_t\widehat U''(t,\cdot)$ is
uniformly bounded on compact sets by Lemma~\ref{superlemma}.\vspace*{1pt}

To save notation, we introduce a function $G\dvtx \R^+\times\R^2\to
\R$ by setting
%
\begin{eqnarray}
\label{leG} G(t;x,y): =\frac{\widehat U'(t,x)-\widehat
U'(t,y)}{x-y},
\end{eqnarray}
with $G(t;x,x): =\widehat U''(t,x)$. Note that
$G(t;x,y)=G(t;y,x)$
and that $G$ is $C^2$ in the spatial coordinates by Lemma~\ref
{superlemma}. Recalling the equation of motion for $\partial_t
\widehat U(t,\cdot)$ in~(\ref{EoMdrei}), we can write
%
\begin{eqnarray}
\qquad \partial_t \widehat U(t,x)&=&\dashint\frac{\widehat U'(t,y)\,\dd
\widehat\rho_{\mathrm{fc}}(t,y)}{y-x}
\nonumber\\[-8pt]\\[-8pt]\nonumber
&=& \int\frac{\widehat U'(t,y)-\widehat U'(t,x)}{y-x}\,\dd\widehat\rho
_{\mathrm{fc}}(t,y)+ \widehat
U'(t,x)\dashint\frac{\dd\widehat\rho
_{\mathrm{fc}}(t,y)}{y-x},
\end{eqnarray}
for $x$ inside the support of the measure $\widehat\rho_{\mathrm{fc}}$. Thus,
recalling~(\ref{definitionwidehatU}) and~(\ref{leG}), we obtain
%
\begin{eqnarray}
\partial_t \widehat U(t,x)=\int G(t;x,y)\,\dd\widehat\rho
_{\mathrm{fc}}(y)-\frac{1}{2}\widehat U'(t,x) \bigl(
\widehat U'(t,x)+x \bigr),
\end{eqnarray}
for $x$ inside the support of the measure $\widehat\rho_{\mathrm{fc}}$.

We hence obtain from~(\ref{neue2}) that
%
\begin{eqnarray}
\label{neue36}
&& -\int\widehat g_t\partial_t\widehat{
\psi}_t \,\dd\mu_G\nonumber
\\
&&\qquad = \bigl[\E^{f_t \mu_G}-
\E^{\widehat\psi_t \mu_G} \bigr] \Biggl[\frac
{\beta N}{2}\sum_{i=1}^N\int G'(t;\widehat
\gamma_i,y)\,\dd\widehat\rho_{\mathrm{fc}}(y) (
\lambda_i-\widehat\gamma_i) \Biggr]
\nonumber
\\
&&\quad\qquad{} - \bigl[\E^{f_t \mu_G}-\E^{\widehat\psi_t \mu_G} \bigr] \Biggl[\frac
{\beta N}{4}
\sum_{i=1}^N\widehat U''(t,
\widehat\gamma_i) \bigl(\widehat U'(t,\widehat
\gamma_i)+\widehat\gamma_i \bigr) (\lambda
_i-\widehat\gamma_i) \Biggr]
\\
&&\quad\qquad{} - \bigl[\E^{f_t \mu_G}-\E^{\widehat\psi_t \mu_G} \bigr] \Biggl[\frac
{\beta N}{4}
\sum_{i=1}^N\widehat U'(t,
\widehat\gamma_i) \bigl(\widehat U''(t,
\widehat\gamma_i)+1 \bigr) (\lambda_i-\widehat
\gamma_i) \Biggr]
\nonumber
\\
&&\quad\qquad{} +\caO\bigl(N^{1-2\mathfrak{a}}\bigr),\nonumber
\end{eqnarray}
on $\Omega$, where we denote by $G'(t;x,y)$ the first derivative of
$G(t;x,y)$ with respect to $x$.

Next we return to~(\ref{neue35}). Using the rigidity estimates of the
Lemmas~\ref{rigidityfortimedependentbeta} and~\ref
{lemmamasterrigiditybound}, we find
%
\begin{eqnarray}
\label{laufet}
&& \frac{1}{2}\int\sum_{i=1}^N
\widehat U'(t,\lambda_i) \bigl(\partial
_i\widehat g_t(\bflambda)\bigr)\,\dd\widehat
\omega_t(\bflambda)\nonumber
\\
&&\qquad = \bigl[\E^{f_t \mu_G}-\E^{\widehat\psi_t \mu_G}
\bigr] \Biggl[-\frac
{1}{2}\sum_{i=1}^N
\widehat U''\bigl(t,\widehat\gamma_i(t)
\bigr) \Biggr]
\nonumber
\\
&&\quad\qquad{}  + \bigl[\E^{f_t \mu_G}-\E^{\widehat\psi_t \mu
_G} \bigr]
\\
&&\qquad\qquad{}\times  \Biggl[\frac{\beta N}{4}
\sum_{i=1}^N \widehat U'(t,
\lambda_i) \Biggl(\widehat U'(t,\lambda_i)+
\lambda_i-\frac{2}{N}\sum_{j}^{(i)}
\frac{1}{\lambda_i-\lambda_j} \Biggr) \Biggr]
\nonumber
\\
&&\quad\qquad{}  +\caO\bigl(N^{1/2-\mathfrak{a}} \bigr),
\nonumber
\end{eqnarray}
on\vspace*{2pt} $\Omega$, where we use a Taylor expansion of the first term on the
right-hand side of~(\ref{neue35}). Here we also use that $\widehat
U'$ is three times continuously differentiable with uniformly
bounded derivatives on compact sets. Note that the first term on the
right-hand side of~(\ref{laufet}) vanishes.

Using the definition of $G(t;\cdot,\cdot)$ in~(\ref{leG}), we can
recast~(\ref{laufet})~as
%
\begin{eqnarray}
\label{neue}
&& \frac{1}{2}\int\sum_{i=1}^N
\widehat U'(t,\lambda_i)\bigl(\partial
_i\widehat g_t(\bflambda)\bigr)\,\dd\widehat
\omega_t(\bflambda)\nonumber
\\
&&\qquad =
 \bigl[\E^{f_t \mu_G}-\E^{\widehat\psi_t \mu_G} \bigr] \Biggl[\frac
{\beta N}{4}
\sum_{i=1}^N \widehat U'(t,
\lambda_i) \bigl(\widehat U'(t,\lambda_i)+
\lambda_i \bigr) \Biggr]
\nonumber\\[-8pt]\\[-8pt]\nonumber
&&\quad\qquad{} + \bigl[\E^{f_t \mu_G}-\E^{\widehat\psi_t \mu_G} \bigr] \Biggl[-
\frac{\beta N}{4}\sum_{i=1}^N
\frac{1}{N}\sum_{j}^{(i)}G(t;
\lambda_i,\lambda_j) \Biggr]
\nonumber
\\
&&\quad\qquad{}+\caO\bigl(N^{1/2-\mathfrak{a}} \bigr),\nonumber
\end{eqnarray}
where we use the symmetry $G(t;x,y)=G(t;y,x)$. Expanding the second
term on the right-hand side~(\ref{neue}) to second order in $(\lambda
_i,\lambda_j)$ around $(\widehat\gamma_i,\widehat\gamma_j)$, we obtain
%
\begin{eqnarray}
&&  \bigl[\E^{f_t \mu_G}-\E^{\widehat\psi_t \mu_G} \bigr] \Biggl[-\frac
{\beta N}{4}
\sum_{i=1}^N\frac{1}{N}\sum
_{j}^{(i)}G(t;\lambda_i,
\lambda_j) \Biggr]\nonumber
\\
&&\qquad =
\bigl[\E^{f_t \mu_G}-\E^{\widehat\psi_t \mu_G} \bigr] \Biggl[-
\frac{\beta N}{2}\sum_{i=1}^N \Biggl(
\frac{1}{N}\sum_{j=1}^{N}G'(t;
\widehat\gamma_i,\widehat\gamma_i) \Biggr) (\lambda
_i-\widehat\gamma_i) \Biggr]
\\
&&\quad\qquad{} +\caO\bigl(N^{1/2-\mathfrak{a}}\bigr)+\caO\bigl(N^{1-2\mathfrak
{a}}\bigr),\nonumber
\end{eqnarray}
on $\Omega$, where we use $G(t;x,y)=G(t;y,x)$, $G(t;x,x)=\widehat
U''(t,x)$ and that $G(t;x,y)$ is $C^2$ in the spatial variables. Thus,
also expanding the first term on the right-hand side of~(\ref{neue})
in $\lambda_i$ around $\widehat\gamma_i$, we obtain
%
\begin{eqnarray}
\label{neue37}
&& \frac{1}{2}\int\sum_{i=1}^N
\widehat U'(t,\lambda_i) \bigl(\partial
_i\widehat g_t(\bflambda)\bigr)\,\dd\widehat
\omega_t(\bflambda)\nonumber
\\
&&\qquad=
 \bigl[\E^{f_t \mu_G}-\E^{\widehat\psi_t \mu_G} \bigr] \Biggl[-\frac
{\beta N}{2}
\sum_{i=1}^N \biggl(\frac{1}{N}\sum
_{j}G'(t;\widehat
\gamma_i,\widehat\gamma_j) \biggr) (\lambda
_i-\widehat\gamma_i) \Biggr]
\nonumber
\\
&&\quad\qquad{} + \bigl[\E^{f_t \mu_G}-\E^{\widehat\psi_t \mu_G} \bigr] \Biggl[\frac
{\beta N}{4}
\sum_{i=1}^N \widehat U''(t,
\widehat\gamma_i) \bigl(\widehat U'(t,\widehat
\gamma_i)+\widehat\gamma_i \bigr) (
\lambda_i-\widehat\gamma_i) \Biggr]
\\
&&\quad\qquad{} + \bigl[\E^{f_t \mu_G}-\E^{\widehat\psi_t \mu_G} \bigr] \Biggl[\frac
{\beta N}{4}
\sum_{i=1}^N \widehat U'(t,
\widehat\gamma_i) \bigl(\widehat U''(t,
\widehat\gamma_i)+1 \bigr) (\lambda_i-\widehat
\gamma_i) \Biggr]
\nonumber
\\
&&\quad\qquad{}+\caO\bigl(N^{1/2-\mathfrak{a}} \bigr)+\caO\bigl({N^{1-2\mathfrak
{a}}} \bigr),\nonumber
\end{eqnarray}
on $\Omega$, where we use the rigidity estimates in Lemmas~\ref
{rigidityfortimedependentbeta} and~\ref{lemmamasterrigiditybound}.

Adding up~(\ref{neue36}) and~(\ref{neue37}), we hence obtain
\begin{eqnarray*}
&& \biggl\llvert\int\frac{(\ls- \partial_t)\widehat{\psi}_t}{\widehat
{\psi}_t} f_t \,\d\mu_G
\biggr\rrvert
\\
&&\qquad \le
 \frac{\beta N}{2}\Biggl\llvert\E^{f_t \mu_G} \Biggl[-\sum
_{i=1}^N \Biggl(\frac{1}{N}\sum
_{j=1}^NG'(t;\widehat
\gamma_i,\widehat\gamma_j)
\\
&&\hspace*{117pt}{} -\int G'(t;
\widehat\gamma_i,y)\,\dd\widehat\rho_{\mathrm{fc}}(y) \Biggr) (
\lambda_i-\widehat\gamma_i) \Biggr]
\\
&&\hspace*{19pt}\quad\qquad{} -\E^{\widehat\psi_t \mu_G}
\Biggl[-\sum_{i=1}^N \Biggl(
\frac
{1}{N}\sum_{j=1}^NG'(t;
\widehat\gamma_i,\widehat\gamma_j)
\\
&&\hspace*{128pt}{}-\int
G'(t;\widehat\gamma_i,y)\,\dd\widehat
\rho_{\mathrm{fc}}(y) \Biggr) (\lambda_i-\widehat
\gamma_i) \Biggr]\Biggr\rrvert
\nonumber
\\
&&\quad\qquad{}+\caO\bigl(N^{1/2-\mathfrak{a}}\bigr)+\caO\bigl(N^{1-2\mathfrak
{a}}\bigr),
\end{eqnarray*}
on $\Omega$. To finish the proof we observe that for all $\widehat
\gamma_i$,
\begin{eqnarray*}
\frac{1}{N}\sum_{j=1}^{N}
G'(t;\widehat\gamma_i,\widehat\gamma
_j)=\int G'(t;\widehat\gamma_i,y)\,\dd
\widehat\rho_{\mathrm{fc}}(t,y)+\caO\bigl(N^{-1}\bigr),
\end{eqnarray*}
on $\Omega$, where we use that $\widehat\gamma_{i+1}-\widehat\gamma
_{i}\sim N^{-2/3}\check\alpha_i^{-1}$ ($\check\alpha_i=\min\{
i,N-i+1\}$), and the square root decay of $\widehat\rho_{\mathrm{fc}}(t)$ at
the edges of the support. Thus
%
\begin{eqnarray}
\label{freudig} \int\frac{(\caL-\partial_t)\widehat\psi_t}{\widehat\psi
_t}f_t \,\dd\mu_G&=&
\caO\bigl(N^{1/2-\mathfrak{a}} \bigr)+\caO\bigl(N^{1-\mathfrak
{a}}\bigr),
\end{eqnarray}
for $N$ sufficiently large on $\Omega$, where we use one last time the
rigidity estimates. Using that $N^{1/2-\mathfrak{a}}< N^{1-2\mathfrak
{a}}$, $\mathfrak{a}\in(0,1/2)$, we get from~(\ref
{landauginzburgformula}),~(\ref{landauginzburglemma0}) and~(\ref
{freudig}) the
desired estimate~(\ref{thesuperpropositionbound}).
\end{pf*}
Before we move on to the proof of Corollary~\ref
{lemmaboundonentroypandirichletform}, we give a rough estimate on
$S_{\widehat
\omega_t}(\widehat g_t)$ for $t> 0$.

\begin{lemma} \label{basicentropybound}
There is a constant $m$ such that, for $\tau>0$ and $t\ge\tau$, we have
%
\begin{eqnarray}
\label{basicentropyboundeq} S_{\widehat\omega_t}(\widehat g_t)=S(f_t
\mu_G\mid\widehat{\psi}_t \mu_G) \leq
CN^m
\end{eqnarray}
on $\Omega$, for $N$ sufficiently large. Here the constant $C$ depends
on $\tau$.
\end{lemma}

\begin{pf}
From the definition of the relative entropy in~(\ref
{definitionofrelativentropy}), we have
%
\begin{eqnarray}
\label{basicentropyboundI}
&& S\bigl(f_t \mu_G\llvert\widehat{
\psi}_t \mu_G\bigr)
\nonumber\\[-6pt]\\[-10pt]\nonumber
&&\qquad \le S\bigl(f_t
\mu_G\rrvert\mu_G\bigr)
+\Biggl\llvert \frac{\beta N}{2}\sum_{i=1}^N\int
\widehat U(t,\lambda_i) f_t(\bflambda)\,\d
\mu_G(\bflambda)\Biggr\rrvert
+\log Z_{\widehat\psi_t}.
\end{eqnarray}
Since the potential $\widehat U(t)$ is bounded below, we have (for $N$
sufficiently large on $\Omega$) $\log Z_{{\widehat\psi_t}}\le C\beta
N^2$. Similarly, using the rigidity estimate~(\ref
{masterrigiditybound}), we can bound the second term on the right-hand
side of~(\ref
{basicentropyboundI}) by $CN^2$. To bound the first term on the right
of~(\ref{basicentropyboundI}), we use that $S(f_t \mu_G\mid\mu
_G)\le S(H_t\mid W')\le N^2\max S(h_{ij,t}\mid w'_{ij})+N\max
S(h_{ii,t}\mid
w'_{ii})$, where $(h_{ij,t})$ are the entries of the in~(\ref
{thenewmatrix}) and~$w'_{ij}$ are the entries of the GOE, respectively, GUE,
matrix ${W'}$. By explicit calculations, remembering that the diagonal
entries $(v_i)$ are fixed, one finds $\max S(h_{ij,t}\mid g_{ij})\le CN$
for $t\ge\tau$; see, for example,~\cite{EKYY2}. (Note that we choose
$t>0$; otherwise the relative entropy may be ill defined.)
\end{pf}

To complete the proof of Corollary~\ref
{lemmaboundonentroypandirichletform} we follow the discussion in~\cite{EY}.

\begin{pf*}{Proof of Corollary~\ref{lemmaboundonentroypandirichletform}}
Using an approximation argument, we can assume that $\tilde g_t\in
\mathrm{L}^{\infty}(\dd\tilde\omega_t)$. Using first the
entropy bound~(\ref{thesuperpropositionbound}) and then the
Dirichlet form estimate in~(\ref{dirichletequivalence2}), we obtain
\begin{eqnarray}
\partial_t S_{\widehat\omega_t}(\widehat g_t)&
\le&-4D_{\widehat
\omega_t}(\sqrt{\widehat g_t})+C N^{1-2\mathfrak{a}}
\nonumber
\\
&\le&-2D_{\tilde\omega_t}(\sqrt{\rule{0pt}{7pt}\smash{\tilde g_{t}}})+C
N^{1-2\mathfrak{a}}+C \frac{N^{1-2\mathfrak{a}}}{\tau^2}
\nonumber
\\
&\le&-C\tau^{-1} S_{\tilde\omega_t}(\tilde g_t)+C
\frac
{N^{1-2\mathfrak{a}}}{\tau^2},
\nonumber
\end{eqnarray}
for $N$ sufficiently large on $\Omega$. To get the third line we use
the logarithmic Sobolev inequality~(\ref{logsobolev}) and that, by
assumption, $\tau<1$. Using the entropy estimate~(\ref
{entropyequivalence}), we thus obtain
%
\begin{eqnarray}
\label{tobeintegrated} \partial_t S_{\widehat\omega_t}(\widehat
g_t)&\le&-C\tau^{-1} S_{\widehat\omega_t}(\widehat
g_t)+C \frac{N^{1-2\mathfrak
{a}}}{\tau^2},
\end{eqnarray}
for $N$ sufficiently large on $\Omega$. Integrating~(\ref
{tobeintegrated}) from $\tau$ to $t/2$, we infer
\begin{eqnarray*}
S_{\widehat\omega_{t/2}}(\widehat g_{t/2})\le\mathrm{e}^{-C\tau
^{-1}(t/2-\tau)}S_{\widehat\omega_{\tau}}(
\widehat g_{\tau
})+C\frac{N^{1-2\mathfrak{a}}}{\tau},
\end{eqnarray*}
for $N$ sufficiently large on $\Omega$. Bounding $S_{\widehat\omega
_{\tau}}(\widehat g_{\tau})$ by~(\ref{basicentropyboundI}), we get
\begin{eqnarray*}
S_{\widehat\omega_{t/2}}(\widehat g_{t/2})\le CN^m
\mathrm{e}^{-C\tau
^{-1}(t/2-\tau)}+C\frac{N^{1-2\mathfrak{a}}}{\tau},
\end{eqnarray*}
for $N$ sufficiently large on $\Omega$. Recalling that $t\ge\tau
_0=\tau N^{\varepsilon'}$ and using the monotonicity of the relative
entropy, we obtain the first inequality in~(\ref
{boundonentropyanddirichletform}).

Integrating~(\ref{thesuperpropositionbound}) from $t/2$ to $t$, we obtain
\begin{eqnarray*}
\int_{t/2}^{t} D_{\widehat\omega_s}(\sqrt{\widehat
g_s})\,\dd s &\le&- \int_{t/2}^{t}
\partial_sS_{\widehat\omega_s}(\widehat g_s)\,\dd s+
Ct{N^{1-2\mathfrak{a}}}.
\end{eqnarray*}
Thus, using the above estimate on the relative entropy and the
monotonicity of the Dirichlet form,
\begin{eqnarray*}
D_{\widehat\omega_t} (\sqrt{\widehat g_t})\le C\frac
{N^{1-2\mathfrak{a}}}{t \tau}+C{N^{1-2\mathfrak{a}}}.
\end{eqnarray*}
Recalling that $t\ge\tau_0=N^{\varepsilon'}\tau$, we get the second
inequality in~(\ref{boundonentropyanddirichletform}).
\end{pf*}

\section{Local equilibrium measures}\label{localequilibriummeasures}
The estimates on the relative entropy and the Dirichlet form obtained
in Corollary~\ref{lemmaboundonentroypandirichletform} do not
directly imply that the local statistics of the measures $f_t \mu_G$
and $\widehat\psi_t \mu_G$ agree\vspace*{1pt} in the limit of large $N$.
However, the averaged local gap statistics of $f_t \mu_G$, $\widehat
\psi_t \mu_G$ and $\mu_G$ can be compared (for $1\gg t\gg
N^{-1/2}$) for large $N$ as is asserted in the main theorems of this
section, Theorem~\ref{maintheoremofsection6} and Theorem~\ref
{corollarygaussian}, below. We first state these results and give a
short outline of their proofs in Section~\ref
{subsectionAveragedlocalgapstatisticsforsmalltimes} before going into
the details in
Sections~\mbox{\ref{subsection1ofsection6}--\ref{subsectionproofoftheoremandcorollarysection6}}.

\subsection{Averaged local gap statistics for small times}\label{subsectionAveragedlocalgapstatisticsforsmalltimes}

Recall that we call a symmetric function~$O\dvtx \R^n\to\R$, $n\in
\N$, an $n$-particle observable if $O$ is smooth and compactly
supported. For a given observable~$O$, a time~$t\ge0$, a small
constant $\alpha>0$ and $j\in\llbracket\alpha N,(1-\alpha)
N\rrbracket$, we define an observable $G_{j,n,t}(\bfx)\equiv
G_{j,n}(\bfx)$, by setting
%
\begin{eqnarray}
\label{theaveragedobservable} G_{j,n}(\bfx) &:=& O\bigl(N
\rho_{j}(x_{j+1}-x_j),N\rho
_{j}(x_{j+2}-x_j),\ldots,
\nonumber\\[-8pt]\\[-8pt]\nonumber
&&\hspace*{102pt} N\rho_{j}(x_{j+n+1}-x_{j})\bigr),
\end{eqnarray}
$\bfx=(x_k)_{k=1}^N\in\digamma^{(N)}$, where we set $G_{j,n}=0$ if
$j+n>(1-\alpha) N$. Here $\rho_j$ denotes the density of the measure
$\widehat\rho_{\mathrm{fc}}(t)$ at the classical location of the $j$th
particle at time~$t$, that is, $\rho_j: =\widehat
\rho
_{\mathrm{fc}}(t,\widehat\gamma_j(t))$. We also set
%
\begin{eqnarray}
\label{theaveragedobservableatinfty}
G_{j,n,\mathrm{sc}}(\bfx)
&:=&O\bigl(N\rho
_{\mathrm{sc},j}(x_{j+1}-x_j),N\rho_{\mathrm{sc},j}(x_{j+2}-x_j),
\ldots,
\nonumber\\[-8pt]\\[-8pt]\nonumber
&&\hspace*{112pt} N\rho_{\mathrm{sc},j}(x_{j+n+1}-x_{j})\bigr),
\end{eqnarray}
$\bfx\in\digamma^{(N)}$, where $\rho_{\mathrm{sc},j}$ denotes the density of
the semicircle law at the classical location of the $j$th particle with
respect to the semicircle law.

In the following, we denote constants depending on $O$ by $C_O$. Recall
the definition of the density~$\widehat\psi_t$ in~(\ref
{definitionofbt}). We have the following statement on the averaged
local gap statistics.

%
\begin{theorem}\label{maintheoremofsection6}
Let $n\in\N$ be fixed, and consider an $n$-particle observable~$O$.
Fix a small constant~$\alpha>0$, and consider an interval of
consecutive integers $J\subset\llbracket\alpha N,(1-\alpha
)N\rrbracket$ in the bulk. Then, for any small $\delta>0$, there is a
constant $\mathfrak{f}>0$ such that, for $t\ge N^{-1/2+\delta}$,
%
\begin{eqnarray}
\label{equationinmaintheoremofsectoin6}
&& \biggl\llvert\int\frac
{1}{\llvert J\rrvert }\sum
_{j\in J}G_{j,n}(\bfx) f_t(\bfx)\,\dd
\mu_G(\bfx)-\int\frac{1}{\llvert J\rrvert }\sum
_{j\in J}G_{j,n} (\bfx) {\widehat\psi_t}(
\bfx)\,\dd\mu_G(\bfx)\biggr\rrvert
\nonumber\\[-8pt]\\[-8pt]\nonumber
&&\qquad \le C_ON^{-\mathfrak{f}},
\end{eqnarray}
for $N$ sufficiently large on $\Omega$. The constant $C_O$ depends on
$\alpha$ and $O$, and the constant $\mathfrak{f}$ depends on $\alpha
$ and $\delta$.
\end{theorem}

We can also compare the averaged local gap statistics of $f_t \mu_G$,
with the averaged local gap statistics of the Gaussian unitary,
respectively, orthogonal, ensemble.

\begin{theorem}\label{corollarygaussian}
Under the same assumptions as in Theorem~\ref{theaveragedobservable}
and with similar constants, we have
\begin{eqnarray*}
\biggl\llvert\int\frac{1}{\llvert J\rrvert }\sum_{j\in J}G_{j,n}(
\bfx) f_t(\bfx)\,\dd\mu_G(\bfx)-\int\frac{1}{\llvert J\rrvert }
\sum_{j\in J}G_{j,n,\mathrm{sc}}(\bfx)\,\dd\mu
_G(\bfx)\biggr\rrvert\le C_O N^{-\mathfrak{f}},
\end{eqnarray*}
for $N$ sufficiently large on $\Omega$.
\end{theorem}

The proofs of Theorem~\ref{maintheoremofsection6} and Theorem~\ref
{corollarygaussian} proceed in two steps. We first localize the
measures $f_t \mu_G$ and $\widehat\psi_t,\mu_G$; that is, we study
the statistics of~$\caK$, $1\ll\caK\ll N$, consecutive particles
inside the bulk---the interior particles---with the remaining
particles---the exterior particles---being fixed; for details, see
Section~\ref{subsection1ofsection6}. For most configurations of
the exterior particles (boundary conditions), we can compare the
statistics of the localized versions of $f_t \mu_G$ and $\widehat
\psi_t,\mu_G$. This is accomplished in Proposition~\ref
{frozenstatsproximity} of Section~\ref
{subsectionlecomparisionoflocalmeasures} by using that (1) the
localized $\beta$-ensemble satisfies a
logarithmic Sobolev inequality~(\ref{frozenlog-sobolev}) with
constant $C\caK/N$ and that (2) the localized Dirichlet form can be
controlled by the global Dirichlet form [see~(\ref
{lewhereweusethedformestimate})], the latter being estimated in
Corollary~\ref{lemmaboundonentroypandirichletform}.

In a second step, we use Theorem~4.1 of~\cite{singlegap} that,
roughly speaking, assures that the local gap statistics of localized
$\beta$-ensembles are essentially independent of the boundary
conditions and indeed agree with the local gap statistics of the
Gaussian ensembles. Putting this universality result to work in
Section~\ref{gapuniversalityforlocalmeasures}, we conclude that
the local gap statistics of the localized version of the measure $f_t
\mu_G$ are universal, for $1\gg t\gg N^{-1/2}$ and for most boundary
conditions. Theorems~\ref{maintheoremofsection6} and~\ref
{corollarygaussian} are then proven in Section~\ref
{subsectionproofoftheoremandcorollarysection6} by integrating out the
boundary conditions.

We conclude this subsection with the following two remarks: once the
entropy estimate of Proposition~\ref{thesuperproposition} has been
established, one can apply the methods of~\cite{singlegap} to prove the gap universality in the bulk for deformed
Wigner matrices; see Remark~\ref{remarkforgapuniversality} above
for an explicit statement; we leave the details to the interested readers.

As an alternative to the approach outlined above, one could combine the
approach from~\cite{EY} with Theorem~2.1 in~\cite{BEY} (see
Theorem~\ref{thmbulkuniversalitybeta} above), to prove
Theorems~\ref{maintheoremofsection6} and~\ref{corollarygaussian}.

\subsection{Preliminaries}\label{subsection1ofsection6}

Let $\alpha, \sigma>0$ be two small positive numbers, and choose two
integer parameters~$L$ and~$K$ such that
%
\begin{eqnarray}
\label{conditiononK} L\in\dllbracket\alpha N,(1-\alpha)N\drrbracket
,\qquad K\in
\dllbracket N^{\sigma}, N^{1/4}\drrbracket.
\end{eqnarray}
We denote by $ I_{L,K}: =\llbracket L-K,L+K\rrbracket
$ a set of
$\caK: =2K+1$ consecutive indices in the bulk of the
spectrum. Below
we often abbreviate $I\equiv I_{L,K}$. Recall the definition of the set
$\digamma^{(N)}\subset\R^N$ in~(\ref{Weylchamber}). For $\bflambda
\in\digamma^{(N)}$, we write
%
\begin{eqnarray}
\label{thebflambda} \bflambda=(y_1,\ldots, y_{L-K-1},x_{L-K},
\ldots, x_{L+K},y_{L+K+1},\ldots,y_{N}),
\end{eqnarray}
and we call $\bflambda$ a configuration (of $N$ particles or points on
the real line). Note that on the right-hand side of~(\ref{thebflambda})
the points keep their original indices and are in increasing
order so that
%
\begin{eqnarray}
\label{bfxnotation} \bfx&=&(x_{L-K},\ldots, x_{K+L})\in
\digamma^{(\caK)},
\nonumber\\[-8pt]\\[-8pt]\nonumber
\bfy&=&(y_1,\ldots, y_{L-K-1},y_{L+K+1},\ldots,
y_N)\in\digamma^{(N-\caK)}.
\end{eqnarray}
We refer to $\bfx$ as the interior points or particles and to $\bfy$
as the exterior points or particles.

In the following, we often fix the exterior points and consider the
conditional measures on the interior points: let~$\omega$ be a measure
on~$\digamma^{(N)}$ with a density. Then we denote by $\omega^{\bfy
}$ the measure obtained by conditioning on~$\bfy$; that is,
for~$\bflambda$ in the form of~(\ref{thebflambda}),
\begin{eqnarray*}
\omega^{\bfy}(\dd\bfx)\equiv\omega^{\bfy}(\bfx)\,\dd\bfx
: = \frac{\omega(\bflambda)\,\dd\bfx}{\int\omega(\bflambda
)\,\dd\bfx
}=\frac{\omega(\bfx,\bfy)\,\dd\bfx}{\int\omega(\bfx,\bfy)\,\dd
\bfx},
\end{eqnarray*}
where, with slight abuse of notation, $\omega(\bfx,\bfy)$ stands for
$\omega(\bflambda)$. We refer to the fixed exterior points $\bfy$ as
boundary conditions of the measure $\omega^{\bfy}$. For fixed $\bfy
\in\digamma^{(N-\caK)}$, all $(x_i)$ lie in the open configuration interval
\begin{eqnarray*}
\mathbf{I}\equiv\mathbf{I}_{L,K}: =(y_{L-K-1},y_{L+K+1}).
\end{eqnarray*}
Set $\bar{y} : =(y_{L-K-1}+y_{L+K+1})/2$, and let
%
\begin{eqnarray}
\label{equidistantpoints1} \alpha_j: =\bar y+
\frac{j-L}{\caK+1}\llvert\bfy\rrvert\qquad(j\in I_{L,K})
\end{eqnarray}
denote $\caK$ equidistant points in the interval $\bsI$.

Let $U\in C^4(\R)$ be a ``regular'' potential satisfying~(\ref
{assumption1forbetauniversality}) and~(\ref
{assumption2forbetauniversality}). We then consider the $\beta$-ensemble
%
\begin{eqnarray}
\label{betaensemblesection6} \mu(\dd\bflambda)\equiv\mu_U (\dd\bflambda)
: =\frac
{1}{Z_{U}}\mathrm{e}^{-\beta N \caH(\bflambda) }\,\dd\bflambda\qquad(
\beta>0),
\end{eqnarray}
with [cf.~(\ref{eqnmeasure})]
%
\begin{eqnarray}
\caH(\bflambda) : =\sum_{i=1}^N
\frac{1}{2} \biggl(U(\lambda_i)+\frac{\lambda_i^2}{2} \biggr)-
\frac{1}{N}\sum_{1\le i<j\le
N}\log\llvert
\lambda_j-\lambda_i\rrvert,
\end{eqnarray}
and with $Z_U\equiv Z_U(\beta)$ a normalization. For $K,L$ and $\bfy$
fixed, we can write~$\mu^{\bfy}$ as the Gibbs measure
%
\begin{eqnarray}
\label{legibsmeasures} \mu^{\bfy}(\dd\bfx)=\frac{1}{Z_U^{\bfy}}
\mathrm{e}^{-\beta N\caH
^{\bfy}(\bfx)}\,\dd\bfx,
\end{eqnarray}
where
%
\begin{eqnarray}
\caH^{\bfy}(\bfx)&=&\sum_{i\in I}
\frac{1}{2} V^{\bfy}(x_i)-\frac
{1}{N}\mathop{\sum_{i,j\in I}}_{i<j }\log\llvert
x_j-x_i\rrvert,
\end{eqnarray}
with
%
\begin{eqnarray}
V^{\bfy}(x)&=& U(x)+\frac{x^2}{2}-\frac{2}{N}\sum
_{i\notin I}\log\llvert x-y_i\rrvert\label{equationexternalpotential}
\end{eqnarray}
an external potential and with $Z_U^{\bfy}\equiv Z_U^{\bfy}(\beta)$
a normalization. Following~\cite{singlegap}, we next introduce the
notion of regular external potential:

\begin{definition}\label{definitionregularexternalpotential}
An external potential $V\equiv V^{\bfy}$ of a $\beta$-ensemble of $K$
points in a configuration interval $\mathbf{I}=(a,b)$ is called
$K^{\chi}$-regular if the following bounds hold:
%
\begin{eqnarray}
\llvert\mathbf{I}\rrvert&=&\frac{\caK}{N\rho(\bar y)}+\caO\biggl(\frac
{K^{\chi}}{N}
\biggr)\label{externalpotential1},
\\
V'(x)&=&\rho(\bar y)\log\frac{d_+(x)}{d_-(x)}+\caO\biggl(
\frac
{K^{\chi}}{Nd(x)} \biggr)\label{externalpotential2},
\\
V''(x)&\ge& 1+\inf U''(x)+
\frac{c}{d(x)}, \label{externalpotential3}
\end{eqnarray}
for $x\in\mathbf{I}$, with some $c>0$ and for some small $\chi
>0$, where
\begin{eqnarray*}
d(x): =\min\bigl\{\llvert x-a\rrvert,\llvert x-b\rrvert
\bigr\}
\end{eqnarray*}
denotes the distance to the boundary of $\mathbf{I}$,
\begin{eqnarray*}
d_-(x): =d(x)+\rho(\bar{y})N^{-1}K^{\chi}
\end{eqnarray*}
and
\begin{eqnarray*}
d_+(x): =\max\bigl\{\llvert x-a\rrvert,\llvert
x-b\rrvert\bigr
\}+\rho(\bar{y})N^{-1}K^{\chi}.
\end{eqnarray*}
\end{definition}

The main technical result we use in this section is Theorem~4.1
of~\cite{singlegap}; see Theorem~\ref{theorem31} below. It asserts
that the local gap statistics of $\mu^{\bfy}$ are essentially
independent of $\bfy$ and $U$, provided that $V^{\bfy}$ is $K^{\chi
}$-regular for some small $\chi>0$.

\subsection{Comparison of local measures}\label{subsectionlecomparisionoflocalmeasures}
Fix small $\alpha,\sigma>0$, and let $K$ and $L$ satisfy~(\ref
{conditiononK}). Recall that we denote by~$f_t \mu_G$ the
distribution of the eigenvalues of the matrix in~(\ref{thenewmatrix})
and by~$\widehat\psi_t \mu_G$ the reference $\beta
$-ensemble defined in~(\ref{definitionofpsit}). Following the
discussion in Section~\ref{subsection1ofsection6}, we introduce
the conditioned densities
%
\begin{eqnarray}
f_t^{\bfy} \mu_G^{\bfy}=(f_t
\mu_G)^{\bfy},\qquad\widehat\psi_t^{\bfy}
\mu_G^{\bfy}=(\widehat\psi_t
\mu_G)^{\bfy}.
\end{eqnarray}

Recall that we denote by $\widehat\rho_{\mathrm{fc}}(t)$ the equilibrium
density of $\widehat\psi_t \mu_G $ and by $\widehat\gamma_k\equiv
\widehat\gamma_k(t)$ the classical location of the $k$th particle
with respect to $\widehat\rho_{\mathrm{fc}}\equiv\widehat\rho_{\mathrm{fc}}(t)$;
cf.~(\ref{classicallocation}). Let $\varepsilon_0>0$ and define the set
of ``good'' boundary conditions, $\caR_{L,K}\equiv\caR
_{L,K}(\varepsilon_0,\alpha)$,
%
\begin{eqnarray}
\label{definitionofRsimple} \caR_{L,K}: &=&\bigl\{
{\bflambda}\in
\digamma^{(N)}\dvtx \llvert\lambda_{k}-\widehat
\gamma_{k}\rrvert\leq{N^{-1+\varepsilon_0}},\forall k\in\dllbracket
\alpha N,(1-\alpha) N\drrbracket\setminus I_{L,K}\bigr\}\hspace*{-25pt}
\nonumber\\[-8pt]\\[-8pt]\nonumber
&&{} \cap\bigl\{ {\bflambda}\in\digamma^{(N)}\dvtx \llvert
\lambda_{k}-\widehat\gamma_{k}\rrvert\le
N^{-2/3+\varepsilon_0},\forall k\in\llbracket1, N\rrbracket\bigr\}.
\end{eqnarray}
The next result compares the local statistics of $ f_t^{\bfy} \mu
_G^{\bfy}$ and $\widehat\psi_t^{\bfy} \mu_G^{\bfy}$ for $\bfy
\in\caR_{L,K}$. Recall that $\mathfrak{a}$ stands for any number in
$(0,1/2)$.

\begin{proposition}\label{frozenstatsproximity}
Fix small constants $\alpha,\sigma>0$ [see~(\ref{conditiononK})]
and $\varepsilon_0>0$; see~(\ref{definitionofRsimple}). Let $K$
satisfy~(\ref{conditiononK}), and let $O$ be an $n$-particle
observable. Let $\varepsilon'>0$, and choose $\tau$ satisfying $1\gg
\tau>N^{-2\mathfrak{a}}$. Then, for any $t\ge N^{\varepsilon'}\tau$
and any constant $\mathfrak{c}\in(0,1)$, there is a set of
configurations $\mathcal{G}\equiv\mathcal{G}_{L,K}(\varepsilon
_0,\alpha) \subset\caR_{L,K}(\varepsilon_0,\alpha)$, with
%
\begin{eqnarray}
\label{probabilityestimateongoodboundarycondition} \bbP^{f_t \mu
_G}(\mathcal{G}) \geq1 - \frac{N^{-\mathfrak
{c}}}{2},
\end{eqnarray}
such that
%
\begin{eqnarray}
\label{foralmostallboxes} \biggl\llvert\int O(\bfx) \bigl(f_t^{{\bfy}}(
\bfx) -{\widehat\psi_t^{\bfy}}(\bfx) \bigr)
\mu_G^{\bfy}(\dd\bfx) \biggr\rrvert\leq C_O
\sqrt{K} N^{\mathfrak{c}-\mathfrak{a}}\tau^{-1},
\end{eqnarray}
$t\ge N^{\varepsilon'}\tau$, for $N$ sufficiently large on~$\Omega$.
The constant~$C_O$, depends only on~$\varepsilon'$,~$\alpha$ and~$O$.

Moreover, there is $\upsilon>0$, such that
%
\begin{eqnarray}
\label{localftrigidity} \bbP^{f_t^{\bfy}\mu_G^{\bfy}} \bigl( \bigl\{
\bigl\llvert x_k
- \widehat\gamma_k(t) \bigr\rrvert< N^{-1+\varepsilon_0}, k\in
I_{L,K} \bigr\} \bigr) \geq1 - \mathrm{e}^{-\upsilon(\varphi_N)^{\xi}},
\end{eqnarray}
$t\ge N^{\varepsilon'}\tau$, for $N$ sufficiently large on~$\Omega$,
with $\xi=A_0\log\log N/2$; see~(\ref{eqxi}).
\end{proposition}

\begin{pf}
We follow\vspace*{1pt} closely the proof of Lemma~6.4 in~\cite{singlegap}. Let
$\tau$ satisfy $1\gg\tau> N^{-2\mathfrak{a}}$, and choose $t\ge
N^{\varepsilon'}\tau$. We estimate
%
\begin{eqnarray}
\label{frozenstatsbound} \biggl\llvert\int O(\bfx) \bigl(f_t^{{\bfy}}(
\bfx) -{\widehat\psi_t^{\bfy}}(\bfx) \bigr)
\mu_G^{{\bfy}}(\dd\bfx)\biggr\rrvert& \leq &C_O
\bigl\llVert f_t^{{\bfy}} \mu_G^{{\bfy}}
- {\widehat\psi_t^{\bfy}} \mu_G^{{\bfy}}
\bigr\rrVert_{1}
\nonumber\\[-8pt]\\[-8pt]\nonumber
&\leq& C_O\sqrt{S_{{\widehat\psi
_t^{\bfy}} \mu_G^{\bfy}}\bigl(\widehat
g_t^{{\bfy}}\bigr)},
\end{eqnarray}
where we use~(\ref{entropyinequality}) and set $\widehat g_t:=
f_t/{\widehat\psi_t}$. For $\bfy\in\caR_{L,K}$, we consider the
locally constrained measure ${\widehat\psi_t^{\bfy}}\mu_G^{\bfy}$,
explicitly given by
\begin{eqnarray*}
{\widehat\psi_t^{\bfy}}\mu_G^{\bfy}(
\dd{\bfx}) = \frac
{1}{Z^{\bfy}}\mathrm{e}^{-N\beta\widehat\caH^{\bfy}(t,\bfx)}\,\dd\bfx,
\end{eqnarray*}
with
\begin{eqnarray*}
\widehat\caH^{\bfy}(t,\bfx)&=&\sum_{k\in I}
\biggl(\frac{\widehat
U(t,x_k)}{2}+\frac{x_k^2}{4} \biggr)
\\
&&{}-\frac{1}{N} \mathop{\sum_{k,l\in I}}_{k<l}
\log\llvert x_k-x_l\rrvert-\frac
{1}{N}\mathop{\sum_{k\in I}}_{l\notin I}\log
\llvert x_k-y_l\rrvert .
\end{eqnarray*}
Here $I\equiv I_{L,K}$. From~(5.20) of~\cite{singlegap}, we know that
%
\begin{eqnarray}
\nabla^2_{\bfx} \widehat\caH^{\bfy}(t,\bfx) \geq
cN/K \qquad(\bfy\in\caR_{L,K}),
\end{eqnarray}
for some $c>0$ independent of $t$. Here, $\nabla^2_{\bfx}$ denotes
the Hessian with respect the variables $\bfx$. Thus the Bakry--\'
{E}mery criterion yields the logarithmic Sobolev inequality
%
\begin{eqnarray}
\label{frozenlog-sobolev} S_{{\widehat\psi_t^{\bfy}} \mu_G^{{\bfy
}}}\bigl(\widehat g_t^{{\bfy}}
\bigr) \leq\frac{CK}{N} D_{{\widehat\psi_t^{\bfy}} \mu_G^{{\bfy
}}} \Bigl(\sqrt{\widehat
g_t^{\bfy}} \Bigr)\qquad(\bfy\in\caR_{L,K}),
\end{eqnarray}
where the constant $C$ can be chosen independent $t$.

For $k\in\llbracket1,N\rrbracket$, denote by $D_{{\widehat\psi_t}
\mu_G, k}$ the Dirichlet form of the particle $k$, that is,
$D_{{\widehat\psi_t} \mu_G, k} (f) : =\frac
{1}{2N} \int\llvert \partial_k
f\rrvert ^2 {\widehat\psi_t} \mu_G$, and by $ D_{{\widehat\psi_t^{\bfy
}} \mu_G^{{\bfy}},k}$ its conditioned analogue (with $k\in
I_{L,K}$). Using the notation of~(\ref{thebflambda}), we may write
\begin{eqnarray*}
\E^{f_t \mu_G}D_{{\widehat\psi_t^{\bfy}} \mu_G^{{{\bfy
}}}} \Bigl(\sqrt{\widehat
g_t^{{{\bfy}}}} \Bigr)=\int D_{\widehat
\psi_t^{\bfy} \mu_G^{\bfy}} \Bigl(\sqrt
{\widehat g_t^{\bfy
}} \Bigr)f_t(
\bflambda) \mu_G(\dd\bflambda),
\end{eqnarray*}
and we can bound
%
\begin{eqnarray}
\label{lewhereweusethedformestimate} \E^{f_t \mu_G} D_{{\widehat\psi
_t^{\bfy}} \mu_G^{{{\bfy
}}}} \Bigl(\sqrt{\widehat
g_t^{{{\bfy}}}} \Bigr)& =&\mathbb{E}^{f_t\mu_G} \sum
_{k\in I} D_{{\widehat\psi_t^{\bfy}} \mu
_G^{{\bfy}},k} \Bigl(\sqrt{
\widehat g_t^{{\bfy}}} \Bigr)
\nonumber
\\
& \leq& D_{{\widehat\psi_t} \mu_G}(\sqrt{\widehat g_t})
\\
& \leq& CN^{1-2\mathfrak{a}} \tau^{-2},\nonumber
\end{eqnarray}
for $N$ sufficiently large, where we use Corollary~\ref
{lemmaboundonentroypandirichletform} in the last line. Thus Markov's inequality
implies, for $\mathfrak{c}>0$, that there exists a set of
configurations $\mathcal{G}^1\subset\mathcal{R}$, with $\mathbb
{P}^{f_t \mu_G}(\mathcal{G}^1)\ge1-N^{-\mathfrak{c}}$, such that,
for $\bfy\in\caG^1$,
%
\begin{eqnarray}
\label{goodintervalestimate} D_{{\widehat\psi_t^{\bfy}} \mu_G^{{\bfy}}}
\Bigl(\sqrt{\widehat
g_t^{{\bfy}}} \Bigr) \leq C N^{2\mathfrak{c}}
N^{1-2\mathfrak
{a}}\tau^{-2}
\end{eqnarray}
holds for $N$ sufficiently large on $\Omega$. Substituting~(\ref
{goodintervalestimate}) into~(\ref{frozenlog-sobolev}) and then into~(\ref{frozenstatsbound}), we find that
\begin{eqnarray*}
\biggl\llvert\int O(\bfx) \bigl(f_t^{{\bfy}} -{\widehat
\psi_t^{\bfy}} \bigr) \mu_G^{{\bfy}}(
\dd\bfx)\biggr\rrvert& \leq& C_O \sqrt{K} N^{\mathfrak{c}}N^{-\mathfrak{a}}
\tau^{-1},
\end{eqnarray*}
on $\Omega$ for $N$ sufficiently large. This proves~(\ref{foralmostallboxes}).

To prove~(\ref{localftrigidity}) note that the rigidity estimates
of Lemma~\ref{masterrigiditybound} imply
\begin{eqnarray*}
&& \bbE^{f_t \mu_G}  \bigl[ \bbP^{f_t^{\bfy} \mu_G^{\bfy}} \bigl( \bigl\{
\bigl\llvert
x_k - \widehat{\gamma}_k(t)\bigr\rrvert>
N^{-1+\varepsilon}, k\in I \bigr\} \bigr) \bigr]
\\
&&\qquad =
\bbP^{f_t \mu
_G} \bigl( \bigl\{\bigl\llvert x_k -
\widehat{\gamma}_k(t)\bigr\rrvert> N^{-1+\varepsilon}, k\in I \bigr
\} \bigr) \le\mathrm{e}^{-\upsilon(\varphi
_N)^{\xi}},
\end{eqnarray*}
for\vspace*{1pt} some $\upsilon>0$, where we have chosen $\xi=A_0\log\log N/2$.
By Markov's inequality we conclude that there is a set of
configurations,~$\caG^{2}$, such that~(\ref{localftrigidity})
holds with $(\xi,\upsilon)$-high probability. Finally, set $\caG
: =
\caG^1\cap\caG^2$, and note that $\caG$ satisfies~(\ref
{probabilityestimateongoodboundarycondition}).
\end{pf}

\subsection{Gap universality for local measures}\label{gapuniversalityforlocalmeasures}
In Section~\ref{subsectionlecomparisionoflocalmeasures}, we show
that the local gap statistics of the measure $f_t^{\bfy} \mu_G^{\bfy
}$ agree with those of $\widehat\psi_t^{\bfy} \mu_G^{\bfy}$ for
boundary conditions $\bfy$ in the set $\caR_{L,K}$. In this
subsection, we are going to show that the local statistics of $\widehat
\psi_t^{\bfy} \mu_G^{\bfy}$ are essentially independent of the
precise form of $\bfy$, as is asserted by the main theorem of~\cite
{singlegap}. Recall the notion of external potential introduced
in~(\ref{equationexternalpotential}).

%
\begin{theorem}[(Gap universality for local measures, Theorem~4.1
in~\cite{singlegap})]\label{theorem31}
Let $L,\widetilde L$ and $\caK=2K+1$ satisfy~(\ref{conditiononK}) with
$\alpha,\sigma>0$. Consider two boundary conditions $\bfy$, $\tilde
\bfy$ such that the configuration intervals coincide, that is,
%
\begin{eqnarray}
\label{assumption1thm31} \bsI=(y_{L-K-1},y_{L+K+1})=(\tilde
{y}_{\widetilde L-K-1},
\tilde{y}_{\widetilde L+K+1}).
\end{eqnarray}
Consider two measures $\mu$ and $\tilde\mu$ in the form of~(\ref
{betaensemblesection6}), with possibly two different potentials $U$
and $\widetilde U$, and consider the constrained measures $\mu^{\bfy}$
and $\tilde\mu^{\tilde\bfy}$. Let $\chi>0$, and assume that the
external potentials $V^{\bfy}$ and $\tilde V^{\tilde\bfy}$
[see~(\ref{equationexternalpotential})] are $K^{\chi}$-regular; see
Definition~\ref{definitionregularexternalpotential}. In particular,
assume that $\bsI$ satisfies
%
\begin{eqnarray}
\label{assumption2thm31} \llvert\bsI\rrvert=\frac{\caK}{N\varrho
_U(\bar y)}+\caO\biggl(
\frac{K^{\chi
}}{N} \biggr).
\end{eqnarray}
Assume further that
%
\begin{eqnarray}
\label{assumption3thm31} \max_{j\in I_{L,K}}\bigl\llvert\E^{\mu^{\bfy}}
x_j-\alpha_j\bigr\rrvert+\max_{j\in I_{
\tilde L,K}}
\bigl\llvert\E^{\tilde\mu^{\tilde{\bfy}}} x_j-\alpha_j\bigr
\rrvert\le CN^{-1}K^{\chi}.
\end{eqnarray}
Let $p\in\Z$ satisfy $\llvert p\rrvert \le K-K^{1-\chi'}$, for some
small $\chi
'>0$. Fix $n\in\N$. Then there is a constant $\chi_0$, such that if
$\chi,\chi'<\chi_0$, then for any $n$-particle observable $O$, we have
\begin{eqnarray*}
&& \bigl\llvert\mathbb{E}^{\mu^{\bfy}}O \bigl(N(x_{L+p+1}-x_{L+p}),
\ldots, N(x_{L+p+n}-x_{L+p}) \bigr)
\\
&&\qquad{} -\mathbb{E}^{\tilde\mu
^{\tilde\bfy}}O
\bigl(N(x_{L+p+1}-x_{L+p}),\ldots, N(x_{L+p+n}-x_{L+p})
\bigr) \bigr\rrvert\le C_O K^{-\mathfrak{b}},
\end{eqnarray*}
for some constant $\mathfrak{b}>0$ depending on $\sigma$, $\alpha$,
and for some constant $C_O$ depending on $O$. This holds for $N$
sufficiently large [depending on the $\chi,\chi',\alpha$ and $C$
in~(\ref{assumption3thm31})].
\end{theorem}

Recall that the measure $\widehat\psi_t^{\bfy} \mu_G^{\bfy}$ can
be written as the Gibbs measure
%
\begin{eqnarray}
{\widehat\psi_t^{\bfy}}\mu_G^{\bfy}(
\dd{\bfx}) = \frac
{1}{Z_{\widehat\psi_t}^{\bfy}}\mathrm{e}^{-N\beta\caH^{\bfy
}(t,\bfx)}\,\dd\bfx,
\end{eqnarray}
where
%
\begin{eqnarray}
\caH^{\bfy}(t,\bfx)&=&\sum_{i\in I}
\frac{1}{2} V^{\bfy
}(t,x_i)-\frac{1}{N}
\mathop{\sum_{i,j\in I}}_{i<j }\log\llvert
x_j-x_i\rrvert,
\end{eqnarray}
with the external potential
%
\begin{eqnarray}
\label{lepotentialexterior} V^{\bfy}(t,x)&=& \widehat U(t,x)+\frac{x^2}{2}-
\frac{2}{N}\sum_{i\notin I}\log\llvert
x-y_i\rrvert.
\end{eqnarray}
Using Theorem~\ref{theorem31} we first show that the local
statistics of ${\psi_t^{\bfy}}\mu_G^{\bfy}$ are virtually
independent of $\bfy$; that is, we apply Theorem~\ref{theorem31}
with $\mu^{\bfy}=(\widehat\psi_t \mu_G)^{\bfy}$ and $\tilde\mu
^{\tilde\bfy}=(\widehat\psi_t \mu_G)^{\tilde\bfy}$.

We first check the regularity assumption of the external
potential~$V^{\bfy}$. Recall the definition of $K^{\chi}$-regular
potential in Definition~\ref{definitionregularexternalpotential}.

\begin{lemma}\label{potentialisregular}
Fix small constants $\alpha,\sigma>0$; see~(\ref{conditiononK}).
Let $\chi>0$, and consider $\bfy\in\caR_{L,K}(\chi\sigma/2,\alpha
/2)$. Then, on the event~$\Omega$, the external potential~$V^{\bfy
}(t,x)$ in~(\ref{lepotentialexterior}) is $K^{\chi}$-regular on
$\mathbf{I}=(y_{L-K-1},y_{L+K+1})$.
\end{lemma}

The proof of Lemma~\ref{potentialisregular} follows almost verbatim
the proof of Lemma~4.5 in the Appendix~A of~\cite{singlegap}, and we
therefore omit it here.\vadjust{\goodbreak}

To check that assumption~(\ref{assumption3thm31}) of Theorem~\ref
{theorem31} holds, we use the following result. Recall the set of
configurations $\caG$ of Proposition~\ref{frozenstatsproximity}.

\begin{lemma}\label{frozenbetaparticleproximity}
Under the assumptions of Proposition~\ref{frozenstatsproximity} the
following holds. Let $\bfy\in\caG$. Then, for all $k\in I_{L,K}$,
%
\begin{eqnarray}
\bigl\llvert\bbE^{f_t^{\bfy}\mu_G^{\bfy}} x_k - \bbE^{{\widehat\psi
_t^{\bfy}} \mu_G^{\bfy}}
x_k \bigr\rrvert&\leq& C \frac{KN^{2\mathfrak
{c}}}{N} N^{-\mathfrak{a}}
\tau^{-1}\qquad\bigl(t\ge\tau N^{\varepsilon
'}\bigr),\label{frozenbetaparticleproximityequation}
\end{eqnarray}
for $N$ sufficiently large on $\Omega$.
\end{lemma}
\begin{pf}
We follow the proof of Lemma~6.5 of~\cite{singlegap}. Fix $t\ge\tau
N^{\varepsilon'}$, where $1\gg\tau\ge N^{-2\mathfrak{a}}$. Let $\bfy
\in\caG$. Denote by $\mathcal{L}_t^{\bfy}$ the generator associated
to the Dirichlet form $D_{{\widehat\psi_t^{\bfy}} \mu_G^{\bfy}}$,
that is,
\begin{eqnarray*}
\int f \mathcal{L}_t^{\bfy}g {\widehat
\psi_t^{\bfy}}\,\dd\mu_G^{\bfy} = -
\frac{1}{\beta N} \sum_{i \in I} \int
\partial_i f \partial_i g {\widehat
\psi_t^{\bfy}}\,\dd\mu_G^{\bfy}
\qquad(I\equiv I_{L,K}).
\end{eqnarray*}
Let $q_s$ be the solution of the evolution equation $\partial_s q_s =
\mathcal{L}_t^{\bfy} q_s$, $s\ge0$, with initial condition $q_0
: =
\widehat g_t^{\bfy} = f_t^{\bfy} / {\widehat\psi_t^{\bfy}}$. Note
that\vspace*{1pt} $q_s$ is a density with respect to the reversible measure,
${\widehat\psi_t^{\bfy}} \mu_G^{\bfy}$, of this dynamics. Hence,
we can write
\begin{eqnarray*}
\bigl\llvert\bbE^{f_t^{\bfy}\mu_G^{\bfy}} x_k - \bbE^{{\widehat\psi
_t^{\bfy}} \mu_G^{\bfy}}
x_k \bigr\rrvert&=& \biggl\llvert\int_0^\infty
\dd s \int x_k \mathcal{L}_t^{\bfy}
q_s {\widehat\psi_t^{\bfy}}\,\dd
\mu_G^{\bfy} \biggr\rrvert
\\
& =& \biggl\llvert\frac{1}{\beta N} \int_0^{\infty}
\dd s \int\partial_k q_s {\widehat
\psi_t^{\bfy
}}\,\dd\mu_G^{\bfy}
\biggr\rrvert.
\end{eqnarray*}
Recall that ${\widehat\psi_t^{\bfy}}\mu^{\bfy}_G$ satisfies the
logarithmic Sobolev inequality~(\ref{frozenlog-sobolev}) with
constant $\tau_K: =CK/N$, provided that $\bfy\in
\caR_{L,K}$.
Thus, upon using Cauchy--Schwarz and the exponential decay of the
Dirichlet form $D_{{\widehat\psi_t^{\bfy}} \mu_G^{\bfy}}(\sqrt
{q_s})$, we obtain for some~$\upsilon',c> 0$,
\begin{eqnarray*}
\bigl\llvert\bbE^{f_t^{\bfy}\mu_G^{\bfy}} x_k - \bbE^{{\widehat\psi
_t^{\bfy}} \mu_G^{\bfy}}
x_k \bigr\rrvert= \biggl\llvert\frac{1}{\beta N} \int
_0^{N^{\upsilon'} \tau_K} \dd s \int\partial_k
q_s {\widehat\psi_t^{\bfy}}\,\dd
\mu_G^{\bfy} \biggr\rrvert+\caO\bigl(\mathrm
{e}^{-cN^{\upsilon'}} \bigr).
\end{eqnarray*}
Using
\begin{eqnarray*}
\llvert\partial_k q_s\rrvert= 2\llvert
\sqrt{q_s} \partial_k\sqrt{q_s}\rrvert
\leq R (\partial_k \sqrt{q_s})^2 +
R^{-1} q_s,
\end{eqnarray*}
where $R>0$ is a free parameter, we obtain
\begin{eqnarray*}
&& \biggl\llvert\frac{1}{\beta N} \int_0^{N^{\upsilon'}\,\dd s \tau_K}
\int\partial_k q_s {\widehat\psi_t^{\bfy}}
\,\dd\mu_G^{\bfy} \biggr\rrvert
\\
&&\qquad \leq R \biggl[\int
_0^{N^{\upsilon'} \tau_K} \,\dd s\, D_{f_t^{\bfy}\mu
_G}(
\sqrt{q_s}) \biggr]
 + \frac
{1}{2}R^{-1}N^{-1+\upsilon'}\tau_K
\\
&&\qquad  \leq R S_{{\widehat\psi_t^{\bfy}} \mu_G^{\bfy}}\bigl(\widehat
g_t^{\bfy}\bigr)
+ \frac{1}{2}R^{-1}N^{-1+\upsilon'}\tau_K
\\
&&\qquad  \leq C R\tau_K D_{{\widehat\psi_t^{\bfy}} \mu_G^{\bfy}}\Bigl(\sqrt
{\widehat
g_t^{\bfy}}\Bigr)+ \frac{1}{2}R^{-1}N^{-1+\upsilon'}
\tau_K,
\end{eqnarray*}
where in the second line we use that the time integral of the Dirichlet
form is bounded by the initial entropy (see, e.g.,~Theorem~2.3 in~\cite
{EY}) and in the final line we used the logarithmic Sobolev
inequality~(\ref{frozenlog-sobolev}). Optimizing over~$R$, we get
\begin{eqnarray*}
\bigl\llvert\bbE^{f_t^{\bfy}\mu_G^{\bfy}} x_k - \bbE^{{\widehat\psi
_t^{\bfy}} \mu_G^{\bfy}}
x_k \bigr\rrvert&\leq& C \tau_K \Bigl(N^{-1+\upsilon'}D_{ {\widehat\psi
_t^{\bfy}} \mu_G^{\bfy}}
\Bigl(\sqrt{\widehat g_t^{\bfy}} \Bigr)
\Bigr)^{1/2} + \caO\bigl(\mathrm{e}^{-c N^{\upsilon'}} \bigr)
\nonumber
\\
&\le& C \frac{KN^{\upsilon'/2}}{N^{3/2}} \Bigl(D_{ {\widehat\psi
_t^{\bfy}} \mu_G^{\bfy}} \Bigl(\sqrt{\widehat
g_t^{\bfy}} \Bigr) \Bigr)^{1/2} + \caO\bigl(
\mathrm{e}^{-c N^{\upsilon'}} \bigr),
\end{eqnarray*}
where we used that $\tau_K= CK/N$.
Using (\ref{goodintervalestimate}) we finally obtain
\begin{eqnarray*}
\bigl\llvert\bbE^{f_t^{\bfy}\mu_G^{\bfy}} x_k - \bbE^{{\widehat\psi
_t^{\bfy}} \mu_G^{\bfy}}
x_k \bigr\rrvert\leq C \frac{KN^{2\mathfrak
{c}}}{N}N^{-\mathfrak{a}}
\tau^{-1} +\caO\bigl(\mathrm{e}^{-c
N^{\upsilon'}} \bigr),
\end{eqnarray*}
for $N$ sufficiently large on $\Omega$.
\end{pf}

%
\begin{lemma}\label{lemmafrozengapuniversality}
Fix small constants $\alpha,\sigma>0$. Fix $\varepsilon'>0$ and $t\ge
\tau N^{\varepsilon'}$, where $\tau$ satisfies $1\gg\tau\ge
N^{-2\mathfrak{a}}$. Fix $n\in\N$, and consider an $n$-particle
observable $O$. Let $\chi',\chi>0$, with $\chi',\chi<\chi_0$,
where $\chi_0$ is the constant in Theorem~\ref{theorem31}. Then the
following holds.

Assume that $0<\mathfrak{a}<1/2$, $0<\mathfrak{c}<1$, $N^{-2\mathfrak
{a}}\le\tau\ll1$ and $K\in\llbracket N^{\sigma},N^{1/4}\rrbracket
$ are chosen such that
%
\begin{eqnarray}
\label{choiceoftheconstants} \frac{KN^{2\mathfrak{c}}}{N} N^{-\mathfrak
{a}}\tau^{-1} \le
\frac
{K^{\chi}}{N}.
\end{eqnarray}
Let $p$ be an integer satisfying $\llvert p\rrvert \le K-K^{1-\chi
'}$. Let $\bfy\in
\mathcal{G}_{L,K}(\frac{\chi^2 \sigma}{2},\frac{\alpha}{2})$. Then,
for the observable $G$, as defined in~(\ref{theaveragedobservable}),
we have
%
\begin{equation}
\label{equationfrozengapuniversality} \qquad\biggl\llvert\int G_{L+p,n}(\bfx
) \bigl(f_t^{\bfy}
\,\dd\mu_G^{\bfy}-{\widehat\psi_t} \,\dd
\mu_G\bigr) \biggr\rrvert\leq C_OK^{-\mathfrak{b}}+C_O
\sqrt{K} N^{\mathfrak{c}}N^{-\mathfrak{a}}\tau^{-1},
\end{equation}
for $N$ sufficiently large on $\Omega$, where the constant $C_O$
depends on $O$ and $\varepsilon'$, and the constant $\mathfrak{b}>0$
depends on~$\alpha$ and~$\sigma$.
\end{lemma}

\begin{pf} We follow~\cite{singlegap}. Fix $t\ge\tau N^{\varepsilon'}$
and $\chi>0$. Let $\bfy\in\mathcal{G}_{L,K}(\frac{\chi^2 \sigma
}{2}, \alpha) \subset\mathcal{G}_{L,K}( \frac{\chi\sigma}{2},
\alpha)$. Then by Proposition~\ref{frozenstatsproximity} and the
assumption in~(\ref{choiceoftheconstants}),
\begin{eqnarray*}
\bigl\llvert\bbE^{f_t^{\bfy}\mu_G^{\bfy}} x_k - \widehat{
\gamma}_k(t) \bigr\rrvert\le C K^{\chi}N^{-1},
\end{eqnarray*}
for all $k\in I\equiv I_{L,K}$. Further, from Lemma~\ref
{frozenbetaparticleproximity} and the assumption in~(\ref
{choiceoftheconstants}) we get
%
\begin{eqnarray}
\label{superacc} \bigl\llvert\bbE^{{\widehat\psi_t^{\bfy}} \mu_G^{\bfy
}} x_k - \widehat{
\gamma}_k(t) \bigr\rrvert\leq C K^{\chi}N^{-1},
\end{eqnarray}
for all $k\in I$. Recall from~(\ref{equidistantpoints1}) that we
denote by $ \bar{y}: =\frac
{1}{2}(y_{L-K-1}+y_{L+K+1})$ the
midpoint of the configuration interval~$\bsI$ and that $(\alpha_k)$
denote $2K+1$ equidistant points in~$\bsI$. As shown in Lemmas~4.5 and~5.2 of~\cite{singlegap}, we have
\begin{eqnarray*}
\bigl\llvert\widehat{\gamma}_k(t) - \alpha_k \bigr
\rrvert\leq CK^{\chi} N^{-1},
\end{eqnarray*}
for all $k\in I$, provided that $\bfy\in\caG_{L,K}(\frac{\chi
\sigma}{2}, \alpha)$. We hence obtain
%
\begin{eqnarray}
\label{improvedacc1} \bigl\llvert\bbE^{{\widehat\psi_t^{\bfy}} \mu
_G^{\bfy}} x_k - \alpha
_k \bigr\rrvert\leq C K^{\chi}N^{-1},
\end{eqnarray}
for $N$ sufficiently large on $\Omega$.

Proposition~\ref{frozenstatsproximity} implies that there is $C_O$
such that
%
\begin{eqnarray}
\label{whatwegetfromtheproposition}
\biggl\llvert\int G_{L+p,n}(\bfx)
\bigl(
f_t^{\bfy}\,\dd\mu_G^{\bfy} - {
\widehat\psi_t^{\bfy}}\,\dd\mu_G^{\bfy}
\bigr) \biggr\rrvert\leq C_O \sqrt{K} N^{\mathfrak{c}}N^{-\mathfrak{a}}
\tau^{-1}
\nonumber\\[-8pt]\\[-10pt]
\eqntext{\bigl(t\ge N^{\varepsilon'}\tau\bigr),}
\end{eqnarray}
for $\bfy\in\caG_{L,K}(\frac{\chi\sigma}{2}, \alpha)$, $N$
sufficiently large on $\Omega$.

For $\alpha,\varepsilon_0,\varsigma_1>0$ and a $\beta$-ensemble $\mu$
on $\digamma^{(N)}$, define a set of particle configurations $\mathcal
{R}_{\mu}^{*}\equiv\mathcal{R}_{\mu}^{*}(\varepsilon_0,\alpha) $ by
\begin{eqnarray*}
\mathcal{R}^{*}_{\mu} : = \bigl\{\bfy\in
\digamma^{(N-\caK)}: \bbP^{\mu^{\bfy}} \bigl(\llvert x_k -
\gamma_k\rrvert> N^{-1+\varepsilon_0} \bigr) \le\mathrm{e}^{-\sklfrac{1}{2} N^{\varsigma
_1}},\forall k\in I_{L,K} \bigr\},
\end{eqnarray*}
where $\gamma_k$ denotes the classical location of the $k$th particle
with respect to the equilibrium measure of $\mu$.

As in the proof of Proposition~\ref{frozenstatsproximity}, it
follows from Markov's inequality and the rigidity bound for the $\beta
$-ensemble ${\widehat\psi_t} \mu_G$ in Lemma~\ref
{rigidityfortimedependentbeta} that we can choose $\mathcal
{R}^{*}_{{\widehat
\psi_t} \mu_G}\subset\caR_{L,K}$ and that $\bbP^{{\widehat\psi
_t} \mu_G}(\mathcal{R}^{*}_{{\widehat\psi_t} \mu_G})\geq1-
c\mathrm{e}^{-\sklfrac{1}{2} N^{\varsigma_1}}$, for some $c>0$, possibly
after decreasing $\varsigma_1$ by a small amount. For ${\tilde\bfy}
\in\mathcal{R}_{{\widehat\psi_t} \mu_G}^{*}(\frac{\chi^2
\sigma}{2}, \frac{\alpha}{2})$, Lemma 5.1 of~\cite{singlegap}
implies that
%
%
\begin{eqnarray}
\label{improvedacc2} \bigl\llvert\bbE^{{\widehat\psi_t^{\tilde\bfy}}
\mu_G^{\tilde\bfy
}}x_k-
\alpha_k \bigr\rrvert\leq CK^{\chi} N^{-1},
\end{eqnarray}
for $N$ sufficiently large on $\Omega$. Thus together with~(\ref
{improvedacc1}), we have on $\Omega$
%
\begin{eqnarray}
\label{tochecktheassumptionintheorem} \bigl\llvert\bbE^{{\widehat\psi
_t^{\tilde\bfy}} \mu_G^{\tilde\bfy
}}x_k-
\alpha_k \bigr\rrvert+\bigl\llvert\bbE^{{\widehat\psi_t^{\bfy}} \mu
_G^{\bfy}}
x_k - \alpha_k \bigr\rrvert\leq C K^{\chi}N^{-1},
\end{eqnarray}
for $N$ sufficiently large, for all $\bfy\in\caG(\frac{\chi\sigma
}{2}, \alpha)$ and all ${\tilde\bfy} \in\mathcal{R}_{{\widehat
\psi_t} \mu_G}^{*}(\frac{\chi^2 \sigma}{2}, \frac{\alpha}{2})$.\vspace*{1pt}

We now apply Theorem~\ref{theorem31}: let $\tilde\bfy$ and $\bfy$
be as above. By\vspace*{1pt} the scaling argument of Lemma~5.3 in~\cite
{singlegap}, we can assume that the two configuration intervals
$\tilde{\mathbf{I}}$ and $\mathbf{I}$ agree, so that
assumption~(\ref{assumption1thm31}) of Theorem~\ref{theorem31}
holds. Moreover, by Lemma~\ref{potentialisregular} we know that
$V^\bfy$ and $V^{\tilde\bfy}$ are $K^{\chi}$-regular external
potentials. The assumption in~(\ref{assumption3thm31}) of
Theorem~\ref{theorem31} is satisfied by~(\ref
{tochecktheassumptionintheorem}). Thus Theorem~\ref{theorem31} implies that
there is $\mathfrak{b}>0$, depending on~$\sigma$ and~$\alpha$, such that
%
\begin{eqnarray}
\label{whatwegetfrom31} \biggl\llvert\int G_{L+p,n}(\bfx) \bigl
({\widehat
\psi_t^{\bfy}} \,\dd\mu_G^{\bfy} - {
\widehat\psi_t^{\tilde\bfy}}\,\dd\mu_G^{\tilde\bfy
}
\bigr)\biggr\rrvert\leq C_O K^{-\mathfrak{b}},
\end{eqnarray}
for $N$ sufficiently large on $\Omega$. Since estimate~(\ref
{whatwegetfrom31}) holds for all $\tilde\bfy\in\mathcal{R}_{{\widehat
\psi_t} \mu_G}^{*}$, and since $\bbP^{{\widehat\psi_t} \mu
_G}(\mathcal{R}_{{\widehat\psi_t} \mu_G}^{*})\geq1- \mathrm
{e}^{-\sklfrac{1}{2} N^{\varsigma_1}}$, we can integrate over $\tilde
\bfy
$ to find that
\begin{eqnarray*}
\biggl\llvert\int G_{L+p,n}(\bfx) \bigl({\widehat\psi_t^{\bfy}}
\,\dd\mu_G^{\bfy} - {\widehat\psi_t} \,\dd
\mu_G \bigr)\biggr\rrvert\leq C_O K^{-\mathfrak{b}},
\end{eqnarray*}
for $N$ sufficiently large on $\Omega$. In combination with~(\ref
{whatwegetfromtheproposition}), this yields~(\ref
{equationfrozengapuniversality}).
\end{pf}

\subsection{Proof of Theorems~\texorpdfstring{\protect\ref{maintheoremofsection6}}{6.1} and~\texorpdfstring{\protect\ref{corollarygaussian}}{6.2}}\label{subsectionproofoftheoremandcorollarysection6}
Lemma~\ref{lemmafrozengapuniversality} compares the local
statistics of the locally-constrained measure $f_t^{\bfy} \mu
_G^{\bfy}$ with the $\beta$-ensemble ${\widehat\psi_t} \mu_G$. In
order to compare with local statistics of the measure $f_t \mu_G$
with $\widehat\psi_t \mu_G$, we next integrate out the boundary
conditions $\bfy$.

%
\begin{lemma}\label{lemmafinalestimate}
Under the assumptions of Lemma~\ref{lemmafrozengapuniversality} the
following holds. Let $J\subset\llbracket\alpha N,(1-\alpha
)N\rrbracket$ be an interval of consecutive integers in the bulk. Then
%
\begin{eqnarray}
\label{finalestimate}
&& \biggl\llvert\int\frac{1}{\llvert J\rrvert
}\sum _{j\in J}G_{j,n}(\bfx) (f_t\,\dd\mu
_G-{\widehat\psi_t}\,\dd\mu_G)\biggr
\rrvert
\nonumber\\[-8pt]\\[-8pt]\nonumber
&&\qquad \le C_O \bigl(N^{-\mathfrak{c}}+ K^{-\mathfrak{b}}+K^{-\chi'/2}
\bigr)
+C_O \sqrt{K} N^{\mathfrak{c}}N^{-\mathfrak{a}}
\tau^{-1},
\end{eqnarray}
for $N$ sufficiently large on $\Omega$.
\end{lemma}

\begin{pf}
For a small $\chi'>0$ as in Lemma~\ref{lemmafrozengapuniversality}, set
${\tilde K}: =K-K^{1-\chi'/2}$.
We first assume
that $J$ is such that $\llvert J\rrvert \le2\tilde K+1$. We then
choose $L$ such
that $J\subset I_{L,\tilde K}\subset I_{L,K}$. Recall the set of
configurations $\caG$ in Proposition~\ref{frozenstatsproximity}.
Using the conditioned measure $f_t^{\bfy} \mu_G^{\bfy}$ we estimate
%
\begin{eqnarray}
&& \E^{f_t \mu_G} \biggl[\frac{1}{\llvert J\rrvert }\sum_{j\in J}G_{j,n}
\biggr]
\nonumber\\[-8pt]\\[-8pt]\nonumber
&&\qquad =\E^{f_t \mu_G} \biggl[\frac{1}{\llvert J\rrvert }\int\sum
_{j\in J}G_{j,n}(\bfx)f_t^{\bfy}
\,\dd\mu_G^{\bfy}\lone(\caG) \biggr]+\caO
\bigl(N^{-\mathfrak{c}} \bigr),
\end{eqnarray}
where we used~(\ref{probabilityestimateongoodboundarycondition}).
Next, using Lemma~\ref{lemmafrozengapuniversality} we obtain on
$\Omega$
\begin{eqnarray*}
&& \frac{1}{\llvert J\rrvert }\int\sum_{j\in J}G_{j,n}(
\bfx)f_t^{\bfy}(\bfx)\,\dd\mu_G^{\bfy}(
\bfx)\lone (\caG)
\\
&&\qquad =\frac{1}{\llvert J\rrvert }\int\sum_{j\in J}
G_{j,n} \widehat\psi_t\,\dd\mu_G
 +\caO\bigl( K^{-\mathfrak{b}} \bigr)+\caO\bigl( \sqrt{K}
N^{\mathfrak
{c}}N^{-\mathfrak{a}}\tau^{-1} \bigr),
\end{eqnarray*}
on $\Omega$. For the special case $\llvert J\rrvert \le2\tilde K+1$, this
yields~(\ref{finalestimate}).

If $\llvert J\rrvert \ge\tilde K+1$, there are $L_a\in\llbracket
\alpha N,(1-\alpha
)N\rrbracket$, with $a\in\llbracket1,M_0\rrbracket$, such that the
intervals $I_{L_a,K}=\llbracket L_a-K,L_a+K\rrbracket$ are
nonintersecting with the properties that $J\subset\bigcup_{a=1}^{M_0}
I_{L_a,K}$ and $J\cap I_{L_{a},K}\neq\varnothing$, for all
$a\in\llbracket1,M_0\rrbracket$. Note that $M_0\le\frac
{\llvert J\rrvert }{K}+2$. For simplicity of notation we abbreviate
$I^{(a)}\equiv
I_{L_a,K}=\llbracket L_a-K,L_a+K\rrbracket$ and $\tilde I^{(a)}\equiv
\llbracket L_a-\tilde K,L_a+\tilde K\rrbracket$. We also label the
interior and exterior points of a configuration $\bflambda\in\digamma
^{(N)}$ accordingly,
\begin{eqnarray*}
\bfx^{(a)}=(x_{L_a-K},\ldots, x_{K+L_a})\in
\digamma^{(\caK)},
\end{eqnarray*}
respectively,
\begin{eqnarray*}
\bfy^{(a)}=(y_1,\ldots, y_{L_a-K-1},y_{L_a+K+1},
\ldots, y_N)\in\digamma^{(N-\caK)};
\end{eqnarray*}
cf. (\ref{bfxnotation}). We let $\caG^{(a)}\equiv\caG
_{L_a,K}(\varepsilon_0,\alpha)\subset\caR_{L_a,K}(\varepsilon_0,\alpha
)$ denote the set of configurations obtained in Proposition~\ref
{frozenstatsproximity}. Using this notation we can write
%
\begin{eqnarray}\label{eq635}
&& \E^{f_t \mu_G} \biggl[\frac{1}{\llvert J\rrvert }\sum
_{j\in J}G_{j,n} \biggr]\nonumber
\\
&&\qquad =
\frac{1}{\llvert J\rrvert }\sum_{{ a\dvtx I^{(a)}\cap J\neq
\varnothing}} \E^{f_t \mu_G}
\biggl[\int\sum_{j\in I^{(a)}\cap
J}G_{j,n}\bigl(
\bfx^{(a)}\bigr)f_t^{\bfy^{(a)}}\,\dd
\mu_G^{\bfy^{(a)}}\lone\bigl(\caG^{(a)}\bigr) \biggr]
\\
&&\quad\qquad{} +\caO\bigl( N^{-\mathfrak{c}} \bigr),\nonumber
\end{eqnarray}
on $\Omega$, where the first summation on the right-hand side is over
indices $a\in\llbracket1,M_0\rrbracket$ such that the intervals
$(I^{(a)})$ satisfy $I^{(a)}\cap J\neq\varnothing$. Here, we also use
the probability estimate on $\caG^{(a)}$ in~(\ref
{probabilityestimateongoodboundarycondition}). In~(\ref{eq635}) we may
further restrict, for each $a$, the summation over the index $j$ from
$I^{(a)}$ to $\tilde I^{(a)}$ at an expense of an error term of order $
\llvert I^{(a)}\setminus\tilde I^{(a)}\rrvert \le K^{1-\chi'/2}$.
Then summing
over $a\in\llbracket1,M_0\rrbracket$, with $M_0\sim\llvert J\rrvert
/K$, we get
%
\begin{eqnarray}
\label{onthewaytothefinalestimate2}
\qquad && \E^{f_t \mu_G} \biggl[\frac
{1}{\llvert J\rrvert }\sum
_{j\in J}G_{j,n} \biggr]\nonumber
\\
&&\qquad =
\frac{1}{\llvert J\rrvert }\sum_{{ a\dvtx I^{(a)}\cap J\neq
\varnothing}} \E^{f_t \mu_G}
\biggl[\int\sum_{j\in\tilde
I^{(a)}\cap J}G_{j,n}\bigl(
\bfx^{(a)}\bigr)f_t^{\bfy^{(a)}}\,\dd
\mu_G^{\bfy
^{(a)}}\lone\bigl(\caG^{(a)}\bigr) \biggr]
\\
&&\qquad\quad{} +\caO\bigl( N^{-\mathfrak{c}} \bigr)+\caO\bigl( K^{-\chi
'/2}\bigr),\nonumber
\end{eqnarray}
on\vspace*{1pt} $\Omega$. Since for each choice of the index $a$ the term in the
expectation on the right-hand side of~(\ref
{onthewaytothefinalestimate2}) can be dealt with as in the case
$\llvert J\rrvert \le2\tilde K+1$
above, this completes the proof of~(\ref{finalestimate}) for general $J$.
\end{pf}

We can now give the proof of Theorem~\ref{maintheoremofsection6}.

\begin{pf*}{Proof of Theorem~\ref{maintheoremofsection6}}
Let $\alpha>0$. We first choose the constants $\mathfrak{a}\in
(0,1/2)$, $\mathfrak{c}\in(0,1)$ and $\varepsilon'>0$, and the
parameter $K\in\llbracket N^{\sigma},N^{1/4}\rrbracket$
appropriately:\vspace*{1pt} let $\delta>0$ be a small constant. Then we set
$\mathfrak{a}\equiv1/2-\delta$, $\mathfrak{c}\equiv\delta/4$,
$K\equiv N^{\delta/4}$, $\varepsilon'\equiv\delta$, $\sigma=\delta
/8$. Note first that for this choice of $K$ condition~(\ref
{conditiononK}) is satisfied. Second, for sufficiently small $\delta
>0$, we
observe that
\begin{eqnarray*}
KN^{2\mathfrak{c}} N^{-\mathfrak{a}}\tau^{-1}=N^{{3\delta
}/{4}}N^{-\mathfrak{a}}
\tau^{-1}\le{K^{\chi}},
\end{eqnarray*}
holds, for example, for $\tau\ge N^{\delta} N^{-\mathfrak{a}}$ and
$\chi>0$ (with $\chi<\chi_0$). Thus~(\ref{choiceoftheconstants})
is satisfied with the above choices.

Hence, for $t\ge N^{2\delta}\tau$, Lemma~\ref{lemmafinalestimate}
yields, for some $\mathfrak{b}>0$,
\begin{eqnarray*}
&& \biggl\llvert\int\frac{1}{\llvert J\rrvert }\sum_{j\in J}G_{j,n}(
\bfx) (f_t\,\dd\mu_G-{\widehat\psi_t}\,\dd
\mu_G)\biggr\rrvert
\\
&&\qquad \le C_O K^{-\mathfrak{b}}+C_ON^{-\mathfrak
{c}}+C_OK^{-\chi'/2}
+C_O \sqrt{K} N^{\mathfrak{c}}N^{-\mathfrak{a}}
\tau^{-1},
\end{eqnarray*}
for $N$ sufficiently large on $\Omega$. Thus, choosing $\tau\ge
N^{\delta}N^{-\mathfrak{a}}$, there is a constant $\mathfrak{f}>0$
such that~(\ref{equationinmaintheoremofsectoin6}) holds. This
completes the proof of Theorem~\ref{maintheoremofsection6}.
\end{pf*}

Next, we sketch the proof of Theorem~\ref{corollarygaussian}.

\begin{pf*}{Proof of Theorem~\ref{corollarygaussian}}
The proof of Theorem~\ref{corollarygaussian} is almost identical to
the proof of Theorem~\ref{maintheoremofsection6}. In fact, it
suffices to establish Lemma~\ref{lemmafrozengapuniversality}
with~$\mu_G$ replacing ${\widehat\psi_t}\mu_G$ on the left-hand
side\vspace*{1pt} of~(\ref{equationfrozengapuniversality}). This can be
accomplished by applying Theorem~\ref{theorem31} with $\mu_G$
instead of ${\widehat\psi_t} \mu_G$: let $\tilde\bfy\in\caR
^{*}_{\mu_G}(\chi^2\sigma/2,\alpha/2)$, and let $\bfy\in\caG
(\chi^2\sigma/2,\alpha/2)$. Using the arguments of Proposition~5.2
in~\cite{singlegap}, we can rescale $\mu_G$ such that~(\ref
{assumption1thm31}) and~(\ref{assumption2thm31}) are satisfied
for $\bfy$ and $\tilde\bfy$. It is also straightforward to check
that the external potentials leading to $\mu_G^{\tilde\bfy}$,
$\tilde\bfy\in\caR^{*}_{\mu_G}(\chi^2\sigma/2,\alpha/2)$, are
$K^{\chi}$-regular. By Lemma~5.1 of~\cite{singlegap} we obtain
\begin{eqnarray*}
\bigl\llvert\bbE^{\mu_G^{\tilde\bfy}}x_k- \alpha_k \bigr
\rrvert\leq CK^{\chi} N^{-1}.
\end{eqnarray*}
Hence, using estimate~(\ref{superacc}), we conclude that
assumption~(\ref{assumption3thm31}) is also satisfied. Thus
Theorem~\ref{theorem31} yields
%
\begin{eqnarray}
\label{whatwegettothesecond} \biggl\llvert\int G_{L+p,n}(\bfx)
{\widehat
\psi_t^{\bfy}} \,\dd\mu_G^{\bfy} -
\int G_{{L}+p,n,\mathrm{sc}}(\bfx)\,\dd\mu_G^{\tilde\bfy}\biggr\rrvert
\leq C_O K^{-\mathfrak{b}},
\end{eqnarray}
for $N$ sufficiently large on $\Omega$. We refer to the proof of
Proposition~5.2 in~\cite{singlegap} for more details.

Since $\caR^{*}_{\mu_G}(\chi^2\sigma/2,\alpha/2)$ has
exponentially high probability under $\mu_G$, we can integrate over
$\tilde\bfy$ to find
\begin{eqnarray*}
\biggl\llvert\int G_{L+p,n}(\bfx) {\widehat\psi_t^{\bfy}}
\,\dd\mu_G^{\bfy} -\int G_{{L}+p,n,\mathrm{sc}}(\bfx)\,\dd
\mu_G\biggr\rrvert\leq C_O K^{-\mathfrak{b}},
\end{eqnarray*}
for $N$ sufficiently large on $\Omega$.

The proof of Theorem~\ref{corollarygaussian} is now completed in the
same way as the proof of Theorem~\ref{maintheoremofsection6}.
\end{pf*}

\section{From gap statistics to correlation functions}\label{Fromgapstatisticstocorrelationfunctions}

In this section, we translate our results on the averaged local gap
statistics into results on averaged correlation functions. Since this
procedure is fairly standard (see, e.g.,~\cite{ESYY}), we refrain from
stating all proofs in detail. We first need to slightly generalize the
setup of Section~\ref{localequilibriummeasures}.

Fix $n\in\N$, let $O$ be an $n$-particle observable and consider an
array of increasing positive integers,
%
\begin{eqnarray}
\label{arrayofinterges} \bsm=(m_1,m_2,\ldots, m_n)\in
\N^{n}.
\end{eqnarray}
Let $\alpha>0$. We define for $j\in\llbracket\alpha N,(1-\alpha
)N\rrbracket$ and $t\ge0$ an observable $G_{j,\bsm,t}\equiv
G_{j,\bsm}$ by
%
\begin{eqnarray}
\label{generalizedobservable}
G_{j,\bsm}(\bfx) &:=& O\bigl(N
\rho_j (x_{j+m_1}-x_j), N\rho
_j(x_{j+m_2}-x_j),\ldots,
\nonumber\\[-8pt]\\[-8pt]\nonumber
&&\hspace*{119pt} N\rho_j(x_{j+m_n}-x_j)\bigr),
\end{eqnarray}
where $\rho_j\equiv\widehat\rho_{\mathrm{fc}}(t,\widehat\gamma_j(t))$
denotes the density of the measure $\widehat\rho_{\mathrm{fc}}(t)$ at the
classical location of the $j$th~particle,~$\widehat\gamma_j(t)$, with
respect to the measure~$\widehat\rho_{\mathrm{fc}}(t)$. We set $G_{j,\bsm}=0$
if $j+m_n\ge(1-\alpha)N$. Similarly, we define $G_{j,\bsm,\mathrm{sc}}$ by
replacing $\rho_j$ by the density of the standard semicircle law at
the classical locations of the $j$th particle with respect to the
semicircle law; cf.~(\ref{theaveragedobservableatinfty}). The
following theorem generalizes Theorem~\ref{corollarygaussian}.

\begin{theorem}\label{generalizedgapdistribution}
Let $n\in\N$ be fixed, and let $O$ be an $n$-particle observable. Fix
small constants $\alpha,\delta>0$, and consider an interval of
consecutive integers $J\subset\llbracket\alpha N,(1-\alpha
)N\rrbracket$ in the bulk. Then there are constants $\mathfrak
{f},\delta'>0$ such that the following holds. Let $\bsm\in\N^{n}$
be an array of increasing integers [see~(\ref{arrayofinterges})]
such that $m_n\le N^{\delta'}$, and consider the observable $G_{j,\bsm
}$, respectively, $G_{j,\bsm,\mathrm{sc}}$; see~(\ref{generalizedobservable}).
Assume that $t\ge N^{-1/2+\delta}$, then
\begin{eqnarray*}
\biggl\llvert\int\frac{1}{\llvert J\rrvert }\sum_{j\in J}G_{j,\bsm}(
\bfx) f_t(\bfx)\,\dd\mu_G(\bfx)-\int\frac{1}{\llvert J\rrvert }
\sum_{j\in J}G_{j,\bsm,\mathrm{sc}}(\bfx)\,\dd
\mu_G(\bfx)\biggr\rrvert\le C_O N^{-\mathfrak{f}},
\end{eqnarray*}
for $N$ sufficiently large on $\Omega$. The constant $C_O$ depends on
$\alpha$ and $O$, and the constants $\mathfrak{f}$ and $\delta'$
depend on $\alpha$ and $\delta$.
\end{theorem}

Theorem~\ref{generalizedgapdistribution} is proven in the same way
as Theorem~\ref{corollarygaussian}. We remark that~$\delta'$ is
chosen such that~$N^{\delta'}\ll K $; that is,~$m_n$ is much smaller
than the size of the interval~$I_{L,K}$.

For $n\ge1$, define the $n$-point correlation function, $\varrho
_{f_t,n}^N$, by
\begin{eqnarray*}
\varrho^N_{f_t,n}(x_1,\ldots,x_n)
: =\int_{\R
^{N-n}} (f_t\mu
_G)^{\#} \,\dd x_{n+1}\cdots\,\dd x_{N},
\end{eqnarray*}
where $(f_t \mu_G)^{\#}$ denote the symmetrized versions of $f_t \mu
_G$. Similarly, we denote by
\begin{eqnarray*}
\varrho^N_{G,n}(x_1,\ldots,x_n)
: =\int_{\R^{N-n}} \mu_G^{\#}
\,\dd x_{n+1}\cdots\,\dd x_{N},
\end{eqnarray*}
the $n$-point correlation functions of the Gaussian ensembles;
see~(\ref{symmetrization}) with \mbox{$U\equiv0$}.

Recall that we denote by $\widehat L_{\pm}(t)$, respectively, $L_\pm
(t)$, the endpoints of the support of the measure~$\widehat\rho
_{\mathrm{fc}}(t)$, respectively, the measure~$\rho_{\mathrm{fc}}(t)$.
Recall that the two densities~$f_t$ and ${\widehat\psi_t}$ are both
conditioned on~$V$; that is, the entries~$(v_i)$ of~$V$ are considered
fixed. We have the following result on the averaged correlation
functions of $f_t \mu_G$ and ${\widehat\psi_t} \mu_G$.

\begin{theorem}\label{theoremtranslation}
Fix $n\in\N$, and choose an $n$-particle observable $O$. Fix a small
$\delta>0$, and let $t \ge N^{-1/2+\delta}$. Let $\tilde\alpha
>0$ be a small constant, and consider two energies $E\in
[L_-(t)+\tilde\alpha,L_+(t)-\tilde\alpha]$ and $E'\in
[-2+\tilde\alpha,2-\tilde\alpha]$. Then we have, for any
$\varepsilon>0$ and for $b\equiv b_N$ satisfying $\tilde\alpha/2\ge b_N>0$,
%
\begin{eqnarray}
\label{equationtheoremtranslation}
&& \biggl\llvert\int_{\R^n}\,\dd
\alpha_1\cdots\,\dd\alpha_n O(\alpha_1,\ldots,
\alpha_n)\nonumber
\\
&&\hspace*{18pt}{}\times\biggl[\int_{E-b}^{E+b}
\frac
{\dd x}{2b}\frac{1}{{[\rho_{\mathrm{fc}}(t,E)]}^n} \varrho^N_{f_t,n}
\biggl(x+\frac{\alpha_1}{N\rho_{\mathrm{fc}}(t,E)},\ldots,x+\frac{\alpha
_n}{N\rho_{\mathrm{fc}}(t,E)} \biggr)
\nonumber\\[-8pt]\\[-8pt]\nonumber
&&\hspace*{39pt}{} -\int
_{E'-b}^{E'+b}\frac{\dd x}{2b}\frac{1}{{[\rho
_{\mathrm{sc}}(E')]}^n}
\varrho^N_{G,n} \biggl(x+\frac{\alpha_1}{N\rho
_{\mathrm{sc}}(E')},\ldots,x+
\frac{\alpha_n}{N\rho_{\mathrm{sc}}(E')} \biggr) \biggr]\biggr\rrvert\hspace*{-10pt}
\nonumber
\\
&&\qquad\le C_ON^{2\varepsilon} \bigl(b^{-1}
N^{-1+\varepsilon} +N^{-\mathfrak{f}}+N^{-c{\alpha_0}} \bigr),\nonumber
\end{eqnarray}
for $N$ sufficiently large on $\Omega$. Here $\mathfrak{a}$ is the
constant in the rigidity estimate~(\ref{masterrigiditybound}), and
$\mathfrak{f}$ is the constant in Theorem~\ref
{generalizedgapdistribution}. Moreover, $\rho_{\mathrm{fc}}(t,E)$ stands for the
density of
the ($N$-independent) measure $\rho_{\mathrm{fc}}(t)$ at the energy $E$. The
constant $C_O$ depends on $O$ and $\tilde\alpha$. Further,
$\alpha_0$ is the constant appearing in Assumption~\ref
{assumptionmuVconvergence}. The constant $c$ depends on the measure~$\nu$.
\end{theorem}

Theorem~\ref{theoremtranslation} follows from Theorem~\ref
{corollarygaussian}. This is an application of Section~7 in~\cite
{ESYY}. The
validity of Assumption~IV in~\cite{ESYY} is a direct consequence of
the local law in Theorem~\ref{thmstrong}. Further, we remark that the
parameter $b_N$ in Theorem~\ref{theoremtranslation} and the interval
of consecutive integers $J$ in Theorem~\ref{generalizedgapdistribution}
are related by $J=\{i\dvtx \widehat\gamma_i(t)\in[
E-b_N,E+b_N]\}$, where $\widehat\gamma_i(t)$ are the classical
locations with respect to the measure $\widehat\rho_{\mathrm{fc}}(t)$. This
explains, up to minor technicalities, $b_N\gg N^{-1}$. Then Section~7
of~\cite{ESYY} yields~(\ref{equationtheoremtranslation}) formulated
in terms of $\widehat\rho_{\mathrm{fc}}(t)$ instead of $\rho_{\mathrm{fc}}(t)$.
Using~(\ref{minilocallaw}) and the smoothness of~$O$, we can replace
$\widehat\rho_{\mathrm{fc}}$ by $\rho_{\mathrm{fc}}$ at the expense of an error of
size $C_ON^{-c\alpha_0}$. This
eventually
gives~(\ref{equationtheoremtranslation}) with $\rho_{\mathrm{fc}}(t)$.

\section{Proofs of main results}\label{Proofsofmainresults}

Theorem~\ref{theoremtranslation} shows that the averaged local
correlation functions of ensembles of the form
\begin{eqnarray*}
H_t=\mathrm{e}^{-(t-t_0)/2}V+\mathrm{e}^{-t/2}W+\bigl(1-
\mathrm{e}^{-t}\bigr)^{1/2}{W'},
\end{eqnarray*}
with some small $t_0\ge0$, and with $ W'$ a GUE/GOE matrix independent
of~$W$ and~$V$, can be compared with the averaged local correlation
functions of the GUE, respectively, GOE, for times satisfying $t\gg
N^{-1/2}$. In this section, we explain how this can be used to prove
the universality at time $t=0$.

\subsection{Green function comparison theorem}

We start with a Green function comparison theorem. Assume that we are
given two complex Hermitian or real symmetric Wigner matrices, $X$ and
$Y$, both satisfying the assumptions in Definition~\ref
{assumptionwigner}. Let~$V$ be a real random or deterministic diagonal matrix
satisfying Assumptions~\ref{assumptionmuV} and~\ref
{assumptionmuVconvergence}. Consider the deformed Wigner matrices
%
\begin{eqnarray}
\label{HXHY} H^X: =V+X,\qquad H^Y
: =V+Y,
\end{eqnarray}
of size $N$. The main theorem of this subsection, Theorem~\ref
{corollaryofgreenfunctioncomparison}, states that the correlation
functions of the two matrices $H^X$ and $H^Y$, when conditioned on~$V$,
are identical on scale~$1/N$ provided that the first four moments
of~$X$ and~$Y$ almost match. Theorem~\ref
{corollaryofgreenfunctioncomparison} is a direct consequence of the
Green function comparison
Theorem~\ref{thegreenfunctioncomparisontheorem}.

Denote the Green functions of $H^X$, $H^Y$, respectively, by
\begin{eqnarray*}
G^X(z): =\frac{1}{H^X-z}, \qquad
G^Y(z): =\frac{1}{H^Y-z}\qquad(z\in\C
\setminus\R),
\end{eqnarray*}
and set $m_N^X(z): =N^{-1}\Tr G^X(z)$,
$m_N^Y(z): =N^{-1}\Tr
G^Y(z)$. From Theorem~\ref{thmstrong}, we know that, for all $z\in
\caD_L$ [see~(\ref{eqDL})], with $L\ge40\xi$,
%
\begin{eqnarray}
\bigl\llvert m_N^X(z)-\widehat m_{\mathrm{fc}}(z)
\bigr\rrvert\le(\varphi_N)^{c\xi}\frac
{1}{N\eta}
\end{eqnarray}
and
%
\begin{eqnarray}
\label{eqofthmstrong} \bigl\llvert G_{ij}^X(z)-
\delta_{ij} \widehat g_{i}(z)\bigr\rrvert\le(
\varphi_N)^{c\xi
} \biggl(\sqrt{\frac{\im\widehat m_{\mathrm{fc}}(z)}{N\eta}}+
\frac{1}{N\eta
} \biggr),
\end{eqnarray}
with $(\xi,\upsilon)$-high probability on $\Omega$ for some
$\upsilon>0$ and $c>0$, where
\begin{eqnarray*}
\widehat g_i(z): =\frac{1}{v_i-z-\widehat
m_{\mathrm{fc}}(z)}\qquad(z\in
\C\setminus\R).
\end{eqnarray*}
Here, $\widehat m_{\mathrm{fc}}$, is the Stieltjes transform of the
measure~$\widehat\rho_{\mathrm{fc}}$, which agrees with~$\widehat\rho
_{\mathrm{fc}}^{\vartheta}$ for the choice~$\vartheta=1$ and with~$\widehat
\rho_{\mathrm{fc}}(t)$ for the choice $t=t_0$. The identical estimates hold
true when $X$ is replaced by $Y$.

Recall\vspace*{1pt} that we denote by $\widehat L_{\pm}$ the endpoints of the
support of~$\widehat\rho_{\mathrm{fc}}$, and that we denote by $\widehat
\kappa_E\equiv\widehat\kappa$ the distance of~$E\in[\widehat
L_-,\widehat L_+]$ to the endpoints~$\widehat L_{\pm}$. Adapting the
Green function theorem of~\cite{EYY1} we obtain the following theorem.

\begin{theorem}[(Green function comparison theorem)]\label
{thegreenfunctioncomparisontheorem}
Assume that $X$ and $Y$ satisfy the assumptions in Definition~\ref
{assumptionwigner}, and let~$V$ satisfy Assumptions~\ref
{assumptionmuVconvergence} and~\ref{assumptionmuV}. Assume further that the
first two moments of $X=(x_{ij})$ and $Y=(y_{ij})$ agree and that the
third and forth moments satisfy
%
%
\begin{eqnarray}
\label{matching1} \bigl\llvert\E\bar{x}_{ij}^{p}
{x}_{ij}^{3-p}-\E\bar{y}_{ij}^{p}
{y}_{ij}^{3-p}\bigr\rrvert\le N^{-\delta-2}\qquad\bigl(p
\in\llbracket0,3\rrbracket\bigr),
\end{eqnarray}
respectively,
%
\begin{eqnarray}
\label{matching2} \bigl\llvert\E\bar{x}_{ij}^{q}
{x}_{ij}^{4-q}-\E\bar{y}_{ij}^{q}
{y}_{ij}^{4-q}\bigr\rrvert\le N^{-\delta}\qquad\bigl(q\in
\llbracket0,4\rrbracket\bigr),
\end{eqnarray}
for some given $\delta>0$.

Let $\varepsilon> 0$ be arbitrary, and let $ N^{-1-\varepsilon} \leq\eta
\leq N^{-1}$. Fix $N$-independent integers $k_1,\ldots,k_n$ and
energies\vspace*{1pt} $E^1_j,\ldots,E^{k_j}_j$, $j=1,\ldots,n$, with $\widehat
\kappa> \tilde\alpha$ for all $E^k_j$ with some fixed
$\tilde{\alpha}>0$. Define $z^k_j : =E^k_j \pm
\mathrm{i} \eta$, with
the sign arbitrarily chosen. Suppose that $F$ is a smooth function such
that for any multi-index ${\sigma} = (\sigma_1,\ldots,\sigma_n)$,
with $1\le\llvert {\sigma}\rrvert \le5$, and any $\varepsilon'>0$ sufficiently
small, there is a $C_0 > 0$ such that
\begin{eqnarray*}
\max\Bigl\{ \bigl\llvert\partial^{\sigma} F(x_1,
\ldots,x_n)\bigr\rrvert\dvtx \max_j \llvert
x_j\rrvert\leq N^{\varepsilon'} \Bigr\} & \leq& N^{C_0\varepsilon'},
\\
\max\Bigl\{ \bigl\llvert\partial^{\sigma} F(x_1,
\ldots,x_n)\bigr\rrvert\dvtx \max_j \llvert
x_j\rrvert\leq N^{2} \Bigr\} & \leq& N^{C_0},
\end{eqnarray*}
for some $C_0$.

Then there exists a constant $C_1$, depending on $\sum_m k_m$, $C_0$
and the constants in~(\ref{eqC0}), such that for any~$\eta$ with
$N^{-1-\varepsilon} \le\eta\le N^{-1}$,
%
\begin{eqnarray}
&&\Biggl\llvert\bbE F \Biggl(\frac{1}{N^{k_1}} \Tr\prod
_{j=1}^{k_1} G^X\bigl(z^1_j
\bigr),\ldots,\frac{1}{N^{k_n}} \Tr\prod_{j=1}^{k_n}
G^X\bigl(z^n_j\bigr) \Biggr) \nonumber
\\
&&\quad{} - \bbE F
\Biggl(\frac{1}{N^{k_1}} \Tr\prod_{j=1}^{k_1}
G^Y\bigl(z^1_j\bigr),\ldots,
\frac
{1}{N^{k_n}} \Tr\prod_{j=1}^{k_n}
G^Y\bigl(z^n_j\bigr) \Biggr) \Biggr\rrvert
\\
&&\qquad\leq C_1 N^{-1/2 + C_1 \varepsilon} + C_1
N^{-1/2 + \delta+ C_1 \varepsilon},\nonumber
\end{eqnarray}
for $N$ sufficiently large on $\Omega$.
\end{theorem}

Theorem~\ref{thegreenfunctioncomparisontheorem} is proven in the
same way as Theorem~2.3 in~\cite{EYY} with the following
modifications. Fix some labeling of $\{ (i,j)\dvtx 1\leq i \leq j \leq
N\}$
by $\llbracket1,\gamma(N)\rrbracket$, with $\gamma(N): =
N(N+1)/2$, and write the $\gamma$th element of this labeling as
$(i_\gamma,j_\gamma)$. Starting with $W^{(0)}\equiv X$, inductively
define $W^{(\gamma)}$ by replacing\vspace*{1pt} the $(i_\gamma,j_\gamma
),(j_\gamma,i_\gamma)$ entries of $W^{(\gamma-1)}$ by the
corresponding entries of $Y$. Moreover set $H^{(\gamma)}: =
V+W^{(\gamma)}$. Thus we have $H^{(0)}=H^X$, $H^{(\gamma(N))}=H^Y$,
and $H^{(\gamma)} - H^{(\gamma-1)}$ is zero\vadjust{\goodbreak} in all but two entries
for every $\gamma$. In short, we use a Lindeberg-type replacement
strategy: we successively replace the entries of the matrix $X$ by
entries of the matrix $Y$. Note, however, that the entries of the
matrix~$V$ are not changed.

The main technical input in the proof of Theorem~2.3 in~\cite{EYY1} is
estimate~(2.21) in that publication. For the case at hand the
corresponding estimate reads as follows: let $\xi$ satisfy~(\ref
{eqxi}). Then, for all $\delta>0$, and any $N^{-1/2} \gg y\ge
N^{-1+\delta}$, we have
%
\begin{eqnarray}
\label{theestimateofthegreenfunctioncomparison}
&&\bbP\biggl( \max
_{\gamma\leq\gamma(N)} \max
_k \sup_{E\dvtx \widehat\kappa_E \geq
\tilde{\alpha} }\biggl\llvert\biggl(
\frac{1}{H^{(\gamma)} - E -
\ii y} \biggr)_{kk} \biggr\rrvert\geq N^{2\delta}
\biggr)
\nonumber\\[-8pt]\\[-8pt]\nonumber
&&\qquad \le\mathrm{e}^{{-\upsilon}(\varphi_N)^{\xi}},
\end{eqnarray}
on $\Omega$ for $N$ sufficiently large, where $\upsilon>0$ depends
only on $\tilde\alpha$, $\delta$ and the constants in~(\ref
{eqC0}). Estimate~(\ref{theestimateofthegreenfunctioncomparison})
follows easily from the local law in~(\ref{eqofthmstrong}), the
stability bound~(\ref{stabilitybound}) and
Lemma~\ref{lemmahatmfc}. The rest of the proof of Theorem~\ref
{thegreenfunctioncomparisontheorem} is identical to the proof in~\cite
{EYY1}. (The matching conditions in~(\ref{matching1}) are weaker than
in~\cite{EYY1}, but the proof carries over without any changes.)

Lindeberg's replacement method was applied in random matrix theory
in~\cite{C} to compare traces of Green functions. This idea was also
used in~\cite{TV1} in the proof of the ``four moment theorem''
that compares individual eigenvalue distributions. The four-moment
matching conditions~(\ref{matching1}) and~(\ref{matching2})
appeared first in~\cite{TV1} with $\delta=0$. The ``Green function
comparison theorem'' of~\cite{EYY1} compares Green functions at fixed
energies. Since the approach in~\cite{TV1} requires additional
difficult estimates due to singularities from neighboring eigenvalues,
we follow the method of~\cite{EYY1}, where difficulties stemming from
such resonances are absent. For deformed Wigner matrices with
deterministic potential the approach of~\cite{TV1} was recently
followed in~\cite{OV} where a ``four moment theorem'' was established.
It allows one to compare local correlation functions of the matrices
$V+W$ and $V+{W'}$ for fixed~$V$, where~$W$ and~$W'$ are real symmetric or
complex Hermitian Wigner matrices, provided that the moments of the
off-diagonal entries of~$W$ and~${W'}$ match to fourth
order.

The Green function comparison theorem leads directly to the equivalence
of local statistics for the matrices~$H^Y$ and~$H^X$.

%
\begin{theorem}\label{corollaryofgreenfunctioncomparison}
Assume that $X$, $Y$ are two complex Hermitian or two real symmetric
Wigner matrices satisfying assumptions in Definition~\ref
{assumptionwigner}. Assume further that $X$ and $Y$ satisfy the matching
conditions~(\ref{matching1}) and~(\ref{matching2}), for some
$\delta>0$. Let~$V$ be a deterministic real diagonal matrix satisfying
the Assumptions~\ref{assumptionmuVconvergence} and~\ref{assumptionmuV}.
Denote by $\varrho^N_{H^X,n}, \varrho^N_{H^Y,n}$ the $n$-point
correlation functions of the eigenvalues with respect to the
probability laws of the matrices $H^X$, $H^Y$, respectively. Then, for
any energy $E$ in the interior of the support of $\rho_{\mathrm{fc}}$ and any
$n$-particle\vadjust{\goodbreak} observable~$O$, we have
\begin{eqnarray*}
&& \lim_{N\to\infty} \int_{\bbR^k} \,\dd
\alpha_1\cdots\,\dd\alpha_n O(\alpha_1,\ldots,
\alpha_n)
\\
&&\hspace*{41pt}{}\times \biggl[\varrho^N_{H^X,n} \biggl(E+
\frac{\alpha_1}{N},\ldots,E+\frac{\alpha_n}{N} \biggr)
-\varrho^N_{H^Y,n} \biggl(E+\frac{\alpha_1}{N},\ldots,E+
\frac{\alpha_n}{N} \biggr) \biggr]
\\
&&\qquad = 0,
\end{eqnarray*}
for any fixed $n\in\N$.
\end{theorem}

Notice that this comparison theorem holds for any fixed energy $E$ in
the bulk. The proof of~\cite{EYY1} applies almost verbatim. The only
technical input in the proof is the local law for $m_N^X$, respectively,
$m_N^Y$, on scales $\eta\sim N^{-1+\varepsilon}$, which we have
established in Theorem~\ref{thmstrong}; see also~(\ref{eqofthmstrong}).

\subsection{Proof of Theorem~\texorpdfstring{\protect\ref{theorem1}}{2.5}}\label{leproofofletheoremun}
In the remaining subsections,~\ref{leproofofletheoremun} and~\ref
{subsectionproofrandomV}, we complete the proofs of our main results
in Theorems~\ref{theorem1}~and~\ref{theorem2}. The proofs
for deterministic and random~$V$ differ slightly. We start with the
case of deterministic~$V$ in this subsection; the random case is
treated in Section~\ref{subsectionproofrandomV}.

\begin{pf*}{Proof of Theorem~\ref{theorem1}}
\label{proofoftheorem1}
Assume that $W=(w_{ij})$ is a complex Hermitian or a real symmetric
Wigner matrix satisfying the assumptions in Definition~\ref
{assumptionwigner}. Let $V=\diag(v_i)$ be a deterministic real diagonal matrix
satisfying Assumptions~\ref{assumptionmuVconvergence} and~\ref
{assumptionmuV}. (Note that the event $\Omega$ then has full
probability.) Set $H=(h_{ij})=V+W$. Let~$E\in\R$ be inside the
support of~$\rho_{\mathrm{fc}}$. Note that by Lemma~\ref{hatmfc},~$E$ is also
contained in the support of~$\widehat\rho_{\mathrm{fc}}$, for~$N$ sufficiently
large. (Here we have $\widehat\rho_{\mathrm{fc}}=\widehat\rho
_{\mathrm{fc}}^{\vartheta=1}$ and similarly for $\rho_{\mathrm{fc}}$.) Fix $\delta
'>0$, and set $t\equiv N^{-1/2+\delta'}$. We first claim that there
exists an auxiliary complex Hermitian or real symmetric Wigner matrix,
$U=(u_{ij})$, satisfying the assumptions in Definition~\ref
{assumptionwigner} such that the following holds: set
%
\begin{eqnarray}
\label{thematrixY} Y: =\mathrm{e}^{-t/2} U+\bigl(1-
\mathrm{e}^{t/2}\bigr){W'},
\end{eqnarray}
where $ W'$ is a GUE/GOE matrix independent of~$W$. Then the moments of
the entries of~$Y$ satisfy
%
\begin{eqnarray}
\label{matching3} \E\bar{y}_{ij}^{p}
{y}_{ij}^{3-p}=\E\bar{w}_{ij}^{p}
{w}_{ij}^{3-p},\qquad\bigl\llvert\E\bar{y}_{ij}^{q}
{y}_{ij}^{4-q}-\E\bar{w}_{ij}^{q}
{w}_{ij}^{4-q}\bigr\rrvert\le Ct,
\end{eqnarray}
for $p\in\llbracket0,3\rrbracket$, $q\in\llbracket0,4\rrbracket$,
where $(w_{ij})$ are the entries of the Wigner matrix~$W$.

Assuming the existence of such a Wigner matrix $U$, we choose
$t_0\equiv t$ and set
\begin{eqnarray*}
H_t &:=&\mathrm{e}^{-(t-t_0)/2}V+\mathrm
{e}^{-t/2}U+\bigl(1-\mathrm{e}^{-t}\bigr)^{1/2}{W'}
\\
&=&V+\mathrm{e}^{-t/2}U+\bigl(1-\mathrm{e}^{-t}
\bigr)^{1/2}{W'}.
\end{eqnarray*}
Then the matrices $H_t$ and $H=V+W$ satisfy the matching
conditions~(\ref{matching1}) and~(\ref{matching2}) of Theorem~\ref
{thegreenfunctioncomparisontheorem} (with, say, $\delta
=1/4-2\delta'$). This follows from~(\ref{matching3}). Thus
Theorem~\ref{corollaryofgreenfunctioncomparison} implies that the
correlation functions of $H_t$ and $H$ agree in the limit of large $N$,
that is,
%
\begin{eqnarray}
\label{almostattheuniversality}
&& \lim_{N\to\infty}\int_{\bbR^n} \,\dd
\alpha_1\cdots\,\dd\alpha_n O(\alpha_1,\ldots,
\alpha_n)\nonumber
\\
&&\hspace*{42pt}{}\times\biggl[\frac
{1}{2b}\int_{E-b}^{E+b}
\frac{\dd x}{[\rho_{\mathrm{fc}}(E)]^n}\varrho^N_{H_t,n} \biggl(x+
\frac{\alpha_1}{\rho_{\mathrm{fc}}(E)N},\ldots,x+\frac
{\alpha_n}{\rho_{\mathrm{fc}}(E)N} \biggr)
\nonumber\\[-8pt]\\[-8pt]\nonumber
&&\hspace*{58pt}{}-\frac{1}{2b}\int_{E-b}^{E+b}
\frac{\dd x}{[\rho_{\mathrm{fc}}(E)]^n}\varrho^N_{H,n} \biggl(x+
\frac{\alpha_1}{\rho_{\mathrm{fc}}(E)N},\ldots,x+\frac{\alpha_n}{\rho
_{\mathrm{fc}}(E)N} \biggr) \biggr]\hspace*{-10pt}
\\
&&\qquad  = 0,\nonumber
\end{eqnarray}
where $(\varrho_{H,n}^{N})$ denote the correlation functions of
$H=V+W$ and where $(\varrho_{H_t,n}^N)$ denote the correlation
functions of $H_t$. [In fact,~(\ref{almostattheuniversality}) holds
even without the averages in the energy around $E$.]

On the other hand, for small $\delta>0$, Theorem~\ref
{theoremtranslation} assures that the local correlation functions of
the matrix
$H_t$ agree with the correlation functions of the GUE (resp., GOE), when
averaged over an interval of size $b$, with $1\gg b\ge N^{-\delta}$;
that is, for any $E'$ with $\llvert E'\rrvert <2$,
%
%
\begin{eqnarray}
\label{almostattheuniversality2}
&& \lim_{N\to\infty}\int_{\bbR^n} \,\dd
\alpha_1\cdots\,\dd\alpha_n O(\alpha_1,\ldots,
\alpha_n)\nonumber
\\
&&\hspace*{42pt}{} \times\biggl[\frac
{1}{2b}\int_{E-b}^{E+b}
\frac{\dd x}{[\rho_{\mathrm{fc}}(E)]^n}\varrho^N_{H_t,n} \biggl(x+
\frac{\alpha_1}{\rho_{\mathrm{fc}}(E)N},\ldots,x+\frac
{\alpha_n}{\rho_{\mathrm{fc}}(E)N} \biggr)
\nonumber\\[-8pt]\\[-8pt]\nonumber
&&\hspace*{58pt}{}-\frac{1}{[\rho_{\mathrm{sc}}(E')]^n}\varrho^N_{G,n}
\biggl(E'+\frac
{\alpha_1}{\rho_{\mathrm{sc}}(E')N},\ldots,E'+
\frac{\alpha_n}{\rho_{\mathrm{sc}}(E')N} \biggr) \biggr]
\\
&&\qquad = 0,\nonumber
\end{eqnarray}
where $(\varrho_{G,n}^N)$ denote the correlations functions of the
GUE, respectively, GOE. Combining~(\ref{almostattheuniversality})
and~(\ref{almostattheuniversality2}), we get~(\ref{theorem1equation1}).

Thus to complete the proof we need to show the existence of a Wigner
matrix~$U$ with the properties described above.
For a real random variables~$\zeta$, denote by~$m_k(\zeta)=\E
\zeta^k$, $k\in\N$, its moments.

\begin{lemma}[(Lemma 6.5 in~\cite{EYY1})]\label{existenceofmatching}
Let $m_3$ and $m_4$ be two real numbers such that
\begin{eqnarray*}
m_4 - m_2^2-1\ge0,\qquad m_4\le
C_1,
\end{eqnarray*}
for some constant $C_1$. Let $\zeta_G$ be a Gaussian random variable
with mean $0$ and variance $1$. Then
for any sufficient small $\gamma>0$, depending on $C_1$, there exists
a real random variable $\zeta_\gamma$ with subexponential
decay and independent of $\zeta_G$, such that the first three moments of
\begin{eqnarray*}
\zeta': =(1 - \gamma)^{1/2}
\zeta_\gamma+ \gamma^{1/2} \zeta_G
\end{eqnarray*}
are $m_1 (\zeta' ) = 0$, $m_2 (\zeta') = 1$, $m_3 (\zeta') = m_3$,
and the forth moment $m_4 (\zeta')$ satisfies
\begin{eqnarray*}
\bigl\llvert m_4 \bigl(\zeta' \bigr)-m_4
\bigr\rrvert\le C\gamma,
\end{eqnarray*}
for some $C$ depending on $C_1$.
\end{lemma}

Since the real and imaginary parts of~$W$ are independent, it is sufficient
to match them individually; that is, we apply Lemma~\ref
{existenceofmatching} separately to the real and imaginary parts of $(w_{ij})$.
This completes the proof of Theorem~\ref{theorem1} for deterministic~$V$.
\end{pf*}

\subsection{Proof of Theorem~\texorpdfstring{\protect\ref{theorem2}}{2.6}}\label{subsectionproofrandomV}
Next, we prove Theorem~\ref{theorem2}. Assume that $W=(w_{ij})$ is a
complex Hermitian or a real symmetric Wigner matrix satisfying the
assumption in Definition~\ref{assumptionwigner}. Let $V=\diag(v_i)$
be a random real diagonal matrix satisfying Assumptions~\ref
{assumptionmuVconvergence} and~\ref{assumptionmuV}. Denote by
$\tilde{f}_t \mu_G$ the distribution of the eigenvalues of the matrix
\begin{eqnarray*}
H_t: =V+\mathrm{e}^{-t/2}W+\bigl(1-\mathrm
{e}^{t}\bigr)^{1/2}{W'} \qquad(t\ge0),
\end{eqnarray*}
where~$W'$ is a GUE/GOE matrix independent of~$V$ and~$W$. Let $\E^V$
stand for the expectation with respect to the law of the
entries~$(v_i)$ of~$V$. Recall the definition of the event $\Omega$ in
Definition~\ref{definitionofomegaV}. Following the notation of
Section~\ref{sectionDBM}, $f_t \mu_G\equiv f_t^V \mu_G$ denotes
the density conditioned on~$V$. For an $n$-particle observable $O$ and
for $G_{j,\bsm}$ as in~(\ref{generalizedobservable}), we may write
\begin{eqnarray*}
&& \int\frac{1}{\llvert J\rrvert }\sum_{j\in J}G_{j,\bsm}(
\bfx) \tilde{f}_t(\bfx)\,\dd\mu_G(\bfx)
\nonumber
\\
&&\qquad =\E^V \biggl[ \biggl(\int\frac{1}{\llvert J\rrvert }\sum
_{j\in J}G_{j,\bsm}(\bfx) f^V_t(
\bfx)\,\dd\mu_G(\bfx) \biggr)\lone(\Omega) \biggr]+\caO\bigl(
N^{-\mathfrak
{t}} \bigr),
\end{eqnarray*}
where $\mathfrak{t}>0$ is the constant in~(\ref{eqassumptionmuV2}) of
the Assumptions in~\ref{assumptionmuV}. Here we use the
definition of $\Omega$. Since $(v_i)$ are i.i.d.,~(\ref
{definitionofomegaV}) holds with exponentially high probability.
Estimate~(\ref
{definitionofomegaVIII}) holds with probability large than
$1-N^{-\mathfrak{t}}$ by Assumption~\ref{assumptionmuV}. Hence
$\mathbb{P}^V(\Omega^c)\le cN^{-\mathfrak{t}}$, for some $c>0$ and
$N$ sufficiently large.

Using Theorem~\ref{generalizedgapdistribution}, we find that
%
\begin{eqnarray}
\label{onthewayforrandomV1}
&& \int\frac{1}{\llvert J\rrvert }\sum_{j\in
J}G_{j,\bsm}(\bfx) \tilde{f}_t(\bfx)\,\dd\mu_G(\bfx)
\nonumber\\[-8pt]\\[-8pt]\nonumber
&&\qquad =\int
\frac{1}{\llvert J\rrvert }\sum_{j\in
J}G_{j,\bsm,\mathrm{sc}}(\bfx)\,\dd\mu_G(\bfx)
+\caO\bigl( N^{-\mathfrak{f}} \bigr)+\caO\bigl( N^{-\mathfrak
{t}} \bigr),
\end{eqnarray}
where we use once more the estimates on the event $\Omega$. Here
$\mathfrak{f}>0$ is the constant appearing in Theorem~\ref
{generalizedgapdistribution}.

To establish the equivalent result to Theorem~\ref{theoremtranslation},
we need a local deformed semicircle law for the setting
when the entries~$(v_i)$ of~$V$ are not fixed. Recall that we denote by
$m_{\mathrm{fc}}$ the Stieltjes transform of the deformed semicircle law $\rho
_{\mathrm{fc}}=\rho_{\mathrm{fc}}^{\vartheta=1}$.

\begin{lemma}[(Theorems~2.10 and 2.21 in~\cite{LS})]\label{locallawrandomV}
Let~$W$ be a complex Hermitian or a real symmetric Wigner matrix
satisfying the assumptions in Definition~\ref{assumptionwigner}.
Let~$V$ be a random real diagonal matrix satisfying Assumptions~\ref
{assumptionmuVconvergence} and~\ref{assumptionmuV}. Set
$H: =
V+W$, $G(z): =(H-z)^{-1}$ and $m_N(z): =N^{-1}\Tr G(z)$, $(z\in\C
^+)$. Let $\xi=A_0\log\log N/2$; see~(\ref{eqxi}). Then there
exists $\upsilon>0$ and $c$ [both depending on the constants in~(\ref
{eqC0}), the constants $A_0$, $E_0$ in~(\ref{eqDL}) and the
measure~$\nu$], such that for $L\ge40\xi$, we have
%
\begin{equation}
\label{locallawrandomVequation} \qquad\bigl\llvert m_N(z)-m_{\mathrm{fc}}(z)\bigr
\rrvert
\le(\varphi_N)^{c\xi} \biggl(\min\biggl\{\frac
{1}{N^{1/4}},
\frac{1}{\sqrt{\kappa_E+\eta}}\frac{1}{\sqrt
{N}} \biggr\}+\frac{1}{N\eta} \biggr),
\end{equation}
and
%
\begin{eqnarray}
\label{lcallawrandomV2} \bigl\llvert G_{ij}(z)-\delta_{ij}
g_{i}(z) \bigr\rrvert\le(\varphi_N)^{c\xi
}
\biggl(\sqrt{\frac{\im m_{\mathrm{fc}}(z)}{N\eta}}+\frac{1}{N\eta} \biggr),
\end{eqnarray}
$i,j\in\llbracket1,N\rrbracket$, with $(\xi,\upsilon)$-high
probability, for all $z=E+\ii\eta\in\caD_L$; see~(\ref{eqDL}).
Here, we have set
\begin{eqnarray*}
{g}_i(z): =\frac{1}{ v_i-z-{m}_{\mathrm{fc}}(z)} \qquad\bigl
(z\in
\C^+,i\in\llbracket1,N\rrbracket\bigr).
\end{eqnarray*}
Moreover, fixing $\alpha>0$, there is $c_1$ [depending on the
constants in~(\ref{eqC0}), the constants $A_0$, $E_0$ in~(\ref
{eqDL}), the measure $\nu$ and $\alpha$], such that
%
\begin{eqnarray}
\label{rigidityrandomV} \llvert\lambda_i-\gamma_i\rrvert
\le(\varphi_N)^{c_1\xi}\frac{1}{\sqrt{N}},
\end{eqnarray}
with $(\xi,\upsilon)$-high probability, for all $i\in\llbracket
\alpha N,(1-\alpha)N\rrbracket$. Here $(\lambda_i)$ denote the
eigenvalues of $H=V+W$, and~$(\gamma_i)$ are their classical locations
with respect the deformed semicircle law $\rho_{\mathrm{fc}}$.
\end{lemma}

Using the local law in~Lemma~\ref{locallawrandomV}, we obtain
from~(\ref{onthewayforrandomV1}) equivalent results to
Theorem~\ref{theoremtranslation}.

\begin{theorem}\label{theoremtranslationrandom}
Fix $n\in\N$, and consider an $n$-particle observable $O$. Fix
\mbox{$\delta>0$}, and let $t \ge N^{-1/4+\delta}$. Let $\tilde\alpha
>0$ be a small constant, and consider two energies $E\in
[L_-(t)+\tilde\alpha,L_+(t)-\tilde\alpha]$ and $E'\in
[-2+\tilde\alpha,2-\tilde\alpha]$. Then, for any $\varepsilon
>0$ and for $b\equiv b_N$ satisfying $\tilde\alpha/2\ge b_N>0$,
we have
\begin{eqnarray*}
&& \biggl\llvert\int_{\R^n}\,\dd\alpha_1\cdots\,\dd
\alpha_n O(\alpha_1,\ldots,\alpha_n)
\\
&&\hspace*{19pt}{} \times
\biggl[\int_{E-b}^{E+b}\frac
{\dd x}{2b}
\frac{1}{[\rho_{\mathrm{fc}}(t,E)]^n} \varrho^N_{\tilde
{f}_t,n} \biggl(x+
\frac{\alpha_1}{N\rho_{\mathrm{fc}}(t,E)},\ldots,x+\frac
{\alpha_n}{N\rho_{\mathrm{fc}}(t,E)} \biggr)
\\
&&\hspace*{41pt}{}-\int
_{E'-b}^{E'+b}\frac{\dd x}{2b}\frac{1}{[\rho
_{\mathrm{sc}}(E')]^n}
\varrho^N_{G,n} \biggl(x+\frac{\alpha_1}{N\rho
_{\mathrm{sc}}(E')},\ldots,x+
\frac{\alpha_n}{N\rho_{\mathrm{sc}}(E')} \biggr) \biggr]\biggr\rrvert
\nonumber
\\
&&\qquad \le C_ON^{2\varepsilon} \bigl(b^{-1}
N^{-1/2+\varepsilon} +N^{-\mathfrak{f}}+N^{-\mathfrak
{t}}+N^{-1/4} \bigr),
\end{eqnarray*}
for $N$ sufficiently large. Here, $\mathfrak{f}>0$ is the constant in
Theorem~\ref{generalizedgapdistribution}. Moreover, $\rho_{\mathrm{fc}}(E)$
stands for the density of the ($N$-independent) measure $\rho_{\mathrm{fc}}$ at
the energy $E$. The constant $C_O$ depends on $O$, $\tilde\alpha$
and the measure $\nu$. The constant~$\mathfrak{f}$ depends on $\delta
$ and $\tilde\alpha$.
\end{theorem}

The proof of Theorem~\ref{theoremtranslationrandom} is an
application of Section~7 in~\cite{ESYY}. The validity of Assumption~IV
in~\cite{ESYY} is a direct consequence of the local law in Lemma~\ref
{locallawrandomV}. Here and also below, we use that the local laws
of~Lemma~\ref{locallawrandomV} are only used on very small
scales~$\eta\sim N^{-1+\varepsilon}$ in the bulk. For such small~$\eta$
the first error term in~(\ref{locallawrandomVequation}) is
negligible compared to the second error term. Also note that the first
term on the right-hand side of the estimate in Theorem~\ref
{theoremtranslationrandom} is bigger than the corresponding term
in~(\ref
{equationtheoremtranslation}). This is due to the weaker rigidity
bounds in case~$V$ is random; see~(\ref{rigidityrandomV}). We
therefore have to impose that $b\gg N^{-1/2}$ in order to have a
vanishing error term in the limit of large $N$. Finally, we mention
that the error term $C_ON^{2\varepsilon}N^{-1/4}$ stems
from replacing $\widehat\rho_{\mathrm{fc}}(t,E)$ by $\rho_{\mathrm{fc}}(t,E)$; see the
comment below Theorem~\ref{theoremtranslation}.

\begin{pf*}{Proof of Theorem~\ref{theorem2}}
The proof Theorem~\ref{theorem2} follows now along the lines of the
proof of Theorem~\ref{theorem1}. First, we check that the Green
function comparison Theorem~\ref{thegreenfunctioncomparisontheorem}
holds true for $H^X=V+X$, respectively, $H^Y=V+Y$ with
random~$V$. This is indeed the case, since the only input we used is
estimate~(\ref{theestimateofthegreenfunctioncomparison}), which
also holds for random~$V$ by the local laws in Lemma~\ref
{locallawrandomV} and the stability estimate~(\ref{stabilitybound}). Note
that we are using that bound~(\ref
{theestimateofthegreenfunctioncomparison}) is only required on
scales~$\eta\ll N^{-1/2}$. Similarly,
we can establish Theorem~\ref{corollaryofgreenfunctioncomparison}
for random~$V$ using the Green function comparison theorem for
random~$V$, the local laws in Lemma~\ref{locallawrandomV} and the
stability estimate~(\ref{stabilitybound}). Finally, we note that the
construction of the
matrix $U$ and $Y$ [see~(\ref{thematrixY})] and the moment matching
in~(\ref{matching3}) do not involve~$V$. We can thus complete the
proof of Theorem~\ref{theorem2} in the same way as the proof of
Theorem~\ref{theorem1}.
\end{pf*}

\section{Edge universality for deformed Wigner matrices}\label{sectionedgeuniversality}

In this section we prove Theorem~\ref{theedgeuniversalitytheorem}.
Its proof is a combination of Corollary~\ref
{lemmaboundonentroypandirichletform} (bounds on the global Dirichlet
form) and the method
of~\cite{BEY}. In fact, the proof of the edge universality is very
similar to the proof of the bulk universality: we first establish the
edge universality for our model with a small Gaussian component (cf.
Section~\ref{localequilibriummeasures} for the bulk), and then
remove the small Gaussian component using Green function comparison and
a moment matching; cf. Section~\ref{Proofsofmainresults} for the bulk.

\subsection{Edge universality with a small Gaussian component}\label{edgeuniverslaitywithasmallgaussiancomponent}
We mainly follow the exposition in Section~3 of~\cite{BEY}. We
consider the local statistics at the lower edge; the upper edge is
treated in exactly the same way.

\subsubsection{Preliminaries}\label{subsubsectionegdePreliminaries}
Recall the definition of the $\beta$-ensemble $\mu_U$ in~(\ref
{eqnmeasure}) for a given potential $U$ that is $C^4$ and ``regular.''
To study the local statistics at the lower edge, we introduce two
auxiliary measures, $\sigma$ and~$\check\sigma$, on $\digamma
^{(N)}$ as follows. By a shift and a rescaling, we can assume that the
equilibrium density, $\varrho_U$, of $\mu_U$ is supported on
$[0,A_+]$, for some $A_+>0$. Fix a small $\varepsilon_0>0$, and set
%
\begin{eqnarray}
\label{lesigmameasure} \sigma(\dd\bflambda): =\frac
{1}{Z_{\sigma
}}
\mathrm{e}^{-\beta N
\caH_\sigma(\bflambda) } \,\dd\bflambda,
\end{eqnarray}
with
%
\begin{eqnarray}
\caH_\sigma(\bflambda)&=&\caH(\bflambda)+\frac{2}{N}\sum
_{i=1}^N\Theta\bigl(N^{2/3-\varepsilon_0}
\lambda_i\bigr),
\nonumber\\[-8pt]\\[-8pt]\nonumber
\Theta(x) &:=& (x+1)^2\lone(x<-1),
\end{eqnarray}
where $\caH$ is given in~(\ref{lehamiltonianwithU}) and where
$Z_{\sigma}\equiv Z_{\sigma}(\beta)$ is a normalization. Similarly,
we introduce
%
\begin{eqnarray}
\label{lesigmacheckmeasure} \check\sigma(\dd\bflambda)
:=\frac{1}{Z_{\check
\sigma}}
\mathrm{e}^{-\beta N \caH_{\check\sigma}(\bflambda) } \,\dd\bflambda,
\end{eqnarray}
with
%
\begin{eqnarray}
\caH_{\check\sigma}(\bflambda)=\caH(\bflambda)+\frac{1}{N}\sum
_{i=1}^N\Theta\bigl(N^{2/3-\varepsilon_0}
\lambda_i\bigr),
\end{eqnarray}
with $Z_{\check\sigma}\equiv Z_{\check\sigma}(\beta)$ a
normalization. The potential $\Theta$ is added to avoid that the
$(x_i)$ deviate too far to the left, yet its influence on the local
statistics at the edge is negligible; see Lemma~4.1 in~\cite{BEY}.
Below, we choose $\beta=1,2$ depending on the symmetry class of our
original matrix.

Following Section~3 of~\cite{BEY}, we choose a small $\delta>\varepsilon
_0$ and an integer $K$ such that $K\in\llbracket N^\delta,
N^{1-\delta}\rrbracket$. Denote by $I=\llbracket1,K\rrbracket$ the
set of the first $K$ indices. For $\bflambda\in\digamma^{(N)}$, we write
%
\begin{eqnarray}
(\lambda_1,\lambda_2,\ldots,\lambda_N)=(x_1,
\ldots, x_K,y_{K+1},\ldots, y_{N}),
\end{eqnarray}
and
%
\begin{eqnarray}
\bfx=(x_1,\ldots, x_K)\in\digamma^{(K)},
\qquad\bfy=(y_{K+1},\ldots, y_N)\in\digamma^{(N-K)};
\end{eqnarray}
cf.~(\ref{thebflambda}) and~(\ref{bfxnotation}). We further denote
$\mathbf{I}: =(-\infty,y_{K+1}]$. For fixed
$\bfy$, we define
the localized measures $\mu_U^\bfy$, $\sigma^\bfy$ and $\check
\sigma^\bfy$ as in Section~\ref{subsection1ofsection6}. (For
simplicity of notation, we do not indicate the $U$ and $\varepsilon_0$
dependences in the measures $\sigma$, $\check\sigma$.)

We introduce the set of ``good'' boundary conditions
%
\begin{eqnarray}
\caR(\varepsilon_0)\equiv\caR: =\bigl\{\bfy\in
\digamma^{(N-K)}: \llvert y_k-\gamma_k\rrvert
\le N^{-2/3+\varepsilon_0}\check{k},k\notin I\bigr\},
\end{eqnarray}
with $\check k=\min\{k,N-k\}$, where $(\gamma_k)$ denote the
classical locations with respect to the equilibrium density.
With\vspace*{1pt} our choices of $\delta$ and $\varepsilon_0$, we have $y_K-y_1\sim
(K/N)^{2/3}$.

\subsubsection{Comparison of the local measures at the edge}
Fix $t>0$. Recall that we denote by $f_t \mu_G$ the distribution of
$\bolds{\lambda}(t)$ under the flow generated by~(\ref
{dysonbrownianmotion}). As in Section~\ref{localequilibriummeasures}, we
fix $(v_i)$, and condition of the event $\Omega$; see Definition~\ref
{definitionoftheeventomega}. Also recall from~(\ref{definitionofbt})
the definition of the time dependent reference $\beta$-ensemble
$\widehat\psi_t \mu_G$, whose equilibrium density is $\widehat
\varrho_{\mathrm{fc}}(t)$. By a simple shift and a scaling, we may assume, for
fixed $t$, that $\supp\widehat\varrho_{\mathrm{fc}}=[0,\widehat L_+(t)]$ and that
%
\begin{eqnarray}
\label{lesqrtaftershiftscaling} \widehat\varrho_{\mathrm{fc}}(t,x)=\frac{1}{\pi
}\sqrt{x}
\bigl(1+O(x)\bigr),
\end{eqnarray}
as $x\searrow0$. This can easily be checked from the proofs of the
Lemmas~\ref{lemmamfcA} and~\ref{lemmahatmfc} in the \hyperref[app]{Appendix}. For
$\bfy\in\caR$, we then introduce the localized measures $\widehat
\psi_t^{\bfy} \mu_G^{\bfy}$ and $f_t^{\bfy}\mu_G^{\bfy}$ in the
obvious way. For technical reasons, we also use the measures~$\sigma$
and~$\check\sigma$, with the choice $U=\widehat U(t)$. (The
Hamiltonians of the measures $\widehat\psi_t \mu_G$ and $\sigma$,
$\check\sigma$, agree up to the confining potential $\Theta$.)

In a first step, we compare the statistics of $\widehat\psi_t^\bfy
\mu_G^\bfy$ and $\sigma^\bfy$. This is the analogue result to
Proposition~\ref{frozenstatsproximity} above, respectively, to
Lemma~5.4 in~\cite{BEY}.

\begin{lemma}\label{edgefrozenstatproximity}
Let $0< \frak{a}<1/2$. Fix small constants $\delta>\varepsilon_0>0$.
Let $K\in\llbracket N^\delta, N^{1-\delta}\rrbracket$, and let $O$
be an $n$-particle observable. Let $\varepsilon'>0$, and choose~$\tau$
satisfying~$1\gg\tau>N^{-2\frak{a}}$. Then, for any~$t\ge
N^{\varepsilon'}\tau$ and any constant $\frak{c}\in(0,1)$, there is a
set of configurations $\caG(\varepsilon_0)\equiv\caG\subset\caR$, with
%
\begin{eqnarray}
\mathbb{P}^{f_t \mu_G}(\caG)\ge1-\frac{N^{-\frak{c}}}{2},
\end{eqnarray}
such that
%
\begin{eqnarray}
\biggl\llvert\int O(\bfx) \bigl(f_t^{{\bfy}}(\bfx)
\mu_G^{\bfy}(\dd\bfx) -\sigma^{\bfy}(\dd\bfx)
\bigr)\biggr\rrvert\leq C_O K^{1/6}N^{1/3}
N^{\mathfrak{c}-\mathfrak{a}}\tau^{-1},
\end{eqnarray}
$t\ge N^{\varepsilon'}\tau$, for $N$ sufficiently large on~$\Omega$.

Moreover, there is $\upsilon>0$, such that
%
\begin{eqnarray}
\bbP^{f_t^{\bfy}\mu_G^{\bfy}} \bigl( \bigl\{\bigl\llvert x_k - \widehat
\gamma_k(t) \bigr\rrvert< N^{-1+\varepsilon_0}, k\in I \bigr\} \bigr)
\geq1 - \mathrm{e}^{-\upsilon(\varphi_N)^{\xi}},
\end{eqnarray}
$t\ge N^{\varepsilon'}\tau$, for $N$ sufficiently large on~$\Omega$,
with $\xi=A_0\log\log N/2$; see~(\ref{eqxi}).
\end{lemma}

\begin{pf}
We follow the proof of Proposition~\ref{frozenstatsproximity} with
some modifications. First, introduce the density $q_t$ by demanding
\begin{eqnarray*}
q_t \sigma=f_t \mu_G.
\end{eqnarray*}
Then we note that, at the lower edge,
\begin{eqnarray*}
\frac{1}{N}\sum_{k\in I^c}\frac{1}{(x-\widehat\gamma_k(t))^2}\ge c
N^{1/3}/K^{1/3},
\end{eqnarray*}
for\vspace*{1pt} $x\ge-N^{-2/3+\varepsilon_0}$ and $\bfy\in\caR$. We thus have
$\nabla_\bfx^2 \caH^{\bfy}_\sigma(\bfx)\ge cN^{1/3}/K^{1/3}$;
cf.~(\ref{legibsmeasures}) for $\caH_\sigma^{\bfy}$. Hence the
logarithmic Sobolev inequality
\begin{eqnarray*}
S_{\sigma^\bfy}\bigl(q_t^\bfy\bigr)\le C
\frac{K^{1/3}}{N^{1/3}}D_{\sigma
^{\bfy}} \Bigl(\sqrt{q_t^{\bfy}}
\Bigr),
\end{eqnarray*}
with the Dirichlet form $D_{\sigma^{\bfy}}(f)=\frac{1}{\beta N}\sum
_{i\in I}\int\llvert \partial_i f(\bfx)\rrvert ^2\sigma^{\bfy}(\dd\bfx)$
holds. To bound the Dirichlet form, we proceed as in~(\ref
{lewhereweusethedformestimate}),
\begin{eqnarray*}
\E^\sigma D_{\sigma^{\bfy}} \Bigl(\sqrt{q_t^{\bfy}}
\Bigr)&\le& D_{\sigma} (\sqrt{q_t} )
\nonumber
\\
&\le&2D_{\widehat\psi_t \mu_G} (\sqrt{\widehat g_t} )+CN^{4/3}
\sum_{i=1}^N\E^{\widehat\psi_t \mu_G}\bigl\llvert
\Theta'\bigl(N^{2/3-\varepsilon_0}x_i\bigr)\bigr\rrvert
^2
\nonumber
\\
&\le& C\frac{N^{1-2\frak{a}}}{\tau^2}+ \mathrm{e}^{-N^{c}},
\end{eqnarray*}
for some\vspace*{1pt} $c>0$, with $\widehat g_t=f_t/\widehat\psi_t$, where we used
the definitions of $ D_{\sigma}$, $D_{\widehat\psi_t \mu_G}$ to
get the second line. The third line follows from Corollary~\ref
{lemmaboundonentroypandirichletform} and Lemma~\ref
{rigidityfortimedependentbeta}.

To complete the proof, we now follow mutatis mutandis the proof of
Proposition~\ref{frozenstatsproximity}. We leave the details aside.
\end{pf}
Eventually, we are going to apply Theorem~3.3 of~\cite{BEY}, which
shows that the statistics of $\sigma^\bfy$ are universal for most
boundary conditions $\bfy$. In order to apply it, we need the analogue
of Lemma~\ref{frozenbetaparticleproximity} above.

%
\begin{lemma}\label{leedgeproximityofxk}
Under the assumptions of Proposition~\ref{edgefrozenstatproximity}
the following holds. Let $\bfy\in\caG$. Then, assuming that
%
\begin{eqnarray}
\label{leedgeconditionK} K^{1/3}N^{-1/3}N^{2\mathfrak{c}-\frak{a}}
\tau^{-1} \le N^{\varepsilon
_0}K^{-1/3}N^{-2/3},
\end{eqnarray}
we get, for all $k\in I$,
%
\begin{eqnarray}
\bigl\llvert\bbE^{f_t^{\bfy}\mu_G^{\bfy}} x_k - \bbE^{\sigma^\bfy}
x_k \bigr\rrvert&\leq& CN^{\varepsilon_0}K^{-1/3}N^{-2/3},\label
{frozenedgebetaparticleproximityequation}
\end{eqnarray}
for $N$ sufficiently large on $\Omega$.
\end{lemma}

\begin{pf}
Replacing the constant $\tau_K=CK/N$ in the logarithmic Sobolev
inequality~\ref{frozenlog-sobolev} by $CK^{1/3}/N^{1/3}$, we can copy
the proof of Lemma~\ref{frozenbetaparticleproximity} (see also
Lemma~5.5 in~\cite{BEY}) almost word by word.
\end{pf}
From~(\ref{frozenedgebetaparticleproximityequation}), we
immediately get, for $\bfy\in\caG$, the estimate
%
\begin{eqnarray}
\bigl\llvert\E^{\sigma^{\bfy}}\bigl(x_k-\widehat
\gamma_k(t)\bigr)\bigr\rrvert\le C N^{\varepsilon_0}N^{-2/3}k^{-1/3},
\end{eqnarray}
provided that~(\ref{leedgeconditionK}) holds.

\subsubsection{Universality of the localized measures at the edge}
In this subsection, we establish the following result.

\begin{lemma}\label{lelemmaedgeuniveralityaftertimet}
Fix an integer $n>0$. Then for any $1/4>\varkappa$ the following holds
on the event $\Omega$. For any $\delta>0$, there is a constant $\frak
{f}>0$ such that, for $t\ge N^{-\delta}$ and for $\Lambda\subset
\llbracket1,N^\varkappa\rrbracket$
with $\llvert \Lambda\rrvert = n$,
%
\begin{eqnarray}
\label{lelemmaedgeuniveralityaftertimetequation}
&& \bigl\llvert\E^{f_t
\mu_G} O \bigl( \bigl(c_t
N^{ 2/3} j^{1/3}\bigl(\lambda_j-\widehat\gamma
_j(t)\bigr) \bigr)_{j\in\Lambda} \bigr) \nonumber
\\
&&\hspace*{18pt}{} - \E^{\mu_{G}} O
\bigl( \bigl(N^{ 2/3} j^{1/3}(\lambda_j-\gamma
_j) \bigr)_{j\in\Lambda} \bigr) \bigr\rrvert
\\
&&\qquad \leq CN^{-\frak{f}},\nonumber
\end{eqnarray}
where $c_t$ depends only on $\widehat\varrho_{\mathrm{fc}}(t)$. Here,
$(\widehat\gamma_j)$ denote the classical locations with respect the
measure $\widehat\varrho_t$, and $(\gamma_j)$ denote the classical
locations with respect the semicircle law $\varrho_{\mathrm{sc}}$.
\end{lemma}

\begin{pf}
We follow the proof of Lemma~5.1 in~\cite{BEY}. We consider the case
$n=1$ only; the general case is proved in the same way. By a shift and
a scaling, we may assume that $c_t=1$ [see~(\ref
{lesqrtaftershiftscaling})], and we may replace $\widehat\gamma_j(t)$
by the $\gamma
_j$. [Here, we implicitly use that we fixed $(v_i)$ and conditioned on
the event $\Omega$.]

We will need two modifications of the set $\caR(\varepsilon_0)$ of
``good'' boundary conditions. Let $\sigma$, $\check\sigma$ be given
by~(\ref{lesigmameasure}), respectively (\ref{lesigmacheckmeasure})
(with a generic potential $U$). Then set
%
\begin{eqnarray}
\caR^*(\varepsilon_0) &:=&\bigl\{\bfy\in\caR(
\varepsilon_0)\dvtx \forall k\in I, \bigl\llvert\E^{\sigma^{\bfy}}x_k-
\gamma_k\bigr\rrvert\le N^{-2/3+\varepsilon
_0}k^{-1/3},
\nonumber\\[-8pt]\\[-8pt]\nonumber
&&\hspace*{87pt} \mathbb{P}^{\check\sigma^{\bfy}}\bigl(x_1\ge\gamma
_1-N^{-2/3+\varepsilon_0}\bigr)\ge1/2\bigr\}.
\end{eqnarray}
We further need the set
%
\begin{eqnarray}
\caR^\#(\varepsilon_0): =\bigl\{\bfy\in\caR(\varepsilon
_0/3)\dvtx \llvert y_{K+1}-y_{K+2}\rrvert\ge
N^{-2/3-\varepsilon_0}K^{-1/3}\bigr\}.
\end{eqnarray}
While the set~$\caR^*(\varepsilon_0)$ incorporates rigidity estimates in
the sense that $\gamma_k$ is a good approximation\vspace*{1pt} in expectation to
$x_k$ and that $x_1$ is not too much on the left, the set $\caR^\#
(\varepsilon_0)$ incorporates a level repulsion estimate. It has no
counterpart in Section~\ref{localequilibriummeasures} above.

We now choose $\frak{a}=1/2-\delta'$, $\frak{c}=\delta'/2$ and
$\tau=N^{-\delta'}$, for some small $1/12>\delta'>0$. With this
choice, we have for $K\le N^{1/4-6\delta'}$, that
%
\begin{eqnarray}
K^{1/6}N^{1/3} N^{\mathfrak{c}-\mathfrak{a}}\tau^{-1}&\le&
N^{-\varepsilon_0},
\end{eqnarray}
respectively,
%
\begin{eqnarray}
\label{leedgeconditionK2} K^{1/3}N^{-1/3}N^{2\mathfrak{c}-\frak{a}}
\tau^{-1} \le N^{\varepsilon
_0}K^{-1/3}N^{-2/3},
\end{eqnarray}
for a small $\varepsilon_0>0$ (with $\delta'>\varepsilon_0$).

Then, from Lemma~\ref{edgefrozenstatproximity} we have, for $\bfy
\in\caG$,
%
\begin{equation}
\label{leedgeontheway1} \biggl\llvert\int O\bigl(N^{2/3}j(x_j-
\gamma_j)\bigr) \bigl(f_t^{\bfy} \,\dd\mu
_G^{\bfy}-\dd\sigma^{\bfy} \bigr) \biggr\rrvert
\le C_O N^{-\chi}\qquad(j\in\Lambda),
\end{equation}
for some $\chi>0$. Here, the measure $\sigma$ is given by~(\ref
{lesigmameasure}) with the potential $\widehat U(t)$.

Let $\tilde\sigma$ denote the measure given by~(\ref{lesigmameasure})
with the potential $U\equiv0$. For $\tilde\bfy\in\caR
(\varepsilon_0)$ (where the classical locations are taken with respect
the semicircle law), we introduce the localized measure $\tilde
\sigma^{\tilde\bfy}$. We now apply Theorem~3.3 of~\cite{BEY}:
for $\bfy\in\caR^\#(\varepsilon_0)\cap\caR^*(\varepsilon_0)$,
respectively, $\tilde\bfy\in\caR^\#(\varepsilon_0)\cap\caR
^*(\varepsilon_0)$, we have
%
\begin{eqnarray}
\label{leedgeontheway2} \biggl\llvert\int O\bigl(N^{2/3}j(x_j-
\gamma_j)\bigr) \bigl(\dd\sigma^{\bfy}-\dd\tilde
\sigma^{\tilde\bfy} \bigr)\biggr\rrvert\le C_ON^{-\chi},
\end{eqnarray}
for\vspace*{1pt} sufficiently large $N$, by choosing $\chi>0$ sufficiently small.
From Lemma~4.1 of~\cite{BEY}, we know that $\mathbb{P}^{\tilde
\sigma}(\caR^\#(\varepsilon_0)\cap\caR^*(\varepsilon_0))\ge1-N^{-c}$,
for some $c>0$. We further know from Lemma~4.1 of~\cite{BEY} that, for
any bounded observable $O$, $\llvert \E^{\tilde\sigma}O-\E^{\mu
_G}O\rrvert \le C_O\exp(-N^{c})$, $c>0$, where $\mu_G$ denotes the GUE/GOE.
Thus, integrating out the boundary conditions $\tilde\bfy$ and
replacing $\tilde\sigma$ with $\mu_G$, we get from~(\ref
{leedgeontheway1}) and~(\ref{leedgeontheway2}),
%
\begin{eqnarray}
\label{leedgeontheway3} \biggl\llvert\int O\bigl(N^{2/3}j(x_j-
\gamma_j)\bigr) \bigl(f_t^{\bfy} \,\dd\mu
_G^{\bfy}-\dd\mu_G \bigr)\biggr\rrvert\le
C_ON^{-\chi},
\end{eqnarray}
for sufficiently small $\chi>0$, where $\bfy\in\caG(\varepsilon
_0)\cap\caR^\#(\varepsilon_0)\cap\caR^*(\varepsilon_0)$. Once we have
established that
%
\begin{eqnarray}
\label{leedgeprobabilityestimate} \mathbb{P}^{f_t \mu_G} \bigl(\caG
(\varepsilon_0)
\cap\caR^\# (\varepsilon_0)\cap\caR^*(\varepsilon_0) \bigr)
\ge1-N^{-c},
\end{eqnarray}
for some $c>0$, we integrate out the boundary condition $\bfy$
in~(\ref{leedgeontheway3}), and we get~(\ref{lelemmaedgeuniveralityaftertimetequation}) for $n=1$.

To prove~(\ref{leedgeprobabilityestimate}) we follow the two steps
of the proof of~(5.23) in~\cite{BEY}. In a first step, one controls
the probability of $\caR^\#(\varepsilon_0)$ using the rigidity estimates
for $f_t \mu_G$ (see Lemma~\ref{rigidityofeigenvalues}), the level
repulsion estimates for the measure $\sigma^{\bfy}$ in Theorem~3.2
of~\cite{BEY}, Lemma~\ref{edgefrozenstatproximity} and the
condition~(\ref{leedgeconditionK2}). In a second step, one shows
that $\caG(\varepsilon_0)\subset\caR^*(\varepsilon_0)$. This follows
from~(\ref{leedgeproximityofxk}) and the arguments given in the
proof of Lemma~5.1 of~\cite{BEY}. In this way~(\ref
{leedgeprobabilityestimate}) can be established; we leave the details
to the
interested reader.
\end{pf}

\subsection{Removal of the Gaussian component}\label{removaloftheGaussiancomponent}
In this subsection we prove Theorem~\ref{theedgeuniversalitytheorem}.
We use the following version of the Green function comparison
theorem at the edge. It is the counterpart to Theorem~\ref
{thegreenfunctioncomparisontheorem} above.

\begin{theorem}\label{GFTattheedge}
Suppose we have two Wigner matrices $X$ and $Y$ satisfying the
conditions in Definition~\ref{assumptionwigner}. Set $H^X: =V+X$,
$H^Y : =V+Y$; see~(\ref{HXHY}). Denote by $\mathbb{P}^{X}$,
$\mathbb{P}^Y$ the probability distributions of $X$,$Y$. Then
on~$\Omega$ the following holds true. For any $\varepsilon>0$, there is
$\delta>0$ [depending on~$\varepsilon$ and the constants $C_0$,
$\vartheta$ in~(\ref{eqC0})], such that
%
\begin{eqnarray*}
&& \mathbb{P}^{X} \bigl(N^{2/3}(\lambda_1-\widehat
\gamma_1)\le s-N^{-\varepsilon} \bigr)-N^{-\delta}\le
\mathbb{P}^{Y} \bigl(N^{2/3}(\lambda_1-\widehat
\gamma_1)\le s \bigr)
\\
&&\qquad \le\mathbb{P}^X \bigl(N^{2/3}(\lambda_1-
\widehat\gamma_1)\le s+N^{-\varepsilon} \bigr)+N^{-\delta},
\qquad\qquad s\in\R,
\end{eqnarray*}
for $N$ sufficiently large, where $(\widehat\gamma_k)$ denote the
classical locations of the measure $\widehat\varrho_{\mathrm{fc}}\equiv
\widehat\varrho_{\mathrm{fc}}^{\theta}$, with $\theta=1$. Analogous results
hold for the joint distributions of the eigenvalues $\lambda
_{i_1},\lambda_{i_2},\ldots,\lambda_{i_p}$, as long as $\llvert
i_p\rrvert \le
N^{\varepsilon}$.
\end{theorem}

Theorem~\ref{GFTattheedge} is proven exactly in the same way as
Theorem~2.4 of~\cite{EYY} for the Wigner case $V=0$. It suffices to
note that the entries of $V$ are fixed in Theorem~\ref{GFTattheedge}
and that the only input needed in the proof are the local laws
for the Green functions of $H^X$ and $H^Y$, which have been established
in Theorem~\ref{thmstrong} above.

Given Theorem~\ref{GFTattheedge}, we now complete the proof of
Theorem~\ref{theedgeuniversalitytheorem}. Following the arguments
in Section~\ref{leproofofletheoremun}, we construct an auxiliary
Wigner matrix $U$ such that the first two moments of the matrix
%
\begin{eqnarray}
H_t=V+\mathrm{e}^{-t/2}U+\bigl(1-\mathrm{e}^{-t}
\bigr)^{1/2}W'
\end{eqnarray}
with $t=N^{-\delta}$ ($\delta>0$ as in Lemma~\ref
{lelemmaedgeuniveralityaftertimet}, and $W'$ an independent GUE/GOE
matrix) and
the matrix $H=V+W$ match. By Lemma~\ref
{lelemmaedgeuniveralityaftertimet} the edge statistics of $H_t$ are
universal. By
Theorem~\ref{GFTattheedge} the eigenvalue statistics of $H_t$ and
$H$ at the edge agree for large $N$. The existence of such $U$ is
assured by Lemma~\ref{existenceofmatching}. We have thus established
that there is a small $\chi>0$ such that
%
\begin{eqnarray}\label{sames}
&& \bigl\llvert\E^{f_0 \mu_G} O \bigl( \bigl(c_0N^{ 2/3}
j^{1/3}(\lambda_j-\widehat\gamma_j)
\bigr)_{j\in\Lambda} \bigr)
\nonumber\\[-8pt]\\[-8pt]\nonumber
&&\hspace*{8pt}{} - \E^{\mu_{G}} O \bigl( \bigl(N^{ 2/3}
j^{1/3}(\lambda_j-\gamma_j)
\bigr)_{j\in\Lambda} \bigr) \bigr\rrvert
\leq C N^{-\chi},\nonumber
\end{eqnarray}
for $N$ sufficiently large on $\Omega$, where $\mu_G$ is the GUE/GOE.

Finally, we use Assumptions~\ref{assumptionmuVconvergence} and~\ref
{assumptionmuV} as well as a simple moment bound to average over
$\widehat\nu$ (the empirical distribution of $V$) in~(\ref{sames}).
This completes the proof of Theorem~\ref{theedgeuniversalitytheorem}.

\setcounter{theorem}{0}
\setcounter{equation}{0}
%
\begin{appendix}\label{app}
\section*{Appendix}
In this \hyperref[app]{Appendix} we prove the auxiliary results used in Sections~\ref
{sectionlocallaw} and~\ref{sectionbetaensemble}: Lemmas~\ref{lemmamfc},
\ref{lemmahatmfc} and~\ref{superlemma}. We start
with a more extended version of Lemma~\ref{lemmamfc}. Recall
from~(\ref{definitionTheta}) that we denote $\Theta_{\varpi
}=[0,1+\varpi']$, $\varpi'=\varpi/10$. Also recall the definition of
the domain $\caD'$ of the spectral parameter $z$ in~(\ref{definitionofD}).

%
\begin{lemma}\label{lemmamfcA}
Let $\nu$ satisfy Assumption~\ref{assumptionmuV} for some $\varpi
>0$. Then the following holds true for any $\vartheta\in\Theta
_{\varpi}$. There are $L_-^{\vartheta},L_+^{\vartheta}\in\R$, with
$L_-^{\vartheta}<0<L_+^{\vartheta}$, such that $\supp\rho
_{\mathrm{fc}}^{\vartheta}=[L_-^{\vartheta},L_+^{\vartheta}]$, and there is a
constant $C>1$ such that, for all $\vartheta\in\Theta_{\varpi}$,
%
\begin{eqnarray}
\label{thesquarerootA} C^{-1}\sqrt{\kappa_E}\le
\rho_{\mathrm{fc}}^{\vartheta}(E)\le C \sqrt{\kappa_E} \qquad
\bigl(E\in\bigl[L_-^{\vartheta},L_+^{\vartheta}\bigr]\bigr),
\end{eqnarray}
where $\kappa_E$ denotes the distance of $E$ to the endpoints of the
support of $\rho_{\mathrm{fc}}^{\vartheta}$, that is,
%
\begin{eqnarray}
\label{definitionofkappaE} \kappa_E: =\min\bigl\{\bigl
\llvert
E-L_-^{\vartheta}\bigr\rrvert, \bigl\llvert E-L_+^{\vartheta}\bigr
\rrvert
\bigr\}.
\end{eqnarray}
The Stieltjes transform, $m^{\vartheta}_{\mathrm{fc}}$, of $\rho^{\vartheta
}_{\mathrm{fc}}$ has the following properties:
\begin{longlist}[(5)]
\item[(1)] for all $z=E+\ii\eta\in\C^+$,
%
\begin{eqnarray}
\label{behaviorofmfcA} \im m^{\vartheta}_{\mathrm{fc}}(z)\sim\cases{ \sqrt{\kappa+
\eta}, &\quad$E\in[L_-,L_+]$,
\vspace*{3pt}\cr
\displaystyle\frac{\eta}{\sqrt{\kappa+\eta}}, &
\quad$E\in[L_-,L_+]^c$;}
\end{eqnarray}
\item[(2)] there exists a constant $C>1$ such that for all $z\in
\caD'$ and for all $x\in I_{\nu}$,
%
\begin{eqnarray}
\label{stabilityboundA'} C^{-1}\le\bigl\llvert\vartheta x
-z-m^{\vartheta}_{\mathrm{fc}}(z)
\bigr\rrvert\le C;
\end{eqnarray}
\item[(3)] there exists a constant $C>1$ such that for all $z\in
\caD'$,
%
\begin{eqnarray}
C^{-1}\sqrt{\kappa+\eta}\le\biggl\llvert1-\int\frac{\dd\nu
(v)}{(\vartheta v-z-m^{\vartheta}_{\mathrm{fc}}(z))^2}
\biggr\rrvert\le C\sqrt{\kappa+\eta};
\end{eqnarray}
\item[(4)] there are constants $C>1$ and $c_0>0$ such for all
$z=E+\ii\eta\in\caD'$ satisfying $\kappa_E+\eta\le c_0$,
%
\begin{eqnarray}
C^{-1}\le\biggl\llvert\int\frac{\dd\nu(v)}{(\vartheta v-z-m^{\vartheta
}_{\mathrm{fc}}(z))^3}\biggr\rrvert\le C;
\end{eqnarray}
moreover, there is $C>1$, such that for all $z\in\caD'$,
%
\begin{eqnarray}
\biggl\llvert\int\frac{\dd\nu(v)}{(\vartheta v-z-m^{\vartheta
}_{\mathrm{fc}}(z))^3}\biggr\rrvert\le C.
\end{eqnarray}
\end{longlist}
The constants in statements $(1)$--$(4)$ can be chosen uniformly in
$\vartheta\in\Theta_{\varpi}$.
\end{lemma}

\begin{pf}
We follow the proofs in~\cite{S1,LS}. Let $\vartheta\in\Theta
_{\varpi}$. Set $\zeta=z+m_{\mathrm{fc}}^{\vartheta}(z)$, and let
%
\begin{eqnarray}
F(\zeta): =\zeta-\int\frac{\dd\nu(v)}{
\vartheta v-\zeta}\qquad\bigl(\zeta\in
\C^+\bigr).
\end{eqnarray}
Then the functional equation~(\ref{lambdamfc}) is equivalent to
$z=F(\zeta)$. As is argued in~\cite{S1}, a~point $E\in\R$ is inside
the support of the measure $\rho_{\mathrm{fc}}^{\vartheta}$ if and only if
$\zeta_E=E+m_{\mathrm{fc}}^{\vartheta}(E)$ satisfies $\im F(\zeta_E)=0$ and
$\im\zeta_E> 0$. Accordingly, the endpoints of the support are
characterized as the solutions of
%
\begin{eqnarray}
\label{equationforzetapmA} H(\zeta): =\int\frac{\dd\nu
(v)}{(\vartheta
v-\zeta)^2}=1\qquad(
\zeta\in\R).
\end{eqnarray}
Note that $H(\zeta)$ is a continuous function outside $\vartheta
I_{\nu}\equiv\{x\dvtx x=\vartheta y, y\in I_{\nu}\}$ which is
decreasing as $\llvert \zeta\rrvert $ increases. Since $\vartheta\in
\Theta
_{\varpi}=[0,1+\varpi']$, with $\varpi'=\varpi/10$, we obtain from
Assumption~\ref{assumptionmuV} that $H(\zeta)\ge1+\varpi/2$, for
all $\zeta\in\vartheta I_{\nu}$. It thus follows that there are
only two solutions, $\zeta_{\pm}^{\vartheta}\in\R\setminus
\vartheta I_{\nu}$, to $H(\zeta)=1$, $\zeta\in\R$. In
particular, $\zeta_-^{\vartheta}<0$, $\zeta_+^{\vartheta} >0$, and
there is a constant $\mathfrak{g}>0$, depending only on $\nu$, such that
%
\begin{eqnarray}
\label{stabilityboundA} \inf_{\vartheta\in\Theta_{\varpi}}\operatorname
{dist}\bigl(\bigl\{
\zeta_\pm^{\vartheta}\bigr\}, \vartheta I_{\nu}\bigr)\ge
\mathfrak{g}.
\end{eqnarray}

As argued in~\cite{S1,LS}, the set $\gamma: =\{
\zeta\in\C^+\dvtx
\im F(\zeta)=0,\im\zeta> 0 \}$ is, for each fixed $\vartheta\in
\Theta_{\varpi}$, a finite curve in the upper half plane that is the
graph of a continuous function which only connects to the real line at
$\zeta_\pm^{\vartheta}$.

Since $\operatorname{dist}(\{\zeta_{\pm}^{\vartheta}\}, \vartheta
I_{\nu}\}\ge\mathfrak{g}>0$, $F(\zeta)$ is analytic in a
neighborhood of~$\zeta_{\pm}^{\vartheta}$. Thus for $\zeta$ in a
neighborhood of~$\zeta_+^{\vartheta}$, we may write
\begin{eqnarray*}
F(\zeta)=F\bigl(\zeta_+^{\vartheta}\bigr)+F'\bigl(
\zeta_+^{\vartheta}\bigr) \bigl(\zeta-\zeta_+^{\vartheta}\bigr)+
\frac{F''(\zeta_{+}^{\vartheta})}{2}\bigl(\zeta-\zeta_{+}^{\vartheta}
\bigr)^2+\caO\bigl(\bigl(\zeta-\zeta_{+}^{\vartheta}
\bigr)^3\bigr).
\end{eqnarray*}
Note that $F'(\zeta_+^{\vartheta})=0$ by the definition of $\zeta
_{+}^{\vartheta}$. Moreover, we know that $\im F(\zeta)=0$,
for~$\zeta$ in a real neighborhood of~$\zeta_{+}^{\vartheta}$, but
we also have $\im F(\zeta)=0$, for~$\zeta\in\gamma\cup\bar
\gamma$. Thus $F''(\zeta_{+}^{\vartheta})\neq0$. We can therefore
invert $F(\zeta)=z$ in a neighborhood of~$\zeta_{+}$ to obtain
%
\begin{eqnarray}
\label{invertedequation} \zeta(z)=F^{(-1)}(z)=\zeta_{+}^{\vartheta
}+c_{+}^{\vartheta}
\sqrt{z-L_+^{\vartheta}} \Bigl(1+\caA_{+}^{\vartheta} \Bigl(
\sqrt{z-L_{+}^{\vartheta}} \Bigr) \Bigr)
\end{eqnarray}
[with the convention $\im F^{(-1)}(z)\ge0$], where $L_+^{\vartheta}$
is defined by $\zeta_{+}^{\vartheta}=L_+^{\vartheta
}+m_{\mathrm{fc}}(L_+^{\vartheta})$. Here, $c^{\vartheta}_+>0$ is a real
constant, and $\caA_+^{\vartheta}$ is an analytic function that is
real-valued on the real line and that satisfies $\caA^{\vartheta
}_{+}(0)=0$. Recalling that $\zeta(z)=z+m_{\mathrm{fc}}^{\vartheta}(z)$ and
taking the limit~$\eta\to0$ we obtain~(\ref{thesquarerootA}), for
fixed $\vartheta$. To achieve uniformity in~$\vartheta$, we use the
(uniform) stability bound~(\ref{stabilityboundA}) and the
(pointwise) positivity of $\llvert F''(\zeta_+^{\vartheta})\rrvert $: we
differentiate~(\ref{equationforzetapmA}) with respect to~$\vartheta
$ and observe that $\partial_{\vartheta}H(\zeta,\vartheta)\mid_{\zeta
=\zeta_{+}^{\vartheta}}\neq0$, for all $\vartheta\in\Theta
_{\varpi}$, since $F''(\zeta_{+}^{\vartheta})\neq0$. Thus by the
implicit function theorem, $\zeta_
{+}^{\vartheta}$ is a $C^1$ function of $\vartheta\in\Theta_{\varpi
}$. Next, we observe that $F''(\zeta)$ is an analytic function
of~$\zeta$, for~$\zeta$ away from~$\vartheta I_{\nu}$. Thus,
using once more~(\ref{stabilityboundA}), we can bound $\llvert
F''(\zeta
_{+}^{\vartheta})\rrvert \ge c$, for some $c>0$, uniformly in
$\vartheta\in
\Theta_{\varpi}$. In fact, $F^{(n)}(\zeta_{+}^{\vartheta})$, $n\in
\N$ are all continuous functions of $\vartheta\in\Theta_{\varpi}$,
and we can bound them uniformly in $\vartheta$ for each $n\in\N$.
Repeating the same argument for $\zeta$ close to $\zeta_-^{\vartheta
}$, we complete the proof of~(\ref{thesquarerootA}).

Statement $(2)$ follows from~(\ref{stabilityboundA}) for $z$ close
to the edges. For $z$ away from the edges, Assumption~\ref
{assumptionmuV} assures that the curve $\gamma$ stays away from the
real line
for all $\vartheta\in\Theta_{\varpi}$ as is readily checked. This
implies the stability bound for that region.

For the proofs of the remaining statements, we refer to the Appendix
of~\cite{LS}.
\end{pf}

Next, we prove Lemma~\ref{lemmahatmfc}.

\begin{pf*}{Proof of Lemma~\ref{lemmahatmfc}}
It follows from Assumption~\ref{assumptionmuV} that on $\Omega$ for
all $N$ sufficiently large,
%
\begin{eqnarray}
\label{assumptionsforproof} \inf_{x\in I_{\widehat\nu}}\frac{1}{N}\sum
_{i=1}^N\frac
{1}{(\vartheta v_i-x)^2}\ge1+\varpi/2,
\end{eqnarray}
for all $\vartheta\in\Theta_{\varpi}=[0,1+\varpi/10]$. The
analogous statements of Lemma~\ref{lemmamfcA}, holding on $\Omega$
for $N$ sufficiently large, follow in the same way as in the proof of
that lemma. To get uniformity in $N$, it suffices to check that the
analogous expression to~(\ref{stabilityboundA}) holds uniformly in
$N$, for $N$\vspace*{1pt} sufficiency large: by~(\ref{assumptionsforproof}) there
are two real solutions $\widehat\zeta_{\pm}^{\vartheta}$ to
$\widehat H(\zeta): =\frac{1}{N}\sum_{i=1}^N\frac
{1}{(\vartheta
v_i-\zeta)^2}=1$ that both lie outside of the interval~$\vartheta
I_{\widehat\nu}$. Thus~(\ref{definitionofomegaV}) and~(\ref
{definitionofomegaVIII}) imply that
%
\begin{eqnarray}
\label{stabilityboundB} \inf_{\vartheta\in\Theta_{\varpi}}\operatorname
{dist}\bigl(\bigl\{\widehat
\zeta_\pm^{\vartheta}\bigr\}, \vartheta I_{\widehat\nu}\bigr)\ge
\mathfrak{g}/2,
\end{eqnarray}
on $\Omega$ for all $N$ sufficiently large. Then we can bound
\begin{eqnarray*}
\widehat F''(\zeta)=-\frac{2}{N}\sum
_{i=1}^N \frac{1}{(\vartheta
v_i-\zeta)^3},
\end{eqnarray*}
evaluated at $\widehat\zeta_{\pm}^{\vartheta}$, uniformly below
in~$\vartheta$ and~$N$, for~$N$ sufficiently large, implying the
uniformity in $N$ of the constants in statements (1)--(4).

Next we prove~(\ref{minilocallaw}). For simplicity we drop
$\vartheta$ from the notation and work on~$\Omega$. As above, set
$\zeta=z+m_{\mathrm{fc}}(z)$ and $\widehat\zeta=z+\widehat m_{\mathrm{fc}}(z)$. From\vspace*{1pt}
the definitions of $F$, $\widehat F$ and equations~(\ref
{lambdamfc}),~(\ref{hatmfc}), we have $\widehat F(\widehat\zeta)=F(\zeta
)=z$, for all $z\in\caD'$. Using the stability bound~(\ref
{stabilityboundB}) and equation~(\ref{definitionofomegaV}) in the
definition of $\Omega$, we get, assuming that $\llvert \widehat\zeta
-\zeta
\rrvert \ll1$,
%
\begin{eqnarray}
\label{miniperturbationexpansion}
&& \bigl[F'(\zeta)+\caO\bigl(N^{-\alpha
_0}\bigr) \bigr] (\widehat\zeta-\zeta)+\frac{F''(\zeta)}{2} (\widehat\zeta-\zeta
)^2
\nonumber\\[-8pt]\\[-8pt]\nonumber
&&\qquad =
o(1) (\widehat\zeta-\zeta)^2+\caO\bigl(N^{-\alpha_0}\bigr),
\end{eqnarray}
uniformly in $\vartheta\in\Theta_{\varpi}$, for all $z\in\caD'$.
From Lemma~\ref{lemmamfcA}, we get $F'(\zeta)\sim\sqrt{\kappa
+\eta}$ and $F''(\zeta)\le C$, for all $z\in\caD'$. We abbreviate
$\Lambda: =\llvert \widehat\zeta-\zeta\rrvert $ in
the following.

We first consider $z=E+\ii\eta\in\caD'$, such that $\kappa_E+\eta
>N^{-\varepsilon}$, for some small $\varepsilon>0$ (with $\varepsilon
<\alpha
_0$). Here $\kappa_E$ is defined in~(\ref{definitionofkappaE}). For
such $z$ we obtain from~(\ref{miniperturbationexpansion}) that $
\Lambda\le CN^{\varepsilon}(\Lambda^2+N^{-\alpha_0})$. Thus either
$\Lambda\le C_0 N^{\varepsilon} N^{-\alpha_0}$ or $C_0N^{-\varepsilon
}\le
\Lambda$, for some constant $C_0$. We now show that $\llvert \Lambda
\rrvert \le
C_0 N^{\varepsilon}N^{-\alpha_0}$, for all $z\in\caD'$ such that
$\kappa_E+\eta\ge N^{-\varepsilon}$. For $z\in\caD'$ with~$\eta=2$,
we have
\begin{eqnarray*}
\widehat\zeta(z)-\zeta(z)=\frac{1}{N}\sum_{i=1}^N
\frac{\widehat
\zeta(z)-\zeta(z)}{(\vartheta v_i-\widehat\zeta(z) )(\vartheta
v_i-\zeta(z))} +\caO\bigl(N^{-\alpha_0}\bigr),
\end{eqnarray*}
where we use~(\ref{definitionofomegaV}). Since $\eta=2$ and $\im
\widehat\zeta, \im\zeta\ge\eta$, we obtain $\Lambda\le\frac
{1}{4}\Lambda+\caO(N^{-\alpha_0})$, that is, $\Lambda(z)\le C
N^{-\alpha_0}$, for $\eta=2$. To extend the conclusion to all~$\eta
$, we use the Lipschitz continuity of $\widehat\zeta(z)$,
respectively, $\zeta(z)$. Differentiating $z= F(\zeta)$, with respect
to $z$ we obtain $\partial_z\zeta=( F'(\zeta))^{-1}$. Thus using
property~$(2)$ of Lemma~\ref{lemmamfcA}, we infer that the Lipschitz
constant of $\zeta(z)$ is, for $z\in\caD'$ satisfying $\kappa
_E+\eta>N^{-\varepsilon}$, bounded above by $N^{\varepsilon/2}$. The same
conclusion also holds for $\widehat\zeta(z)$. Bootstrapping, we obtain
%
\begin{eqnarray}
\label{minilocallawbis} \bigl\llvert\widehat\zeta(z)-\zeta(z)\bigr
\rrvert\le
CN^{\varepsilon}N^{-\alpha_0},
\end{eqnarray}
on\vspace*{1pt} $\Omega$ for $N$ sufficiently large, for all $z\in\caD'$
satisfying~$\kappa_E+\eta> N^{-\varepsilon}$.

In\vspace*{1pt} order to control $\widehat\zeta(z)-\zeta(z)$ for $z=E+\ii\eta
\in\caD'$ with $\kappa_E+\eta\le N^{-\varepsilon}$, $\varepsilon>0$, we
first derive the estimate $\llvert \widehat{L}_\pm^{\vartheta}-L_\pm
^{\vartheta}\rrvert \le C N^{-\alpha_0}$, for some $c>0$, on~$\Omega$. We
recall that~$\widehat{L}_\pm$, respectively, $L_\pm$, are obtained
through the relations
\begin{eqnarray*}
\frac{1}{N}\sum_{i=1}^N
\frac{1}{(\vartheta v_i-\widehat\zeta_\pm
)^2}=1,\qquad\int\frac{\dd\nu(v)}{(\vartheta v-\zeta_\pm
)^2}=1.
\end{eqnarray*}
Then a similar argument as given above shows that $\llvert \widehat
\zeta_\pm
-\zeta_\pm\rrvert \le C N^{-\alpha_0}$ and $\llvert \widehat L_\pm
-L_\pm\rrvert \le
CN^{-\alpha_0}$ on $\Omega$, $N$ sufficiently large. We refer to
Section~4.3 in~\cite{LS2} for details.

Second, following the arguments in the proof of Lemma~\ref{lemmamfcA},
we may write, for $\widehat\zeta$ and $\zeta$ in a neighborhood
of $\zeta_\pm$,
%
\begin{eqnarray}
\label{minilocallawtres} \widehat\zeta(z)-\widehat\zeta_\pm&=&\widehat
c_\pm\sqrt{z-\widehat L_\pm}\bigl(1+\caO(z-\widehat
L_\pm)\bigr),
\nonumber\\[-8pt]\\[-8pt]\nonumber
\zeta(z)-\zeta_\pm&=& c_+\sqrt{z- L_\pm}\bigl(1+\caO(z-
L_\pm)\bigr).
\end{eqnarray}
We therefore get $\llvert \widehat\zeta(z)-\zeta(z)\rrvert \le C\sqrt
{\kappa
_E+\eta}+CN^{-\alpha_0/2}$. Note that the constants can be chosen
uniformly in $\vartheta\in\Theta_{\varpi}$. Choosing, for example,
$\varepsilon=\alpha_0/4$, we get from~(\ref{minilocallawbis})
and~(\ref{minilocallawtres}) the desired inequality~(\ref{minilocallaw}).
\end{pf*}

We now move on to the construction of the potentials $\widehat U$ and
$U$. We first record the following corollary of Lemma~\ref{lemmamfcA}.
Set $\mathrm{B}_{r}(p): =\{z\in\C: \llvert z-p\rrvert
<r\}
$. Recall the
conventions in~(\ref{abuseofnotationI}) and the definition of
$\kappa_E$ in~(\ref{definitionofkappaE}).

%
\begin{corollary}\label{remarkaboutrealandimaginarypartofmfc}
Under the assumptions of Lemma~\ref{lemmamfc} there are constants
$c_+^{\vartheta}$, $r_+>0$, such that for any $E\in\mathrm
{B}_{r_+}(L_+^{\vartheta})\cap\R$,
%
\begin{eqnarray}
\label{realandimaginarypartofmfcI} \im m_{\mathrm{fc}}^{\vartheta}(E)=\cases{
\displaystyle
\sqrt{\kappa_E}\bigl(c_{+}^{\vartheta}+
\caB_{+}^{\vartheta
}({-\kappa_E})\bigr), &\quad$E\le
L_+^{\vartheta}$,
\vspace*{3pt}\cr
0, &\quad$E\ge L_+^{\vartheta}$,}
\end{eqnarray}
and
%
\begin{eqnarray}
\label{realandimaginarypartofmfcII} \re m_{\mathrm{fc}}^{\vartheta}(E)=\cases{
\caC_+^{\vartheta}(-\kappa_E), &\quad$E\le L_+^{\vartheta}$,
\vspace*{3pt}\cr
\displaystyle\sqrt{\kappa_E}\bigl(c_{+}^{\vartheta}+
\caB_{+}^{\vartheta
}({\kappa_E})\bigr)+
\caC_+^{\vartheta}(\kappa_E), &\quad$E\ge L_+^{\vartheta}$,}
\end{eqnarray}
where $\caB_+^{\vartheta},\caC_+^{\vartheta}$ are analytic
functions on $\mathrm{B}_{r_+}(0)$ that are real-valued on $\R$ and
that satisfy $\caB_+^{\vartheta}(0)=0$, $c_+^{\vartheta}+\caB
_+^{\vartheta}> 0$, respectively, ${\caC_+^{\vartheta}}<0$, on
${\mathrm{B}_{r_+}(0)}\cap\R$. Moreover, for all $z\in\mathrm
{B}_{r_+}(L_+^{\vartheta})$, the functions $\caB_+^{\vartheta},\caC
_+^{\vartheta}$, respectively, $\im m_{\mathrm{fc}}^{\vartheta}$, $\re
m_{\mathrm{fc}}^{\vartheta}$, are continuous in $\vartheta\in\Theta_{\varpi}$.

Similar statements hold at the lower edge $L_{-}^\vartheta$.
\end{corollary}

\begin{pf}
Fix $\vartheta\in\Theta_{\varpi}$. As argued in the proof of
Lemma~\ref{lemmamfcA}, the function~$F(\zeta)$ can locally be
inverted around $\zeta_{\pm}^{\vartheta}$; see~(\ref{invertedequation})
above. Thus for $\zeta$ in a neighborhood of $\zeta
_+^{\vartheta}$, we may write
\begin{eqnarray*}
m^{\vartheta}_{\mathrm{fc}}(z)&=&F^{-1}(z)-z = \zeta_{+}^{\vartheta
}-z+c_{+}^{\vartheta}
\sqrt{z-L_+^{\vartheta}} \Bigl(1+\caA_{+}^{\vartheta} \Bigl(
\sqrt{z-L_{+}^{\vartheta}} \Bigr) \Bigr)
\\
&=&c_{+}^{\vartheta}\sqrt{z-L_+^{\vartheta}} \bigl(1+\caB
_+^{\vartheta}\bigl(z-L_+^{\vartheta}\bigr) \bigr)+\caC_{+}^{\vartheta
}
\bigl(z-L_{+}^{\vartheta}\bigr),
\end{eqnarray*}
for $z$ in $\mathrm{B}_{r}(L_+^{\vartheta})$, for some $r>0$, where
$\caB_{+}^{\vartheta}$ and $\caC_{+}^{\vartheta}$ are analytic in a
neighborhood of zero and real-valued on the real line, since $\im
F^{-1}(E)=0$, for $E\in[L_-^{\vartheta},L_+^{\vartheta}]^c$.
Equations~(\ref{realandimaginarypartofmfcI}) and~(\ref
{realandimaginarypartofmfcII}) follow. From the proof of Lemma~\ref
{lemmamfcA}, it is immediate that $c^{\vartheta}_+>0$. Thus
$c^{\vartheta}_++\caB^{\vartheta}_+>0$ in a real neighborhood of
zero. Since $x-L_+^{\vartheta}-m_{\mathrm{fc}}^{\vartheta}(L_+^{\vartheta})<0
$, for all $x\in\vartheta I_{\nu}$, we must have $\caC
_+^{\vartheta}<0$ in a real neighborhood of zero. Since $F(z)$ is
analytic on $\mathrm{B}_{r}(L_+^{\vartheta})$, for all $\vartheta\in
\Theta_{\varpi}$, and since $\zeta_{+}^{\vartheta}$ is a $C^1$
function of $\vartheta$, the functions $\caB_{+}^{\vartheta}$ and
$\caC_{+}^{\vartheta}$ are $C^1$ in $\vartheta\in\Theta_{\varpi
}$. Then it
is clear from~(\ref{stabilityboundA}) that we can choose $r>0$
uniformly in $\vartheta\in\Theta_{\varpi}$. The same arguments
apply for $\zeta$ close to $\zeta_-^{\vartheta}$.
\end{pf}

The analogous result to Corollary~\ref
{remarkaboutrealandimaginarypartofmfc} is stated next. Recall the
notation $\widehat\kappa
_E: =\min\{\llvert E-\widehat{L}_-^{\vartheta}\rrvert,
\llvert E-\widehat
{L}_+^{\vartheta}\rrvert \}$.

%
\begin{corollary}\label{remarkaboutrealandimaginarypartofhatmfc}
Under the assumptions of Lemma~\ref{lemmahatmfc} the following holds
on $\Omega$, for~$N$ sufficiently large. There are constants $\widehat
c_+^{\vartheta}$, $r_+'$, with $r_+\ge r_+'>0$, such that for any
$E\in\mathrm{B}_{r_+'}(L_+^{\vartheta})\cap\R$,
%
\begin{eqnarray}
\im\widehat m_{\mathrm{fc}}^{\vartheta}(E)=\cases{ \displaystyle\sqrt{
\widehat\kappa_E} \bigl(\widehat c_{+}^{\vartheta}+
\widehat\caB_{+}^{\vartheta}({-\widehat\kappa_E})
\bigr), &\quad$E\le\widehat L_+^{\vartheta}$,
\vspace*{3pt}\cr
0, &\quad$E\ge\widehat
L_+^{\vartheta}$,}
\end{eqnarray}
and
%
\begin{eqnarray}
\re\widehat m_{\mathrm{fc}}^{\vartheta}(E)=\cases{ \displaystyle\widehat
\caC_+^{\vartheta}(-\widehat\kappa_E), &\quad$E\le\widehat
L_+^{\vartheta}$,
\vspace*{3pt}\cr
\displaystyle\sqrt{\widehat\kappa_E}
\bigl(\widehat c_+^{\vartheta
}+\widehat\caB_{+}^{\vartheta} ({
\widehat\kappa_E}) \bigr)+ \widehat\caC_+^{\vartheta}(\widehat
\kappa_E), &\quad$E\ge\widehat L_+^{\vartheta}$,}
\end{eqnarray}
where $ \widehat\caB_+^{\vartheta}, \widehat\caC_+^{\vartheta}$
are analytic functions on $\mathrm{B}_{r_+'}(0)$ that are real-valued
on $\R$ and that satisfy
$\widehat\caB_+^{\vartheta}(0)=0$ and~$\widehat c_+^{\vartheta
}+\widehat\caB_+^{\vartheta}> 0$, respectively, ${\widehat\caC
_+^{\vartheta}}<0$, on ${\mathrm{B}_{r'_+}(0)}\cap\R$. Moreover,
the constant $r_+'$ can be chosen independent of $\vartheta\in\Theta
_{\varpi}$ and~$N$, for~$N$ sufficiently large.

Further, the functions $\widehat\caB_+^{\vartheta}, \widehat\caC
_+^{\vartheta}$, respectively, $\im\widehat m_{\mathrm{fc}}^{\vartheta}$,
$\re\widehat m_{\mathrm{fc}}^{\vartheta}$, are continuous functions in
$\vartheta\in\Theta_{\delta}$, for all $z\in\mathrm
{B}_{r_+'}(L_+^{\vartheta})$. There is $c>0$, such that
%
\begin{equation}
\label{convergenceofBandCfunction} \bigl\llvert\widehat\caB
_+^{\vartheta}(z)-
\caB_+^{\vartheta}(z)\bigr\rrvert\le N^{-c\alpha_0/2},\qquad\bigl
\llvert
\widehat\caC_+^{\vartheta}(z)- \caC^{\vartheta}_+(z)\bigr\rrvert\le
N^{-c\alpha_0/2},
\end{equation}
for all $ z\in\mathrm{B}_{r_+'}(L_+^{\vartheta})$ and all
$\vartheta\in\Theta_{\varpi}$, on $\Omega$ for $N$ sufficiently large.

Similar statements hold at the lower edge $\widehat L_{-}^\vartheta$.
\end{corollary}

\begin{pf}
Corollary~\ref{remarkaboutrealandimaginarypartofhatmfc} is
proven in the same way as Corollary~\ref
{remarkaboutrealandimaginarypartofmfc}. The only things to be checked
are that $r_\pm
'>0$ can be chosen uniformly in~$N$,~$N$ sufficiently large, and the
bounds in~(\ref{convergenceofBandCfunction}). The former
statement is an immediate consequence of the stability bound~(\ref
{stabilityboundB}). The latter follows from $z=F(\zeta)=\widehat
F(\widehat\zeta)$, with $\zeta=z+m_{\mathrm{fc}}^{\vartheta}(z)$ and
$\widehat\zeta=z+\widehat m_{\mathrm{fc}}^{\vartheta}(z)$. Then using~(\ref
{definitionofomegaV}), the stability bound~(\ref{stabilityboundB})
and the uniform lower bound on $F''(\zeta_\pm^{\vartheta})$, it is
straightforward to derive estimate~(\ref{convergenceofBandCfunction})
from~(\ref{minilocallaw}).
\end{pf}

Next we prove Lemma~\ref{superlemma}. Recall from~(\ref
{equationforvartheta}) that we chose $\vartheta\equiv\vartheta
(t): =\mathrm
{e}^{-(t-t_0)/2}$.

\begin{pf*}{Proof of Lemma~\ref{superlemma}}
For $c>0$ and a measure $\omega$ on $\R$, we define $\supp_c\omega
: =\supp\omega+[-c,c]$. Recall the constants
$r'_\pm>0$ of
Corollary~\ref{remarkaboutrealandimaginarypartofhatmfc}. Set
$s: =\min\{r'_-, r'_+\}/2$.

We specify the potentials $\widehat U$ and $U$ through their spatial
derivatives $\widehat U'$ and $U'$. For $t\ge0$, we set
\begin{eqnarray*}
\widehat U'(t,x)+x &:=& -2
\dashint_\R\frac{\widehat
\rho
_{\mathrm{fc}}(t,y)}{y-x} \,\dd y,\qquad U'(t,x)+x
: =-2\dashint_\R\frac
{\rho_{\mathrm{fc}}(t,y)}{y-x} \,\dd y,
\end{eqnarray*}
for $x\in\supp\widehat\rho_{\mathrm{fc}}(t)$, respectively, $x\in\supp\rho
_{\mathrm{fc}}(t)$.

For $x\in\R$ satisfying $\llvert x- L_\pm(t)\rrvert \le s$, where
$L_\pm(t)$
denote the endpoints of the support of the measure $\rho_{\mathrm{fc}}(t)$, we set
%
\begin{eqnarray}\label{definitionwidehatUIIappendix}
\widehat U'(t,x)+x
&:=& -2\widehat
\caC^{\vartheta
}_\pm(\widehat k_{\pm}),\qquad\widehat
k_{\pm}\equiv x-\widehat L_\pm(t),
\nonumber\\[-8pt]\\[-8pt]\nonumber
U'(t,x)+x &:=& -2 \caC^{\vartheta}_\pm(k_{\pm}),\qquad
k_{\pm
}\equiv x-L_\pm(t),
\end{eqnarray}
where~$\widehat\caC^{\vartheta}_\pm$ are the functions appearing in
Corollary~\ref{remarkaboutrealandimaginarypartofhatmfc}
with~$\vartheta\equiv\vartheta(t)$, and~$\caC^{\vartheta}_\pm$
are the functions appearing in Corollary~\ref
{remarkaboutrealandimaginarypartofmfc} with $\vartheta\equiv\vartheta
(t)$. From
Lemma~\ref{lemmamfcA}, Corollaries~\ref{remarkaboutrealandimaginarypartofmfc} and~\ref
{remarkaboutrealandimaginarypartofhatmfc}, we conclude that $\widehat U'(t,x)$,
respectively, $U'(t,x)$ are well defined for $x\in\supp_s\rho
_{\mathrm{fc}}(t)$, $t\ge0$, where $s=\min\{ r'_-, r'_+\}/2$.

For $x\notin\supp_{s}\rho_{\mathrm{fc}}(t)$, we define $U'$ as a $C^3$
extension in $x$ such that: (1)~$ U^{(n)}(t,x)$, $\partial_t
U^{(n)}(t,x)$, $n\in\llbracket1,3\rrbracket$, are\vspace*{1pt} continuous
in~$t$;~(2) for all $t\ge0$ and for all $x\notin\supp_s\rho
_{\mathrm{fc}}(t)$, $\llvert U'(t,x)+x\rrvert >\llvert 2\re
m_{\mathrm{fc}}(t,x)\rrvert $ and $\widehat U''(t,x)\ge
-C_U$, for some constant $C_U\ge0$; (3) $U'(t,x)+x\sim x$ for all
$t\ge0$, as $\llvert x\rrvert \to\infty$.\vspace*{2pt} Similarly, we define
$\widehat U(t,x)$
as $C^3$ extensions such that: (1) $\widehat U^{(n)}(t,x)$, $\partial
_t\widehat U^{(n)}(t,x)$, $n\in\llbracket1,3\rrbracket$, are
continuous in $t$;~(2) there is $c>0$ such that $\sup_{t\ge
0}\llvert \widehat U^{(n)}(t,x)-U^{(n)}(t,x)\rrvert \le N^{-c\alpha
_0/2}$, $n\in
\llbracket1,3\rrbracket$, for $N$ sufficiently large on $\Omega$.

We next show that the potential $ U'(t,x)$ is a $C^3$ function in $x$.
For simplicity, we often drop the $t$-dependence from the notation. Let
$\zeta=z+ m_{\mathrm{fc}}(z)$, and recall from the proof of Lemma~\ref
{lemmamfcA} that $\zeta(z)$ satisfies $\zeta(z) = F^{(-1)}(z)$, where
$F(\zeta)=\zeta-\int\frac{\dd\nu(v)}{(\vartheta v-\zeta)}$.
Thus, to prove regularity of $ U'(t,x)$ in $x$ in the support of the
measure $\rho_{\mathrm{fc}}(t)$, it suffices to show that $F'(\zeta) \neq0$
on the curve $\gamma\cap\C^+$ where $\im F = 0$. Recall that on
$\gamma$ we have
%
\begin{equation}
\label{appendixA20} \tilde{H}(\zeta): =\int\frac{\dd\nu
(v)}{\llvert \vartheta v- \zeta
\rrvert ^2} =
1,
\end{equation}
where $\vartheta\equiv\vartheta(t)$. On the other hand, we have
\[
\re{F}'(\zeta)=1-\int\frac{(\vartheta v-\re\zeta)^2-(\im\zeta
)^2}{\llvert \vartheta v-\zeta\rrvert ^4} \,\dd\nu(v).
\]
Thus, on the curve $\gamma$,
\begin{eqnarray*}
\re{F}'(\zeta)&=&\int\frac{\dd\nu(v)}{\llvert \vartheta v- \zeta
\rrvert ^2}-\int\frac{(\vartheta v-\re\zeta)^2-(\im\zeta)^2}{\llvert
\vartheta
v-\zeta\rrvert ^4}
\,\dd\nu(v)
\nonumber
\\
&=&\int\frac{2(\im\zeta
)^2}{\llvert \vartheta v-\zeta\rrvert ^4} \,\dd\nu(v).
\end{eqnarray*}
From~(\ref{sumrule}) we get
%
\begin{equation}
\label{appendixA24} \int\frac{\dd\nu(v)}{\llvert \vartheta v-\zeta
\rrvert ^4}\ge\biggl(\int\frac
{\dd\nu(v)}{\llvert \vartheta v-\zeta\rrvert ^2}
\biggr)^2=1,
\end{equation}
on $\gamma$. Since $F'\neq0$ on $\gamma$, the inverse function
theorem implies that the real part of ${m}_{\mathrm{fc}}(t,x)$ is a smooth
function in the interior of $\supp\rho_{\mathrm{fc}}(t)$, whose derivatives
are continuous in $t$. For $x\in\mathrm{B}_s(L_{\pm}^\vartheta)$,
we already showed in Lemma~\ref{remarkaboutrealandimaginarypartofmfc}
that $\caC_{\pm}^{\vartheta}(x)$ is a smooth function, whose
derivatives are continuous in $t$. Thus we have shown that $U'(t,x)$ is
smooth in $\supp_s\rho_{\mathrm{fc}}(t)$. Outside $\supp_s\rho_{\mathrm{fc}}(t)$,
$U'(t,x)$ is manifestly $C^3$ by definition: it is a $C^3$ extension of
the functions $\caC_\pm(t)$. Thus $\R\ni x\mapsto U'(t,x),
\partial_t U'(t,x)$ are $C^3$ functions for all $t\ge0$.

Clearly, we can bound the derivatives $U^{(n)}(t,x)$, $\partial_t
U^{(n)}(t,x)$, $n\in\llbracket1,3\rrbracket$, uniformly on compact
sets. It is also immediate that $U^{(n)}(t,x)$ are continuous functions
in $t\ge0$. Thus we can bound $U^{(n)}$ uniformly in $t$ and uniformly
in~$x$ on compact sets, for $n\in\llbracket1,3\rrbracket$. For $x\in
\supp_s\rho_{\mathrm{fc}}(t)$, we have $U''(t,x)\ge-C$, for some $C\ge0$.
For $x\notin\supp_s\rho_{\mathrm{fc}}(t)$, a similar bound holds true by
construction. Thus $U'(t,x)$ satisfies~(\ref
{assumption1forbetauniversality}) uniformly in $t\ge0$. Further, since
$U'(t,x)+x\sim x$,
as $\llvert x\rrvert \to\infty$,~(\ref
{assumption2forbetauniversality}) also
holds uniformly in $t\ge0$.

On $\Omega$, we can extend the reasoning above to $\widehat U'(t,x)$,
$\partial_t \widehat U'(t,x)$, for~$N$ sufficiently large. For
example, the arguments in~(\ref{appendixA20})--(\ref{appendixA24})
can be extended to the finite $N$ case by using~(\ref
{definitionofomegaV}) and Lemma~\ref{lemmahatmfc}. Let again $s\equiv
\min\{
r'_-, r'_+\}/2$. Then for~$x\in\supp_s\rho_{\mathrm{fc}}(t)$ we have by
Lemma~\ref{lemmahatmfc} that $\llvert \widehat m_{\mathrm{fc}}(t,x+\ii\eta
)-m_{\mathrm{fc}}(t,x+\ii\eta)\rrvert \le N^{-c\alpha_0}$, for some $c>0$,
on~$\Omega$ for all~$\eta\ge0$ and all $t\ge0$. Together
with~(\ref{convergenceofBandCfunction}) we can conclude that
$\llvert \widehat U'(t,x)-U'(t,x)\rrvert \le N^{-c\alpha_0/2}$
on~$\Omega$, for
$x\in\supp_s\rho_{\mathrm{fc}}(t)$. We also have $\llvert \partial_x\widehat
m_{\mathrm{fc}}(t,x+\ii\eta)-\partial_x m_{\mathrm{fc}}(t,x+\ii\eta)\rrvert \le
CN^{-c\alpha_0}$, for $x$ satisfying $\min\{\llvert x-L_+\rrvert
,\llvert x-L_-\rrvert \}\ge s$,
as can be checked as in the proof of Lemma~\ref{lemmahatmfc}. Hence,
combining this last statement with the regularity of~$ \widehat
C_\pm^{\vartheta}$ claimed in
Lemma~\ref{remarkaboutrealandimaginarypartofhatmfc}, we have
$\llvert \widehat U''(t,x)-U''(t,x)\rrvert \le N^{-c\alpha_0}$,
for~$x\in\supp
_s\rho_{\mathrm{fc}}(t)$, $t\ge0$, on~$\Omega$ for~$N$ sufficiently large.
This conclusion can be extended to arbitrary~$\widehat U^{(n)}$.
Similarly, one checks that~$U^{(n)}(t,x)$, $n\in\llbracket
1,3\rrbracket$ are continuous functions of~$t\ge0$. For $x\notin
\supp_s\rho_{\mathrm{fc}}(t)$, these properties follow directly from the
definition of~$\widehat U'$ above. Thus~$\widehat U'(t,x)$
satisfies~(\ref{assumption1forbetauniversality}) and~(\ref
{assumption2forbetauniversality}) with uniform constants for
all~$t\ge0$ and~$N$ sufficiently large on~$\Omega$.

Finally, the potentials $\widehat U(t)$ and $U(t)$ are ``regular'' as
follows from Lemmas~\ref{lemmamfc} and~\ref{lemmahatmfc}.
\end{pf*}
\end{appendix}


\section*{Acknowledgments}
We thank Paul Bourgade, L\'aszl\'o Erd\H{o}s and Antti Know\-les for
helpful comments. We are grateful to the
Taida Institute for Mathematical Sciences and National Taiwan
Universality for their hospitality during part of this research. We
thank Thomas Spencer and the Institute for Advanced Study for their
hospitality during the academic year 2013--2014.


%

\printaddresses
\end{document}